
%
%
%
\def\unredoffs{} \def\redoffs{\voffset=-.31truein\hoffset=-.59truein}
\def\speclscape{\special{ps: landscape}}
%
%
%
%
\newbox\leftpage \newdimen\fullhsize \newdimen\hstitle \newdimen\hsbody
\tolerance=1000\hfuzz=2pt
\catcode`\@=11 
%
\ifx\answ\bigans\message{(This will come out unreduced.}
\magnification=1200\unredoffs\baselineskip=16pt plus 2pt minus 1pt
\hsbody=\hsize \hstitle=\hsize 
\else\message{(This will be reduced.} \let\l@r=L
\magnification=1000\baselineskip=16pt plus 2pt minus 1pt \vsize=7truein
\redoffs \hstitle=8truein\hsbody=4.75truein\fullhsize=10truein\hsize=\hsbody
\output={\ifnum\pageno=0 
  \shipout\vbox{\speclscape{\hsize\fullhsize\makeheadline}
    \hbox to \fullhsize{\hfill\pagebody\hfill}}\advancepageno
  \else
  \almostshipout{\leftline{\vbox{\pagebody\makefootline}}}\advancepageno
  \fi}
\def\almostshipout#1{\if L\l@r \count1=1 \message{[\the\count0.\the\count1]}
      \global\setbox\leftpage=#1 \global\let\l@r=R
 \else \count1=2
  \shipout\vbox{\speclscape{\hsize\fullhsize\makeheadline}
      \hbox to\fullhsize{\box\leftpage\hfil#1}}  \global\let\l@r=L\fi}
\fi
%
\newcount\yearltd\yearltd=\year\advance\yearltd by -1900

%

%
%

\def\draftmode{\message{ DRAFTMODE }\def\draftdate{{\rm preliminary draft:
\number\month/\number\day/\number\yearltd\ \ \hourmin}}%
\headline={\hfil\draftdate}\writelabels\baselineskip=20pt plus 2pt minus 2pt
 {\count255=\time\divide\count255 by 60 \xdef\hourmin{\number\count255}
  \multiply\count255 by-60\advance\count255 by\time
  \xdef\hourmin{\hourmin:\ifnum\count255<10 0\fi\the\count255}}}
\def\nolabels{\def\wrlabeL##1{}\def\eqlabeL##1{}\def\reflabeL##1{}}
\def\writelabels{\def\wrlabeL##1{\leavevmode\vadjust{\rlap{\smash%
{\line{{\escapechar=` \hfill\rlap{\sevenrm\hskip.03in\string##1}}}}}}}%
\def\eqlabeL##1{{\escapechar-1\rlap{\sevenrm\hskip.05in\string##1}}}%
\def\reflabeL##1{\noexpand\llap{\noexpand\sevenrm\string\string\string##1}}}
\nolabels
%
\global\newcount\secno \global\secno=0
\global\newcount\meqno \global\meqno=1
\def\newsec#1{\global\advance\secno by1\message{(\the\secno. #1)}
\global\subsecno=0\eqnres@t\noindent{\bf\the\secno. #1}
\writetoca{{\secsym} {#1}}\par\nobreak\medskip\nobreak}
\def\eqnres@t{\xdef\secsym{\the\secno.}\global\meqno=1\bigbreak\bigskip}
\def\sequentialequations{\def\eqnres@t{\bigbreak}}\xdef\secsym{}
\global\newcount\subsecno \global\subsecno=0
\def\subsec#1{\global\advance\subsecno by1\message{(\secsym\the\subsecno. #1)}
\ifnum\lastpenalty>9000\else\bigbreak\fi
\noindent{\it\secsym\the\subsecno. #1}\writetoca{\string\quad
{\secsym\the\subsecno.} {#1}}\par\nobreak\medskip\nobreak}
\def\appendix#1#2{\global\meqno=1\global\subsecno=0\xdef\secsym{\hbox{#1.}}
\bigbreak\bigskip\noindent{\bf Appendix #1. #2}\message{(#1. #2)}
\writetoca{Appendix {#1.} {#2}}\par\nobreak\medskip\nobreak}
%
%
\def\eqnn#1{\xdef #1{(\secsym\the\meqno)}\writedef{#1\leftbracket#1}%
\global\advance\meqno by1\wrlabeL#1}
\def\eqna#1{\xdef #1##1{\hbox{$(\secsym\the\meqno##1)$}}
\writedef{#1\numbersign1\leftbracket#1{\numbersign1}}%
\global\advance\meqno by1\wrlabeL{#1$\{\}$}}
\def\eqn#1#2{\xdef #1{(\secsym\the\meqno)}\writedef{#1\leftbracket#1}%
\global\advance\meqno by1$$#2\eqno#1\eqlabeL#1$$}
%
\newskip\footskip\footskip14pt plus 1pt minus 1pt 
\def\footnotefont{\ninepoint}\def\f@t#1{\footnotefont #1\@foot}
\def\f@@t{\baselineskip\footskip\bgroup\footnotefont\aftergroup\@foot\let\next}
\setbox\strutbox=\hbox{\vrule height9.5pt depth4.5pt width0pt}
\global\newcount\ftno \global\ftno=0
\def\foot{\global\advance\ftno by1\footnote{$^{\the\ftno}$}}
%
\newwrite\ftfile
\def\footend{\def\foot{\global\advance\ftno by1\chardef\wfile=\ftfile
$^{\the\ftno}$\ifnum\ftno=1\immediate\openout\ftfile=foots.tmp\fi%
\immediate\write\ftfile{\noexpand\smallskip%
\noexpand\item{f\the\ftno:\ }\pctsign}\findarg}%
\def\footatend{\vfill\eject\immediate\closeout\ftfile{\parindent=20pt
\centerline{\bf Footnotes}\nobreak\bigskip\input foots.tmp }}}
\def\footatend{}
%
%
\global\newcount\refno \global\refno=1
\newwrite\rfile
%
\def\ref{\nref}
\def\nref#1{\xdef#1{[\the\refno]}\writedef{#1\leftbracket#1}%
\ifnum\refno=1\immediate\openout\rfile=refs.tmp\fi
\global\advance\refno by1\chardef\wfile=\rfile\immediate
\write\rfile{\noexpand\item{#1\ }\reflabeL{#1\hskip.31in}\pctsign}\findarg}
\def\findarg#1#{\begingroup\obeylines\newlinechar=`\^^M\pass@rg}
{\obeylines\gdef\pass@rg#1{\writ@line\relax #1^^M\hbox{}^^M}%
\gdef\writ@line#1^^M{\expandafter\toks0\expandafter{\striprel@x #1}%
\edef\next{\the\toks0}\ifx\next\em@rk\let\next=\endgroup\else\ifx\next\empty%
\else\immediate\write\wfile{\the\toks0}\fi\let\next=\writ@line\fi\next\relax}}
\def\striprel@x#1{} \def\em@rk{\hbox{}}
\def\lref{\begingroup\obeylines\lr@f}
\def\lr@f#1#2{\gdef#1{\ref#1{#2}}\endgroup\unskip}

\def\addref#1{\immediate\write\rfile{\noexpand\item{}#1}} 
\def\startrefs#1{\immediate\openout\rfile=refs.tmp\refno=#1}
\def\refs#1{\count255=1[\r@fs #1{\hbox{}}]}
\def\r@fs#1{\ifx\und@fined#1\message{reflabel \string#1 is undefined.}%
\nref#1{need to supply reference \string#1.}\fi%
\vphantom{\hphantom{#1}}\edef\next{#1}\ifx\next\em@rk\def\next{}%
\else\ifx\next#1\ifodd\count255\relax\xref#1\count255=0\fi%
\else#1\count255=1\fi\let\next=\r@fs\fi\next}
%

%
\newwrite\ffile\global\newcount\figno \global\figno=1
\def\fig{fig.~\the\figno\nfig}
\def\nfig#1{\xdef#1{fig.~\the\figno}%
\writedef{#1\leftbracket fig.\noexpand~\the\figno}%
\ifnum\figno=1\immediate\openout\ffile=figs.tmp\fi\chardef\wfile=\ffile%
\immediate\write\ffile{\noexpand\medskip\noexpand\item{Fig.\ \the\figno. }
\reflabeL{#1\hskip.55in}\pctsign}\global\advance\figno by1\findarg}
\def\xfig{\expandafter\xf@g}\def\xf@g fig.\penalty\@M\ {}
\def\figs#1{figs.~\f@gs #1{\hbox{}}}
\def\f@gs#1{\edef\next{#1}\ifx\next\em@rk\def\next{}\else
\ifx\next#1\xfig #1\else#1\fi\let\next=\f@gs\fi\next}
\newwrite\lfile
{\escapechar-1\xdef\pctsign{\string\%}\xdef\leftbracket{\string\{}
\xdef\rightbracket{\string\}}\xdef\numbersign{\string\#}}

\def\writestop{\def\writestoppt{\immediate\write\lfile{\string\pageno%
\the\pageno\string\startrefs\leftbracket\the\refno\rightbracket%
\string\def\string\secsym\leftbracket\secsym\rightbracket%
\string\secno\the\secno\string\meqno\the\meqno}\immediate\closeout\lfile}}
\def\writestoppt{}\def\writedef#1{}
\def\seclab#1{\xdef #1{\the\secno}\writedef{#1\leftbracket#1}\wrlabeL{#1=#1}}
\def\subseclab#1{\xdef #1{\secsym\the\subsecno}%
\writedef{#1\leftbracket#1}\wrlabeL{#1=#1}}
\newwrite\tfile \def\writetoca#1{}
\def\leaderfill{\leaders\hbox to 1em{\hss.\hss}\hfill}
\def\writetoc{\immediate\openout\tfile=toc.tmp
   \def\writetoca##1{{\edef\next{\write\tfile{\noindent ##1
   \string\leaderfill {\noexpand\number\pageno} \par}}\next}}}
%
%
%

\catcode`\@=12 
%
\edef\tfontsize{\ifx\answ\bigans scaled\magstep3\else scaled\magstep4\fi}
\font\titlerm=cmr10 \tfontsize \font\titlerms=cmr7 \tfontsize
\font\titlermss=cmr5 \tfontsize \font\titlei=cmmi10 \tfontsize
\font\titleis=cmmi7 \tfontsize \font\titleiss=cmmi5 \tfontsize
\font\titlesy=cmsy10 \tfontsize \font\titlesys=cmsy7 \tfontsize
\font\titlesyss=cmsy5 \tfontsize \font\titleit=cmti10 \tfontsize
\skewchar\titlei='177 \skewchar\titleis='177 \skewchar\titleiss='177
\skewchar\titlesy='60 \skewchar\titlesys='60 \skewchar\titlesyss='60
\def\titlefont{\def\rm{\fam0\titlerm}
\textfont0=\titlerm \scriptfont0=\titlerms \scriptscriptfont0=\titlermss
\textfont1=\titlei \scriptfont1=\titleis \scriptscriptfont1=\titleiss
\textfont2=\titlesy \scriptfont2=\titlesys \scriptscriptfont2=\titlesyss
\textfont\itfam=\titleit \def\it{\fam\itfam\titleit}\rm}
 \ifx\answ\bigans\else scaled\magstep1\fi
\ifx\answ\bigans\else

 \font\absi=cmmi10 scaled\magstep1
\font\absis=cmmi7 scaled\magstep1 \font\absiss=cmmi5 scaled\magstep1
\font\abssy=cmsy10 scaled\magstep1 \font\abssys=cmsy7 scaled\magstep1
\font\abssyss=cmsy5 scaled\magstep1 
\skewchar\absi='177 \skewchar\absis='177 \skewchar\absiss='177
\skewchar\abssy='60 \skewchar\abssys='60 \skewchar\abssyss='60
\fi
\font\ninerm=cmr9 \font\sixrm=cmr6 \font\ninei=cmmi9 \font\sixi=cmmi6
\font\ninesy=cmsy9 \font\sixsy=cmsy6 \font\ninebf=cmbx9
\font\nineit=cmti9 \font\ninesl=cmsl9 \skewchar\ninei='177
\skewchar\sixi='177 \skewchar\ninesy='60 \skewchar\sixsy='60
\def\ninepoint{\def\rm{\fam0\ninerm}
\textfont0=\ninerm \scriptfont0=\sixrm \scriptscriptfont0=\fiverm
\textfont1=\ninei \scriptfont1=\sixi \scriptscriptfont1=\fivei
\textfont2=\ninesy \scriptfont2=\sixsy \scriptscriptfont2=\fivesy
\textfont\itfam=\ninei \def\it{\fam\itfam\nineit}\def\sl{\fam\slfam\ninesl}%
\textfont\bffam=\ninebf \def\bf{\fam\bffam\ninebf}\rm}
%
%

\hyphenation{anom-aly anom-alies coun-ter-term coun-ter-terms}
\def\inv{^{\raise.15ex\hbox{${\scriptscriptstyle -}$}\kern-.05em 1}}

\def\Dsl{\,\raise.15ex\hbox{/}\mkern-13.5mu D} 
\def\dsl{\raise.15ex\hbox{/}\kern-.57em\partial}

\def\lspace{\ifx\answ\bigans{}\else\qquad\fi}
\def\lbspace{\ifx\answ\bigans{}\else\hskip-.2in\fi} 
\def\boxeqn#1{\vcenter{\vbox{\hrule\hbox{\vrule\kern3pt\vbox{\kern3pt
    \hbox{${\displaystyle #1}$}\kern3pt}\kern3pt\vrule}\hrule}}}
\def\mbox#1#2{\vcenter{\hrule \hbox{\vrule height#2in
        \kern#1in \vrule} \hrule}}  
%

\def\darr#1{\raise1.5ex\hbox{$\leftrightarrow$}\mkern-16.5mu #1}

\def\half{{\textstyle{1\over2}}} 
\def\roughly#1{\raise.3ex\hbox{$#1$\kern-.75em\lower1ex\hbox{$\sim$}}}

%
%


\def\frac#1#2{{#1\over#2}}

\def\half{\frac12}

\def\journal#1&#2(#3){\unskip, #1~\bf #2 \rm(19#3) }
\def\andjournal#1&#2(#3){\sl #1~\bf #2 \rm (19#3) }

\def\bra#1{\left\langle #1\right|}
\def\ket#1{\left| #1\right\rangle}

\def\One{{1\hskip -3pt {\rm l}}}
\catcode`\@=11\def\slash#1{\mathord{\mathpalette\c@ncel{#1}}}
\overfullrule=0pt
\def\steepslash{\c@ncel}
\def\frac#1#2{{#1\over #2}}

\def\:{\!:\!}
\def\inbar{\,\vrule height1.5ex width.4pt depth0pt}
\def\IQ{\relax\,\hbox{$\inbar\kern-.3em{\rm Q}$}}
\def\IB{\relax{\rm I\kern-.18em B}}
\def\IC{\relax\hbox{$\inbar\kern-.3em{\rm C}$}}
\def\IP{\relax{\rm I\kern-.18em P}}
\def\IR{\relax{\rm I\kern-.18em R}}
\def\ZZ{\relax\ifmmode\mathchoice
{\hbox{Z\kern-.4em Z}}{\hbox{Z\kern-.4em Z}}
{\lower.9pt\hbox{Z\kern-.4em Z}}
{\lower1.2pt\hbox{Z\kern-.4em Z}}\else{Z\kern-.4em Z}\fi}
\def\ol{\overline}

\catcode`\@=12

\def\npb#1(#2)#3{{ Nucl. Phys. }{B#1} (#2) #3}
\def\plb#1(#2)#3{{ Phys. Lett. }{#1B} (#2) #3}
\def\pla#1(#2)#3{{ Phys. Lett. }{#1A} (#2) #3}
\def\prl#1(#2)#3{{ Phys. Rev. Lett. }{#1} (#2) #3}
\def\mpla#1(#2)#3{{ Mod. Phys. Lett. }{A#1} (#2) #3}
\def\ijmpa#1(#2)#3{{ Int. J. Mod. Phys. }{A#1} (#2) #3}
\def\cmp#1(#2)#3{{ Comm. Math. Phys. }{#1} (#2) #3}
\def\cqg#1(#2)#3{{ Class. Quantum Grav. }{#1} (#2) #3}
\def\jmp#1(#2)#3{{ J. Math. Phys. }{#1} (#2) #3}
\def\anp#1(#2)#3{{ Ann. Phys. }{#1} (#2) #3}
\def\prd#1(#2)#3{{ Phys. Rev. } {D{#1}} (#2) #3}
\def\ptp#1(#2)#3{{ Progr. Theor. Phys. }{#1} (#2) #3}
\def\aom#1(#2)#3{{ Ann. Math. }{#1} (#2) #3}

\def\noi{\noindent}
\def\bs{\bigskip}
\def\bu{\bs\noindent $\bullet$}

\def\br{\buildrel}
\def\bra{\langle}
\def\ket{\rangle}

\def\C{{\bf C}}

\def\F{{\bf F}}

\def\H{{\bf H}}

\def\L{{\bf L}}

\def\R{{\bf R}}

\def\Z{{\bf Z}}
\def\cA{{\cal A}}
\def\cB{{\cal B}}
\def\cC{{\cal C}}
\def\cD{{\cal D}}
\def\cE{{\cal E}}
\def\cF{{\cal F}}
\def\cG{{\cal G}}

\def\cI{{\cal I}}

\def\cQ{{\cal Q}}

\def\cS{{\cal S}}

\def\cV{{\cal V}}
\def\cW{{\cal W}}

\input amssym
\def\ga{{\goth a}}
\def\gb{{\goth b}}

\def\gg{{\goth g}}
\def\gh{{\goth h}}

\def\gj{{\goth j}}

\def\gs{{\goth s}}
\def\gt{{\goth t}}

\def\gU{{\goth U}}

\def\gX{{\goth X}}

\def\da#1{{\p \over \p s_{#1}}}

\def\cicy#1(#2|#3)#4{\left(\matrix{#2}\right|\!\!
                     \left|\matrix{#3}\right)^{{#4}}_{#1}}

\def\hra{\hookrightarrow}

\def\ra{\rightarrow}

\def\da{\downarrow}

\def\bs{\bigskip}

\def\Box{{\,\lower0.9pt\vbox{\hrule
\hbox{\vrule height 0.2 cm \hskip 0.2 cm
\vrule height 0.2 cm}\hrule}\,}}

\global\newcount\thmno \global\thmno=0
\def\definition#1{\global\advance\thmno by1
\bigskip\noindent{\bf Definition \secsym\the\thmno. }{\it #1}
\par\nobreak\medskip\nobreak}
\def\question#1{\global\advance\thmno by1
\bigskip\noindent{\bf Question \secsym\the\thmno. }{\it #1}
\par\nobreak\medskip\nobreak}
\def\theorem#1{\global\advance\thmno by1
\bigskip\noindent{\bf Theorem \secsym\the\thmno. }{\it #1}
\par\nobreak\medskip\nobreak}
\def\proposition#1{\global\advance\thmno by1
\bigskip\noindent{\bf Proposition \secsym\the\thmno. }{\it #1}
\par\nobreak\medskip\nobreak}
\def\corollary#1{\global\advance\thmno by1
\bigskip\noindent{\bf Corollary \secsym\the\thmno. }{\it #1}
\par\nobreak\medskip\nobreak}
\def\lemma#1{\global\advance\thmno by1
\bigskip\noindent{\bf Lemma \secsym\the\thmno. }{\it #1}
\par\nobreak\medskip\nobreak}
\def\conjecture#1{\global\advance\thmno by1
\bigskip\noindent{\bf Conjecture \secsym\the\thmno. }{\it #1}
\par\nobreak\medskip\nobreak}
\def\exercise#1{\global\advance\thmno by1
\bigskip\noindent{\bf Exercise \secsym\the\thmno. }{\it #1}
\par\nobreak\medskip\nobreak}
\def\remark#1{\global\advance\thmno by1
\bigskip\noindent{\bf Remark \secsym\the\thmno. }{\it #1}
\par\nobreak\medskip\nobreak}
\def\problem#1{\global\advance\thmno by1
\bigskip\noindent{\bf Problem \secsym\the\thmno. }{\it #1}
\par\nobreak\medskip\nobreak}
\def\others#1#2{\global\advance\thmno by1
\bigskip\noindent{\bf #1 \secsym\the\thmno. }{\it #2}
\par\nobreak\medskip\nobreak}
\def\proof{\noindent Proof: }

\def\thmlab#1{\xdef #1{\secsym\the\thmno}\writedef{#1\leftbracket#1}\wrlabeL{#1=#1}}
%
%
\def\newsec#1{\global\advance\secno by1\message{(\the\secno. #1)}
\global\subsecno=0\thmno=0\eqnres@t\noindent{\bf\the\secno. #1}
\writetoca{{\secsym} {#1}}\par\nobreak\medskip\nobreak}
\def\eqnres@t{\xdef\secsym{\the\secno.}\global\meqno=1\bigbreak\bigskip}
\def\sequentialequations{\def\eqnres@t{\bigbreak}}\xdef\secsym{}
%

%
\newcount{\exnum}
\def\prob{\advance\exnum by 1
\bigskip\item{\the\exnum.}\ }
\newcount{\exnum}
\def\next{\advance\exnum by 1
\bigskip\noindent{\the\exnum.}\ }
\def\np{\vfill\eject}

\nopagenumbers

\ref\AK{ F. Akman,  The Semi-infinite Weil Complex of a Graded Lie
Algebra,  PhD Thesis, Yale University, 1993.} \ref\Berezin{F. A.
Berezin, Introduction to Superanalysis, Math. Phys. and Appl.
Math., 9 D. Reidel Publishing Co., Dordrecht-Boston (1987).}
\ref\Bor{R. Borcherds, Vertex operator algebras, Kac-Moody
algebras and the Monster, Proc. Nat. Acad. Sci. USA 83 (1986)
3068-3071.} 
\ref\DVV{R. Dijkgraaf, H. Verlinde and E. Verlinde,
Notes on topological string theory and 2d quantum gravity, in
String Theory and Quantum Gravity (Trieste 1990), World Sci.
Publishing, River Edge, NJ (1991) 91-156.} 
\ref\CartanI{ H.
Cartan, Notions d'alg\`ebre diff\'erentielle; application aux
groupes de Lie et aux vari\'eti\'es o\`u op\`ere un groupe de Lie,
Colloque de Topologie, C.B.R.M., Bruxelles 15-27 (1950).}
\ref\CartanII{ H. Cartan, La Transgression dans un groupe de Lie
et dans un espace fibr\'e principal, Colloque de Topologie,
C.B.R.M., Bruxelles 57-71 (1950).} \ref\EY{T. Eguchi and S.K.
Yang, N=2 superconformal models as topological field theories,
Mod. Phys. Lett. A5 (1990) 1693-1701. }
\ref\DKV{M. Duflo, S. Kumar, and M. Vergne, Sur la Cohomologie \'Equivariante des Vari\'et\'es Diff\'erentiables, Ast\'erisque 215 (1993).}
 \ref\F{ B. Feigin,   The Semi-Infinite Homology of Lie, Kac-Moody
and Virasoro Algebras, Russian Math. Survey 39 (1984) no. 2
195-196.}
 \ref\FF{ B. Feigin and E. Frenkel,   Semi-Infinite Weil Complex
and the Virasoro Algebra,  Comm. Math. Phys. 137 (1991), 617-639.}
 \ref\FGZ{ I.B. Frenkel, H. Garland, G.J. Zuckerman,   Semi-Infinite
Cohomology and String Theory , Proc. Natl. Acad. Sci. USA Vol. 83,
No. 22 (1986) 8442-8446.} 
\ref\FHL{I.B. Frenkel, Y.Z. Huang, J.
Lepowsky, On axiomatic approaches to vertex operator algebras and
modules, Mem. Amer. Math. Soc. 104 (1993), no. 494, viii+64.}
 \ref\FLM{ I.B. Frenkel, J. Lepowsky, A. Meurman,   Vertex
Operator Algebras and the Monster, Academic Press, New York,
1988.}
 \ref\FZ{ I.B. Frenkel, and Y.C. Zhu, Vertex Operator Algebras
Associated to Representations of Affine and Virasoro Algebras ,
Duke Mathematical Journal, Vol. 66, No. 1, (1992), 123-168.}
 \ref\FMS{ D. Friedan, E. Martinec, S. Shenker,   Conformal
Invariance, Supersymmetry and String Theory , Nucl. Phys. B271
(1986) 93-165.} \ref\Getz{E. Getzler, Two-dimensional topological
gravity and equivariant cohomology, Commun. Math. Phys. 163 (1994)
473-489.} 
\ref\GKO{ P. Goddard, A. Kent, and D. Olive,   Virasoro
Algebras and Coset Space Models,  Phys. Lett B 152 (1985) 88-93.}
 \ref\GS{ V. Guillemin and S. Sternberg, Supersymmetry and
Equivariant de Rham Theory, Springer, 1999.}
 \ref\HS{P. Hilton, and U. Stammbach, A Course in Homological
Algebra, 2nd ed., Springer-Verlag, New York 1997.}
 \ref\Kac{V. Kac, Vertex Algebras for Beginners, AMS Univ. Lecture
Series Vol. 10, 2nd corrected ed., 2001.}
\ref\KacPeterson{V. Kac, and D. Peterson, Infinite-dimensional Lie algebras, theta functions and modular forms, Adv. Math. 53 (1984) 125-264.}
\ref\Li{H. Li, Local
systems of vertex operators, vertex superalgebras and modules,
hep-th/9406185.}
\ref\Kalkman{J. Kalkman,  A BRST Model Applied to
Symplectic Geometry, Thesis, Utrecht (1993).}
 \ref\Kostant{ B. Kostant, Graded manifolds, graded
Lie theory and pre-quantization, Lecture Notes in Math. 570,
Springer-Verlag, Berlin (1977) 177-306.}
 \ref\LZI{ B. Lian, G.J. Zuckerman, New perspectives on the BRST-algebraic structure of string theory,
Commun. Math. Phys. 154 (1993) 613-646.}
 \ref\LZII{ B. Lian, G.J. Zuckerman,   Commutative Quantum Operator
Algebras , J. Pure Appl. Algebra 100 (1995) no. 1-3, 117-139.}
\ref\LtypeI{B. Lian, On the classification of simple vertex
operator algebras, Comm. Math. Phys. 163 (1994) No. 2, 307-357.}
 \ref\Linshaw{A. Linshaw, Vertex Algebras  and Invariant Theory, Ph.D. Thesis, Brandeis University, 2005.}
 \ref\MSV{F. Malikov, V. Schectman, and A. Vaintrob, Chiral de Rham Complex, Comm. Math. Phys, 204, (1999) 439-473.}
 \ref\MalikovSchectman{F. Malikov, V. Schectman, Chiral de Rham Complex II, arXiv:math/AG9901065 v4.}
 \ref\MQ{V. Mathai, D. Quillen,  Superconnections, Thom classes and equivariant differential forms, Topology 25 (1986), no. 1, 85-110.}
  \ref\MS{G. Moore, and N. Seiberg, Classical and quantum conformal field theory, Comm. Math. Phys. 123 (1989) 177-254.}
  \ref\Pol{J. Polchinski, String Theory, Vol 2, Cambridge Press 1998, Chap. 11.}
 \ref\Procesi{ C. Procesi,  The Invariant Theory of $n\times
n$ Matrices, Adv. Math. 19, (1976) 306-381.}
\ref\Quillen{D. Quillen, Superconnections and the Chern character, Topology 24,
(1985), no. 1, 89-95.}

\centerline{\titlefont Chiral Equivariant Cohomology I}

\bs
\centerline{Bong H. Lian and Andrew R. Linshaw}
\bs

ABSTRACT. We construct a new equivariant cohomology theory for a certain class of differential vertex algebras, which we call the chiral equivariant cohomology. A principal example of a differential vertex algebra in this class is the chiral de Rham complex of Malikov-Schechtman-Vaintrob of a manifold with a group action. The main idea in this paper is to synthesize the
algebraic approach to classical equivariant cohomology
due to H. Cartan$^1$\footnote{}{Cartan's theory was further developed by Duflo-Kumar-Vergne \DKV~and Guillemin-Sternberg \GS. This paper follows closely the latter approach.},
with the theory of differential vertex algebras, by using an appropriate notion of invariant theory. 
We also construct the vertex algebra analogues of the Mathai-Quillen isomorphism, the Weil and the
Cartan models for equivariant cohomology, and the Chern-Weil map. We give interesting cohomology classes in the new theory that have no classical analogues.

\bs
\noi Keywords: differential vertex algebras, equivariant de Rham theory, invariant theory, semi-infinite Weil algebra, Virasoro algebra.



\np
\headline{\ifodd\pageno\rightheadline\else\leftheadline\fi}
\def\rightheadline{\tenrm\hfil Chiral Equivariant Cohomology I
\hfil\folio}
\def\leftheadline{\tenrm\folio\hfil Bong H. Lian \& Andrew R. Linshaw \hfil}

\np
\centerline{\bf Contents}\nobreak\medskip{\baselineskip=12pt
 \parskip=0pt\catcode`\@=11 \input toc.tmp \catcode`\@=12 \bigbreak\bigskip}  

\newsec{Introduction}


\subsec{Equivariant de Rham theory}

For a topological space $M$ equipped with an action of a compact Lie group $G$,
the $G$-equivariant cohomology of $M$, denoted by $H^*_G(M)$, is defined to be $H^*((M\times \cE)/G))$, where $\cE$ is any contractible topological space on which $G$ acts freely. When $M$ is a smooth manifold on which $G$ acts by diffeomorphisms, there is a de Rham model of $H^*_G(M)$ due to H. Cartan \CartanI\CartanII, and developed further by Duflo-Kumar-Vergne \DKV~and Guillemin-Sternberg \GS. The treatment in \GS~
is simplified considerably by the use of supersymmetry
\Berezin\Kostant\MQ\Quillen\Kalkman, and will be the approach adopted in the present paper.
Guillemin-Sternberg define the equivariant cohomology $H^*_G(A)$ of any $G^*$-algebra
$A$, of which the algebra $\Omega(M)$ of smooth differential forms on $M$ is an example. A $G^*$-algebra is a
commutative superalgebra $A$ equipped with an action of $G$,
together with a compatible action of a certain differential Lie
superalgebra $(\gs\gg,d)$ associated to the Lie algebra $\gg$ of
$G$. Taking $A=\Omega(M)$ gives us the de Rham model of $H^*_G(M)$, and $H^*_G(\Omega(M)) = H^*_G(M)$ by an equivariant version of the de Rham theorem.  A $G^*$-algebra $(A,d)$ is a cochain complex, and the subalgebra of $A$ which is both $G$-invariant and killed by $\gs\gg$ forms a subcomplex
known as the basic subcomplex. $H^*_G(A)$ may be defined to be
$H^*_{bas}(A\otimes W(\gg))$, where $W(\gg) =\Lambda(\gg^*)\otimes
S(\gg^*)$ is the Koszul complex of $\gg$. A change of variables
\GS~shows that $W(\gg)$ is isomorphic to the subcomplex of
$\Omega(EG)$ which is freely generated by the connection one-forms
and curvature two-forms. Here $EG$ is the total space of the
classifying bundle of $G$ and $\Omega(EG)$ is the de Rham complex
of $EG$. This subcomplex is known as the Weil complex, and
$H^*_{bas}(A\otimes W(\gg))$ is known as the {\it Weil model} for
$H^*_G(A)$. Using an automorphism of the space $A\otimes W(\gg)$
called the {\it Mathai-Quillen isomorphism}, one can construct the
{\it Cartan model} which is often more convenient for computational
purposes.

\subsec{A vertex algebra analogue of $G^*$-algebras}

Associated to $G$ is a certain universal differential vertex algebra we call $O(\gs\gg)$, which is analogous to $(\gs\gg,d)$.
An $O(\gs\gg)$-algebra is then a differential vertex algebra $\cA$ equipped with an action of $G$ together with a compatible action of $O(\gs\gg)$. When $G$ is connected, the $G$-action can be absorbed into the $O(\gs\gg)$-action.

Associated to any Lie algebra $\gg$ is a $\Z_{\geq0}$-graded vertex algebra $\cW(\gg)$
known as the {\it semi-infinite Weil complex} of $\gg$ (a.k.a. the $bc\beta\gamma$-system in physics \FMS).
When $\gg$ is finite dimensional, $\cW(\gg)$ is a vertex algebra which contains the classical Weil complex $W(\gg)$
as the subspace of conformal weight zero, and is an example of an $O(\gs\gg)$-algebra. It was studied in \FF\AK~ in the context of semi-infinite cohomology of the loop algebra of $\gg$,
which is a vertex algebra analogue of Lie algebra cohomology and an example of the general theory of semi-infinite cohomology developed in
\F\FGZ. 

In \MSV, Malikov-Schechtman-Vaintrob constructed a sheaf $\cQ_M$
of vertex algebras on any nonsingular algebraic scheme $M$, which
they call the {\it chiral de Rham sheaf}. They also pointed out that the same construction can be done in the analytic and smooth categories (see remark 3.9 \MSV.) In this paper, we will carry out a
construction that is equivalent in the smooth category. 
The space $\cQ(M)$ of global sections of the MSV sheaf
$\cQ_M$ is a $\Z_{\geq 0}$-graded vertex algebra, graded by
conformal weight, which contains the ordinary de Rham algebra
$\Omega(M)$ as the subspace of conformal weight zero. There is a
square-zero derivation $d_\cQ$ on $\cQ(M)$ whose restriction to
$\Omega(M)$ is the ordinary de Rham differential $d$, and the
inclusion of complexes $(\Omega(M),d)\hookrightarrow
(\cQ(M),d_\cQ)$ induces an isomorphism in cohomology. When $M$ is a $G$-manifold, the algebra $\cQ(M)$ is another example of an $O(\gs\gg)$-algebra.

\subsec{Vertex algebra invariant theory}

For any vertex algebra $V$ and any subalgebra $B\subset V$, there is a new subalgebra
$Com(B,V)\subset V$ known as the {\it commutant} of $B$ in $V$.
This construction was introduced in \FZ~ as a vertex algebra abstraction
of a construction in representation theory \KacPeterson~ and in conformal field theory \GKO, known as the coset construction. It may be interpreted either as the vertex algebra analogue of the ordinary
commutant construction in the theory of associative algebras, or
as a vertex algebra notion of invariant theory. The latter
interpretation was developed in \Linshaw, and is the point of view
we adopt in this paper.

\subsec{Chiral equivariant cohomology of an $O(\gs\gg)$-algebra}

Our construction of chiral equivariant cohomology synthesizes the
three theories outlined above. 
We define the chiral equivariant cohomology $\H^*_G(\cA)$ of any $O(\gs\gg)$-algebra
$\cA$ by replacing the main ingredients in the classical Weil
model for equivariant cohomology with their vertex algebra counterparts.
The commutant construction plays the same role that ordinary invariant
theory plays in classical equivariant cohomology. We also construct the chiral analogues of the
Mathai-Quillen isomorphism, the Cartan model for $\H^*_G(\cA)$, and
a vertex algebra homomorphism $\kappa_G:\H^*_G(\C)\rightarrow
\H^*_G(\cA)$ which is the chiral version of the Chern-Weil map. Here $\C$ is the one-dimensional trivial $O(\gs\gg)$-algebra.

Specializing to $\cA=\cQ(M),$ for a $G$-manifold $M$, gives us a chiral
equivariant cohomology theory of $M$ which contains the classical equivariant cohomology. It turns out that there are other interesting differential vertex algebras $\cA$, some of which are subalgebras of $\cQ(M)$, for which $\H_G^*(\cA)$ can be defined and also contains the classical equivariant cohomology $M$. This will be the focus of a separate paper; see further remarks in the last section.

In the case where $G$ is an $n$-dimensional torus $T$, we give a
complete description of $\H^*_T(\C)$. Working in the Cartan model,
we show that for any $O(\gs\gt)$-algebra $\cA$, $\H^*_T(\cA)$ is actually the
cohomology of a much smaller subcomplex of the chiral Cartan
complex, which we call the small chiral Cartan complex. Like the
classical Cartan complex, the small chiral Cartan complex has the
structure of a double complex, and there is an associated
filtration and spectral sequence that computes $\H_T^*(\cA)$.
For non-abelian $G$, we also construct a double complex structure in the Weil and Cartan models, and derive two corresponding spectral sequences.

When $G$ is a simple, connected Lie group, we show that
$\H^*_G(\C)$ contains a vertex operator $\L(z)$ which has no
classical analogue, and satisfies the Virasoro OPE relation. In
particular, $\H^*_G(\C)$ is an interesting non-abelian vertex
algebra. This algebra plays the role of $H_G^*(\C)=S(\gg^*)^G$, the equivariant cohomology of a point, in the classical theory.

{\it Acknowledgement.} We thank J. Levine for discussions and for
his interest in this work, and G. Schwarz for helpful discussions on invariant theory. We thank E. Paksoy and B. Song for helping to correct a number of mistakes in an earlier draft of this paper.
B.H.L.'s research is partially supported by a J.S. Guggenheim Fellowship and an NUS grant.
A.R.L. would like to thank the Department of Mathematics, The National University of Singapore, for its hospitality and financial support during his visit there, where this paper was written.

\newsec{Background}

In this section we discuss the necessary background material in
preparation for the main results to be developed in sections 2-7.
Vertex algebras and modules have been discussed from various
different points of view in \FMS\MS\Bor\FLM\FHL\LZII\Li\Kac. We
will follow the formalism developed in \LZII~and partly in \Li. We
also carry out the construction of the chiral de Rham sheaf in the
smooth category.

\subsec{An interlude on vertex algebras}

Let $V$ be a vector space (always assumed defined over the complex numbers).
Let $z,w$ be formal variables.
By $QO(V)$, we mean the space of all linear maps $V\ra V((z)):=\{\sum_{n\in\Z} v(n) z^{-n-1}|
v(n)\in V,~v(n)=0~for~n>>0\}$. Each element $a\in QO(V)$ can be uniquely represented
as a power series $a=a(z):=\sum_{n\in\Z}a(n)z^{-n-1}\in (End~V)[[z,z^{-1}]]$, though the latter space is clearly much larger than $QO(V)$. We refer to $a(n)$ as the $n$-th Fourier mode of $a(z)$.
If one regards $V((z))$ as a kind of ``$z$-adic'' completion of $V[z,z^{-1}]$, then $a\in QO(V)$ can be thought of as a map on $V((z))$ which is only defined on the dense subset $V[z,z^{-1}]$. {\it When $V$ is equipped with a super vector space structure $V=V^0\oplus V^1$
then an element $a\in QO(V)$ is assumed to be of the shape $a=a_0+a_1$
where $a_i:V^j\ra V^{i+j}((z))$ for $i,j\in\Z/2$.} 

On $QO(V)$ there is a set of non-associative bilinear operations, $\circ_n$,
indexed by $n\in\Z$, which we call the $n$-th circle products. They are defined by
$$
a(w)\circ_n b(w)=Res_z a(z)b(w)~i_{|z|>|w|}(z-w)^n-
Res_z b(w)a(z)~i_{|w|>|z|}(z-w)^n\in QO(V).
$$
Here $i_{|z|>|w|}f(z,w)\in\C[[z,z^{-1},w,w^{-1}]]$ denotes the power series expansion
of a rational function $f$ in the region $|z|>|w|$. Be warned that $i_{|z|>|w|}(z-w)^{-1}\neq
i_{|w|>|z|}(z-w)^{-1}$. As it is customary, we shall drop the symbol $i_{|z|>|w|}$
and just write $(z-w)^{-1}$ to mean the expansion in the region $|z|>|w|$,
and write $-(w-z)^{-1}$ to mean the expansion in $|w|>|z|$.
$Res_z(\cdots)$ here means taking the coefficient of $z^{-1}$ of $(\cdots)$.
It is easy to check that $a(w)\circ_n b(w)$ above is a well-defined element of $QO(V)$.
{\it When $V$ is equipped with a super vector space structure then
the definition of $a\circ_n b$ above is replaced by one with the extra
sign $(-1)^{|a||b|}$ in the second term}. Here $|a|$ is the $\Z/2$ grading of a homogeneous element $a\in QO(V)$.

The circle products are connected through the {\it operator product expansion} (OPE) formula (\LZII, Prop. 2.3): for $a,b\in QO(V)$, we have
\eqn\dumb{
a(z)b(w)=\sum_{n\geq0}a(w)\circ_n b(w)~(z-w)^{-n-1}+:a(z)b(w):
}
where
$$\eqalign{
:a(z)b(w):&=a(z)_-b(w)+(-1)^{|a||b|} b(w)a(z)_+\cr
a(z)_-&=\sum_{n<0}a(n)z^{-n-1},~~~
a(z)_+=\sum_{n\geq0}a(n)z^{-n-1}. }$$ Note that $:a(w)b(w):$ is a
well-defined element of $QO(V)$. It is called the Wick product of
$a$ and $b$, and it coincides with $a\circ_{-1}b$. The other
negative circle products are related to this by \eqn\deriv{
n!~a(w)\circ_{-n-1}b(w)=:({d^n\over dw^n} a(w))~b(w):. } For
$a_1(z),...,a_k(z)\in QO(V)$, it is convenient to define the
$k$-fold iterated Wick product
$$
:a_1(z)a_2(z)\cdots a_k(z):{\br def\over=}:a_1(z)b(z):
$$
where $b(z)=:a_2(z)\cdots a_k(z):$. It is customary to rewrite \dumb~as
$$
a(z)b(w)\sim\sum_{n\geq0}a(w)\circ_n b(w)~(z-w)^{-n-1}.
$$
Thus $\sim$ means equal modulo the term $:a(z)b(w):$. Note that
when $a\circ_n b=0$ for $n>>0$ (which will be the case throughout
this paper later), then formally $a(z)b(w)$ can be thought of as a
kind of meromorphic function with poles along $z=w$. The product $a\circ_n b$
is formally $\oint_C a(z)b(w)(z-w)^n dz$ where $C$ is a small
circle around $w$ (hence the name circle product).

From the definition, we see  that
$$
a(w)\circ_0 b(w)=[a(0),b(w)].
$$
From this, it follows easily that $a\circ_0$ is a (graded) derivation of every circle product \LZI. This property of the zeroth circle product will be used often later.

The set $QO(V)$ is a nonassociative algebra with the operations $\circ_n$
and a unit $1$. We have $1\circ_n a=\delta_{n,-1}a$ for all $n$,
and $a\circ_n 1=\delta_{n,-1}a$ for $n\geq-1$.
We are interested in subalgebras $A\subset QO(V)$, i.e.
linear subspaces of $QO(V)$ containing 1, which are closed under the circle products.
In particular $A$ is closed under formal differentiation 
$$
\partial a(w)={d\over dw} a(w)=a\circ_{-2}1.
$$
We shall call such a subalgebra a {\it circle algebra} (also called a quantum operator algebra in \LZII).

\remark{Fix a nonzero vector $\One\in V$ and let $a,b\in QO(V)$ such that $a(z)_+\One=b(z)_+\One=0$ for $n\geq0$. Then it follows immediately from the definition of the circle products that $(a\circ_p b)_+(z)\One=0$ for all $p$. Thus if a circle algebra $A$ is generated by elements $a$ with the property that $a(z)_+\One=0$, then every element in $A$ has this property. In this case the vector $\One$ determines a linear map
$$
\chi:A\ra V,~~~a\mapsto a(-1)\One=\lim_{z\ra0}a(z)\One
$$
(called the creation map in \LZII)
\thmlab\Based
having the following basic properties:
\eqn\chieqn{
\chi(1)=\One,~~~\chi(a\circ_n b)=a(n)b(-1)\One,~~~\chi(\partial^p a)=p! ~a(-p-1)\One.
}
}\thmlab\RemarkCreation

\definition{We say that $a,b\in QO(V)$ circle commute if $(z-w)^N [a(z),b(w)]=0$ for some $N\geq0$. If $N$ can be chosen to be 0, then we say that $a,b$ commute. We say that $a\in QO(V)$ is a vertex operator if it circle commutes with itself.}

\definition{A circle algebra is said to be commutative if its elements pairwise circle commute.}
\thmlab\VA

Again when there is a $\Z/2$ graded structure, the bracket in the definition above means the super commutator. We will see shortly that the notion of a commutative circle algebra is essentially equivalent to the notion of a vertex algebra (see for e.g. \FLM). An easy calculation gives the following very useful characterization of circle commutativity.

\lemma{Given $N\geq0$ and $a,b\in QO(V)$, we have
$$\eqalign{
& ~~~~~(z-w)^N[a(z),b(w)]=0 \cr &\Longleftrightarrow
[a(z)_+,b(w)]=\sum_{p=0}^{N-1}(a\circ_p b)(w)~(z-w)^{-p-1} \cr
&~~\&~~[a(z)_-,b(w)]=\sum_{p=0}^{N-1}(-1)^p(a\circ_p b)(w)~(w-z)^{-p-1}\cr
&\Longleftrightarrow
[a(m),b(n)]=\sum_{p=0}^{N-1}\left(\matrix{m\cr p}\right)
(a\circ_p b)(m+n-p)~~~\forall m,n\in\Z.
}$$}
\thmlab\Commutator

Using this lemma, it is not difficult to show that for any circle
commuting $a(z),b(z)\in QO(V)$ and $n\in\Z$, we have \eqn\VOcom{
a(z)\circ_n b(z) = \sum_{p\in\Z} (-1)^{p+1}(b(z)\circ_p
a(z))\circ_{n-p-1}1. } Note that this is a finite sum by circle
commutativity and the fact that $c(z)\circ_k 1 =0$ for all
$c(z)\in QO(V)$ and $k\geq 0$.

Many known commutative circle algebras can be constructed as
follows. Start with a set $S\subset QO(V)$ and use this lemma to
verify circle commutativity of the set. Then $S$ generates a
commutative circle algebra $A$ by the next lemma \LZII\Li.

\lemma{Let $a,b,c\in QO(V)$ be such that any two of them circle
commute. Then $a$ circle commutes with all $b\circ_p c$.} \proof
We have \eqn\dumb{
[a(z_1),[b(z_2),c(z_3)]]=[[a(z_1),b(z_2)],c(z_3)]]\pm[b(z_2),[a(z_1),c(z_3)]].
} For $M,N\geq0$, write
$(z_1-z_3)^{M+N}=((z_1-z_2)+(z_2-z_3))^N(z_1-z_3)^M$ and expand
the first factor binomially. For $M,N>>0$, each term
$(z_1-z_2)^i(z_2-z_3)^{N-i}(z_1-z_3)^M$ annihilates either the
left side of $\dumb$ or the right side. Thus $(z_1-z_3)^{M+N}$
annihilates \dumb. Multiplying \dumb~by $(z_2-z_3)^p$, $p\geq0$,
and taking $Res_{z_2}$, we see that
$(z_1-z_3)^{M+N}[a(z_1),(b\circ_pc)(z_3)]=0$. From this, we can
also conclude that $c$ circle commutes with all $b\circ_p a$,
$p\geq0$.

Now consider the case $p<0$.  For simplicity, we write $a=a(z),b=b(w),c=c(w)$. Suppose $(z-w)^N[a,b]=0$. Differentiating $(z-w)^{N+1}[a,b]=0$ with respect to $w$ shows that
$a$ circle commutes with $\partial b$. By \deriv, it remains to show that $a$ circle commutes with the Wick product $:bc:$. We have
$$
[a,:bc:]=[a,b_-]c\pm b_-[a,c]+[a,c]b_+\pm c[a,b_+].
$$
For $M>>0$, $(w-z)^M$ annihilates $[a,b],[a,c]$. In particular, $(w-z)^M[a,b_-]=(w-z)^M[b_+,a]$.
It follows that
$$
(w-z)^M[a,:bc:]=(w-z)^M[b_+,a]c\mp c(w-z)^M[b_+,a].
$$
For $M>>0$, the right side is zero by Lemma \Commutator~because $c$ circle commutes with
all $b\circ_pa$, $p\geq0$, and that $b\circ_pa=0$ for $p>>0$. $\Box$

In the formulation Definition \VA, many formal algebraic notions become immediately clear: a homomorphism is just a linear map that preserves all circle products and 1; a module over a circle algebra $A$ is a vector space $M$ equipped with a circle algebra homomorphism $A\ra QO(M)$, etc. For example, every commutative circle algebra $A$ is itself a faithful $A$-module, called {\it the left regular module}, as we now show. Define
$$
\rho:A\ra QO(A),~~a\mapsto\hat a,~~~~~\hat a(\zeta)b=\sum (a\circ_n b)~\zeta^{-n-1}.
$$

\lemma{For $a,b\in A$, $m,n\in\Z$, we have
$$
[\hat a(m),\hat b(n)]=\sum_{p\geq0}\left(\matrix{m\cr p}\right)\widehat{a\circ_p b}(m+n-p).
$$}
\proof
Applying the left side to a test vector $u\in A$, and using Lemma \Commutator, we have
$$\eqalign{
&\hat a(m)\cdot \hat b(n)\cdot u(z)-\hat b(n)\cdot \hat a(m)\cdot u(z)\cr
&=Res_{z_1}Res_{z_2}[a(z_2),b(z_1)]u(z)(z_2-z)^m(z_1-z)^n\cr
&-Res_{z_2}Res_{z_1}u(z)[a(z_2),b(z_1)](-z+z_2)^m(-z+z_1)^n\cr
&=Res_{z_1}Res_{z_2}\sum_{p\geq0}(a\circ_pb)(z_1)u(z){(-1)^p\over p!}
\left(\partial_{z_2}^p\delta(z_2,z_1)\right)
(z_2-z)^m(z_1-z)^n\cr
&-Res_{z_1}Res_{z_2}\sum_{p\geq0}u(z)(a\circ_pb)(z_1){(-1)^p\over p!}
\left(\partial_{z_2}^p\delta(z_2,z_1)\right)
(-z+z_2)^m(-z+z_1)^n
}$$
where $\delta(z_1,z_2)=(z_1-z_2)^{-1}+(z_2-z_1)^{-1}$.
By doing formal integration by parts and using the fact that $Res_{z_2}z_2^n\delta(z_2,z_1)=z_1^n$, 
the last expression becomes
$$
\sum_pRes_{z_1}\left(\matrix{m\cr p}\right)(a\circ_p b)(z_1)u(z)
(z_1-z)^{m-p+n}
-\sum_pRes_{z_1}\left(\matrix{m\cr p}\right)u(z)(a\circ_p b)(z_1)
(-z+z_1)^{m-p+n}.
$$
This is equal to the right side of our assertion applied on $u$. $\Box$

\theorem{$\rho$ is an injective circle algebra homomorphism.}
\proof
We will consider the case without the $\Z/2$ grading. The argument carries over to superalgebra case with some sign changes, as usual. The map $\rho$ is injective because $\hat a(-1)1=a\circ_{-1}1=a$.
Multiplying the formula in the preceding lemma by $\zeta^{-n-1}$ and summing over $n$, we find
\eqn\dumb{
[\hat a(m),\hat b(\zeta)]=\sum_{p\geq0}\left(\matrix{m\cr p}\right)
\widehat{a\circ_p b}(\zeta)~\zeta^{m-p}.
}
On the other hand, it follows from the OPE formula that for $m\geq0$,
$$
[\hat a(m),\hat b(\zeta)]=\sum_{p\geq0}\left(\matrix{m\cr p}\right)
(\hat a\circ_p\hat b)(\zeta)~\zeta^{m-p}.
$$
Specializing the two preceding formulas to $m=0,1,2,...$, we find that
$$
\widehat{a\circ_p b}=\hat a\circ_p\hat b
$$
for $p\geq0$. This shows that $\rho$ preserves the circle products $\circ_p$, $p\geq0$.
In particular \dumb~becomes
$$
[\hat a(m),\hat b(\zeta)]=\sum_{p\geq0}\left(\matrix{m\cr p}\right)
(\hat a\circ_p\hat b)(\zeta)~\zeta^{m-p}
$$
for all $m\in\Z$. This implies that $\hat a,\hat b$, circle commute, by Lemma \Commutator.

Let $A'$ be the (commutative) circle algebra generated by $\rho(A)$ in $QO(A)$. Since $\hat a(n)1=a\circ_n1=0$ for $a\in A$, $n\geq0$, i.e. $\hat a_+1=0$, it follows that every element $\alpha\in A'$ has $\alpha_+1=0$ by Remark \Based. Consider the creation map $\chi:A'\ra A$, $\alpha\mapsto\alpha(-1)1$, which is clearly surjective because $\chi\circ\rho=id$. We also have $[\partial,\hat a(\zeta)]b={\partial\over\partial\zeta}\hat a(\zeta)b$, where $\partial b(z)={d\over dz}b(z)$. Applying the next lemma to the algebra $A'\subset QO(A)$, the vector $1\in A$, and the linear map $\partial:A\ra A$, we find that $\chi$ is an isomorphism with inverse $\rho$. (In particular this shows that $A'=\rho(A)$, hence $\rho(A)$ is closed under the circle products.)
By Remark \RemarkCreation, we have
$$
\chi(\hat a\circ_n\hat b)=a\circ_n b
$$
for all $n$. Applying $\rho$ to both sides yields that $\hat
a\circ_n\hat b=\widehat{a\circ_nb}$ for all $n$. This shows that
$\rho$ preserves all circle products. $\Box$

\lemma{Let $A\subset QO(V)$ be a commutative circle algebra,
$\One\in V$ a nonzero vector, and $D:V\ra V$ a linear map such
that $D\One=0=a_+\One$ and $[D,a(z)]=\partial a(z)$ for $a\in A$.
If the creation map $\chi:A\ra V$, $a\mapsto a(-1)\One$, is
surjective then it is injective.}
\thmlab\ChiInjective
\proof By assumption, for $a\in
A$, we have $Da(n)\One=-na(n-1)\One$. Thus if $a(-1)\One=0$, then
$a(-2)\One=0$. Likewise $a(n)\One=0$ for all $n<0$. Since
$a_+\One=0$, it follows that $a\One=0$. Since $\chi$ is surjective
it suffices to show that $a(z)b(-1)\One=0$ for arbitrary $b\in A$.
Fix $N\geq0$ with $(z-w)^N[a(z),b(w)]=0$. Then
$$
(z-w)^Na(z)b(w)\One=(z-w)^Nb(w)a(z)\One=0.
$$
Since $b_+\One=0$, we have $b(w)\One\ra b(-1)\One$ as $w\ra0$. This shows that $z^N a(z)b(-1)\One=0$, implying that $a(z)b(-1)\One=0$. $\Box$

The following are useful identities for circle commuting operators
which measure the non-associativity and non-commutativity of the Wick product, and the
failure of the positive circle products to be left and right derivations of the
Wick product.

\lemma{Let $a,b,c$ be pairwise circle commuting, and $n\geq0$. Then we have the identities
$$\eqalign{
&:(:ab:)c:-:abc:=\sum_{k\geq0}{1\over(k+1)!}\left(:(\partial^{k+1}a)(b\circ_k
c): +(-1)^{|a||b|}:(\partial^{k+1}b)(a\circ_k c):\right)\cr
&a\circ_n(:bc:)-:(a\circ_nb)c:-(-1)^{|a||b|}:b(a\circ_nc):=
\sum_{k=1}^n\left(\matrix{n\cr
k}\right)(a\circ_{n-k}b)\circ_{k-1}c\cr 
&(:ab:)\circ_n c=\sum_{k\geq0}{1\over k!}:(\partial^ka)(b\circ_{n+k}c):
+(-1)^{|a||b|}\sum_{k\geq0}b\circ_{n-k-1}(a\circ_k c)\cr
&:ab:-(-1)^{|a||b|}:ba:=\sum_{k\geq0}{(-1)^k\over(k+1)!}\partial^{k+1}(a\circ_kb).
}$$}
\proof\thmlab\Associator
By the preceding theorem, it suffices to
show that $\hat a,\hat b,\hat c$ satisfy these identities. They
can be checked as follows. First, apply the creation map $\chi$ to
both sides and use \chieqn~
and Lemma \Commutator.
The calculations are straightforward, and details are left to the reader. $\Box$

Let $A$ be a commutative circle algebra. A two-sided ideal of circle algebra $A$ is a subspace $I$ invariant under left and right operations by the circle products. In this case, there is a canonical homomorphism
$$
A\ra QO(A/I),~~~ a\mapsto \bar a(\zeta), ~
\bar a(n)(b+I)=\hat a(n)b+I=a\circ_n b+I.
$$
This preserves the circle products, since the preceding theorem says that $A\ra QO(A)$, $a\mapsto\hat a$, is a circle algebra homomorphism. Likewise for $a,b\in A$, we have that $\bar a,\bar b$ circle commute. 
Thus the image $\bar A$ of $A$ in $QO(A/I)$ is a commutative circle algebra, and we have an exact sequence $0\ra I\ra A\ra\bar A\ra 0$. We call $\bar A$ the {\it quotient algebra} of $A$ by $I$.

\theorem{If $A$ is a commutative circle algebra, then $(A,1,\partial,\rho)$ is a vertex algebra in the sense of \FLM~(without grading or Virasoro element).}
\proof
We know that the map $\rho:A\ra QO(A)$, $a\mapsto\hat a$, has the property that $[\partial,\hat a(\zeta)]b={\partial\over\partial\zeta}\hat a(\zeta)b$. Moreover $\partial 1=0$ and that
$\chi:\rho(A)\ra A$, $\hat a\mapsto \hat a(-1)1=a$, is the inverse of $\rho$. So it remains to verify the vertex algebra Jacobi identity:
\eqn\dumb{
Res_\zeta(\widehat{\hat a(\zeta)b})(w)\zeta^n(w+\zeta)^q
=Res_z \hat a(z)\hat b(w)(z-w)^n z^q-Res_z\hat b(w)\hat a(z)(-w+z)^n z^q
}
for $n,q\in\Z$, $a,b\in A$. We will do this in several steps.

{\it Case 1.} $n\in\Z$, $q=0$. The identity
$$
\widehat{a\circ_n b}=\hat a\circ_n\hat b
$$
is nothing but \dumb~in this case. For convenience, we will drop the $\hat{}$ from the notations temporarily.

{\it Case 2.} $n=0$, $q=-1$.
The right side of \dumb~becomes, using Lemma \Commutator,
$$
[a(-1),b]=\sum_{p\geq0}(-1)^p(a\circ_p b)(w)~w^{-p-1},
$$
which agrees with the left side of \dumb.

{\it Case 3.} $n=-1$, $q=-1$. By direct computation, the right side of \dumb~is
$$\eqalign{
\sum_{p\geq0}a(-p-2)b(w)w^p-\sum_{p\geq0}b(w)a(p-1)w^{-p-1}
&=(a_-b-a(-1)b)w^{-1}+(ba_++ba(-1))w^{-1}\cr
&=:ab:w^{-1}-[a(-1),b]w^{-1}.
}$$
This agrees with the left side of \dumb.

{\it Case 4.} $n=0$, $q<0$. Using integration by parts, the first term of the right side of \dumb~becomes
$$
Res_z a(z)b(w) z^q={-1\over q+1}Res_z \partial a(z) b(w) z^{q+1}
$$
for $q<-1$. Likewise for two other terms in \dumb. Thus this case
can be reduced to Case 2.

{\it Case 5.} $n\in\Z$, $q\geq0$. We have
$$\eqalign{
(z-w)^n z^q&=(z-w)^{n+1} z^{q-1}+(z-w)^nz^{q-1} w\cr
(-w+z)^nz^q&=(-w+z)^{n+1} z^{q-1}+(-w+z)^nz^{q-1}w\cr
\zeta^n(w+\zeta)^q&=\zeta^{n+1}(w+\zeta)^{q-1}+\zeta^n(w+\zeta)^{q-1}w
}$$
Using these identities, we can easily reduce this case to Case 1.

{\it Case 6.} $n\geq0$, $q<0$. We have
$$\eqalign{
(z-w)^n z^q&=(z-w)^{n-1} z^{q+1}-(z-w)^{n-1}z^qw\cr
(-w+z)^n z^q&=(-w+z)^{n-1} z^{q+1}-(-w+z)^{n-1}z^qw\cr
\zeta^n(w+\zeta)^q&=\zeta^{n-1}(w+\zeta)^{q+1}-\zeta^{n-1}(w+\zeta)^qw
}$$
Using this identities, we reduce this case to Case 4.

{\it Case 7.} $n<0$, $q=-1$.  Take \dumb~in Case 3, and operate on both sides by ${d\over dw}$
repeatedly. We then get \dumb~for $n<0$, $q=-1$.

{\it Case 8.} $n<0$, $q<0$. Using integration by parts, the first term of the right side of \dumb~becomes
$$
Res_z a(z)b(w) (z-w)^nz^q={-1\over q+1}Res_z \partial a(z) b(w) (z-w)^nz^{q+1}
+{-n\over q+1}Res_z a(z) b(w) (z-w)^{n-1} z^{q+1}.
$$
for $q<-1$. Likewise for two other terms in \dumb. Now we verify that this case reduces to
Case 7.

This completes the proof. $\Box$

\lemma{If $(V,\One,D,Y)$ is a vertex algebra, then $Y(V)\subset
QO(V)$ is a commutative circle algebra.} \proof Write
$a(z)=Y(a,z)$, $b(z)=Y(b,z)$, for $a,b\in V$. The Jacobi identity
implies that $Y(a(p)b,z)=a(z)\circ_p b(z)$ for all $p$, which
shows that $Y(V)$ is closed under the circle products. The Jacobi
identity also implies that the commutator relations in Lemma
\Commutator~hold for $a(z),b(z)$, which shows that $a(z),b(z)$
circle commute. This shows that $Y(V)\subset QO(V)$ is a
commutative circle algebra. $\Box$

\remark{Thus the notion of a vertex algebra is abstractly equivalent to our notion of a commutative
circle algebra. While the former theory emphasizes the quadruple of structures $(V,\One,D,Y)$ satisfying an infinite family of (Jacobi) identities, the latter theory emphasizes the circle products and circle commutativity, and shows that all other structures can be obtained canonically in any given commutative circle algebra. The latter theory will be more convenient for the purposes of this paper. Note that the formal algebraic notions such as modules, ideals, and quotients  for vertex algebras \FHL~ are equivalent to the corresponding notions for commutative circle algebras under this dictionary. We will refer to a commutative circle algebra simply as a vertex algebra throughout the rest of the paper.}


The left regular module guarantees that for any given abstract
vertex algebra $A$, one can always embed $A$ in $QO(A)$ in a
canonical way. It is often convenient to pass between $A$ and its
image $\rho(A)$ in $QO(A)$. For example, we shall often denote the
Fourier modes $\hat a(n)$ simply as $a(n)$. Thus when we say that
a vertex operator $b(z)$ is annihilated by the Fourier mode $a(n)$
of a vertex operator $a(z)$, we mean that $a\circ_n b=0$. Here we
regard $b$ as being an element in the state space $A$, while $a$
operates on the state space, and the map $a\mapsto \hat a$ is the
state-operator correspondence.

Note that every commutative (super) algebra is canonically a
vertex algebra where any two elements strictly (graded) commute.
More generally we shall say that a vertex algebra is abelian if
any two elements pairwise commute. Otherwise we say that the
vertex algebra is non-abelian. If $a,b$ are two vertex operators
which commute, then their Wick product is the ordinary product and
we write $ab$ or $a(z)b(z)$.

\subsec{Examples}

We now give several constructions of known examples of vertex (super) algebras, all of which will be used extensively later.

\others{Example}{Current algebras.}
\thmlab\ExampleI
\noi Let $\gg$ be a Lie algebra equipped with a symmetric
$\gg$-invariant bilinear form $B$, possibly degenerate. The loop
algebra of $\gg$ is defined to be
$$
\gg[t,t^{-1}] = \gg\otimes \C[t,t^{-1}],
$$
with bracket given by $[u t^n,v t^m]=[u,v] t^{n+m}$. The form $B$
determines a 1-dimensional central extension $\hat{\gg}$ of
$\gg[t,t^{-1}]$ as follows:
$$
\hat{\gg}= \gg[t,t^{-1}]\oplus \C\tau,
$$
with bracket
$$
[u t^n,vt^m]=[u,v] t^{n+m} + n B(u,v)\delta_{n+m,0}\tau.
$$
$\hat{\gg}$ is equipped with the $\Z$-grading $deg(ut^n)=n$, and
$deg(\tau)=0$. Let $\gg_\geq$ be the subalgebra of elements of
non-negative degree, and let
$$
N(\gg,B)=\gU\hat\gg\otimes_{\gg_\geq}\C
$$
where $\C$ is the $\gg_\geq$-module in which $\gg[t]$ acts by zero
and $\tau$ by 1. Clearly $N(\gg,B)$ is graded by the non-positive
integers. For $u\in\gg$, denote by $u(n)$ the linear operator on
$N(\gg,B)$ representing $ut^n$, and put
$$
u(z)=\sum_n u(n)z^{-n-1}.
$$
Then for $u,v\in\gg$, we get
$$\eqalign{
[u(z)_+,v(w)]&=B(u,v)(z-w)^{-2}+[u,v](w)(z-w)^{-1}\cr
[u(z)_-,v(w)]&=-B(u,v)(w-z)^{-2}+[u,v](w)(w-z)^{-1}. }$$ It
follows immediately that $(z-w)^2[u(z),v(w)]=0$. Thus the
operators $u(z)\in QO(N(\gg,B))$ generate a vertex algebra
\FZ\LtypeI\LZII, which we denote by $O(\gg,B)$. Consider the
vector $\One=1\otimes 1\in N(\gg,B)$, called the vacuum vector.

\lemma{\LtypeI~The creation map $\chi:O(\gg,B)\ra N(\gg,B)$,
$a(z)\mapsto a(-1)\One$, is an $O(\gg,B)$-module isomorphism.}
\proof
We sketch a proof. By Remark \RemarkCreation, we have $\chi(a\circ_n b)=a(n)\chi(b)$, hence $O(\gg,B)$ is an $O(\gg,B)$-module homomorphism.
Next, $U\hat\gg$ has a derivation defined by $D\tau=0$, $D (ut^n)=-n ut^{n-1}$, and it descends to a linear map on $N(\gg,B)$ such that $[D,u(z)]=\partial u(z)$. This implies that $[D,a]=\partial a$ for all $a\in O(\gg,B)$. Thus to show that $\chi$ is a linear isomorphism, it suffices to show that it is surjective, by Lemma \ChiInjective. But this follows from PBW (see below). $\Box$


It is convenient to identify the spaces $N(\gg,B)$ and $O(\gg,B)$ under this isomorphism.
Obviously $O(\gg,B)$ contains the iterated Wick products
$$
:u^{I_0}\partial u^{I_1}\cdots{\partial^p u}^{I_p}:
$$
where $u^I$ means the {\it symbol} $u_1(z)^{i_1}\cdots
u_d(z)^{i_d}$,
 $\partial u^I$ means the symbol $\partial u_1(z)^{i_1}\cdots \partial u_d(z)^{i_d}$,
 for a given multi-index $I=(i_1,..,i_d)$, and likewise for other multi-index monomials. Here the $u_1,..,u_d$ form a basis of $\gg$. Under the creation map the image of the iterated Wick products above are the vectors, up to nonzero scalars,
$$
u(-1)^{I_0}u(-2)^{I_1}\cdots u(-p-1)^{I_p}\One
$$
which form a PBW basis, indexed by $(I_0,I_1,I_2,...)$, of the
induced module $N(\gg,B)$. Note also that there is a canonical
inclusion of linear spaces $\gg\hra O(\gg,B)$, $u\mapsto u(z)$.

An even vertex operator $J$ is called a current if
$J(z)J(w)\sim\alpha~(z-w)^{-2}$ for some scalar $\alpha$. The
formula for $[u(z)_+,v(w)]$ above implies the more familiar OPE
relation
$$
u(z)v(w)\sim B(u,v)(z-w)^{-2}+[u,v](w)~(z-w)^{-1}.
$$
In particular each $u(z)$ is a current (hence the name current
algebra).
The vertex algebra $O(\gg,B)$ has the following universal property
\LtypeI. Suppose that $A$ is any vertex algebra and $\phi:\gg\ra
A$ is a linear map such that $\phi(u)(z)~\phi(v)(w)\sim
B(u,v)(z-w)^{-2}+\phi([u,v])(w)~(z-w)^{-1}$ for $u,v\in\gg$. Then
there exists a unique vertex algebra homomorphism $O(\gg,B)\ra A$
sending $u(z)$ to $\phi(u)(z)$ for $u\in\gg$. In particular, any
Lie algebra homomorphism $(\gg,B)\ra(\gg',B')$ preserving the
bilinear forms induces a unique vertex algebra homomorphism
$O(\gg,B)\ra O(\gg',B')$ extending $\gg\ra\gg'$. It is also known
\LtypeI~that any Lie algebra derivation $d:(\gg,B)\ra (\gg,B)$
induces a unique vertex algebra derivation (i.e. a graded
derivation of all circle products) ${\bf d}:O(\gg,B)\ra O(\gg,B)$
with $u(z)\mapsto (du)(z)$.

When $\gg$ is a finite-dimensional Lie algebra, $\gg$ possesses a
canonical invariant, symmetric bilinear form, namely, the Killing
form $\kappa(u,v) = Tr\big(ad(u)\cdot ad(v)\big)$. In this case,
the current algebra $O(\gg,\lambda\kappa)$ is said to have a {\it
Schwinger charge} $\lambda$ \Pol.

It is easy to see that if $B_1,B_2$ are bilinear forms on $\gg$,
and $M_1,M_2$ are $O(\gg,B_1)$-, $O(\gg,B_2)$-modules
respectively, then $M_1\otimes M_2$ is canonically an
$O(\gg,B_1+B_2)$-module. In particular, tensor products of
$O(\gg,0)$-modules are again $O(\gg,0)$-modules.

There is a verbatim construction for any Lie super algebra equipped with an invariant form.

\others{Example}{Semi-infinite symmetric and exterior algebras.}
\thmlab\ExampleIV
\noi
Let $V$ be a finite dimensional vector space. Regard
$V\oplus V^*$ as an abelian Lie algebra. Then its loop algebra has a one-dimensional central extension
by $\C\tau$ with bracket
$$
[(x,x')t^n,(y,y')t^m]=(\bra y',x\ket-\bra
x',y\ket)\delta_{n+m,0}\tau,
$$
which is a Heisenberg algebra, which we denote by $\gh=\gh(V)$.
Let $\gb\subset \gh$ be the subalgebra generated by $\tau$,
$(x,0)t^n$, $(0,x')t^{n+1}$, for $n\geq 0$, and let $\C$
be the one-dimensional $\gb$-module on which each $(x,0)t^n$,
$(0,x')t^{n+1}$ act trivially and the central element $\tau$ acts by
the identity. Consider the $\gU\gh$-module $\gU\gh\otimes_{\gb}\C$.
The operators representing $(x,0)t^n,(0,x')t^{n+1}$ on this module are
denoted by $\beta^x(n),\gamma^{x'}(n)$, and the Fourier series
$$
\beta^x(z)=\sum\beta^x(n)z^{-n-1},~~~\gamma^{x'}(z)=\sum\gamma^{x'}(n)z^{-n-1}\in
QO(\gU\gh\otimes_{\gb}\C)
$$
have the properties
$$
[\beta^x_+(z),\gamma^{x'}(w)]=\bra
x',x\ket(z-w)^{-1},~~~[\beta^x_-(z),\gamma^{x'}(w)]=\bra
x',x\ket(w-z)^{-1}.
$$
It follows that $(z-w)[\beta^x(z),\gamma^{x'}(w)]=0$. Moreover the
$\beta^x(z)$ commute; likewise for the $\gamma^{x'}(z)$. Thus the
$\beta^x(z),\gamma^{x'}(z)$ generate a vertex algebra $\cS(V)$.
This algebra was introduced in [FMS], and is known as a fermionic
ghost system, or a $\beta\gamma$-system, or a semi-infinite
symmetric algebra. By using the Lie algebra derivation
$D:\gh\ra\gh$ defined by $(x,0)t^n\mapsto -n(x,0)t^{n-1}$, $(0,x')t^{n+1}\mapsto-n(0,x')t^n$, $\tau\mapsto 0$, one can easily show, as in the case of $O(\gg,B)$, that
the creation map $\cS(V)\ra \gU\gh\otimes_{\gb}\C$, $a(z)\mapsto
a(-1)1\otimes 1$, is a linear isomorphism, and that the
$\beta^x,\gamma^{x'}$ have the OPE relation
$$
\beta^x(z)\gamma^{x'}(w)\sim\bra x',x\ket(z-w)^{-1}.
$$
By the PBW theorem, it is easy to see that the vector space
$\gU\gh\otimes_{\gb}\C$ has the structure of a polynomial algebra
with generators given by the negative Fourier modes
$\beta^x(n),\gamma^{x'}(n)$, $n<0$, which are linear in $x\in V$
and $x'\in V^*$.

We can also regard $V\oplus V^*$ as an odd abelian Lie (super)
algebra, and consider its loop algebra and a one-dimensional
central extension by $\C\tau$ with bracket
$$
[(x,x')t^n,(y,y')t^m]=(\bra y',x\ket+\bra
x',y\ket)\delta_{n+m,0}\tau.
$$
Call this $\Z$-graded algebra $\gj=\gj(V)$, and form the induced
module $\gU\gj\otimes_{\ga}\C$. Here $\ga$ is the subalgebra of
$\gj$ generated by $\tau$, $(x,0)t^n$, $(0,x')t^{n+1}$, for $n\geq
0$, and $\C$ is the one-dimensional $\ga$-module on
which $(x,0)t^n$, $(0,x')t^{n+1}$ act trivially and $\tau$ acts by
$1$. Then there is clearly a vertex algebra $\cE(V)$, analogous to
$\cS(V)$, and generated by odd vertex operators
$b^x(z),c^{x'}(z)\in QO(\gU\gj\otimes_{\ga}\C)$ with OPE
$$
b^x(z)c^{x'}(w)\sim\bra x',x\ket(z-w)^{-1}.
$$
This vertex algebra is known as a bosonic ghost system, or
$bc$-system, or a semi-infinite exterior algebra.
Again the creation map $\cE(V)\ra \gU\gj\otimes_{\ga}\C$,
$a(z)\mapsto a(-1)1\otimes 1$, is a linear isomorphism. As in the
symmetric case, the vector space $\gU\gj\otimes_{\ga}\C$ has the
structure of an odd polynomial algebra with generators given by
the negative Fourier modes $b^x(n),c^{x'}(n)$, $n<0$, which are
linear in $x\in V$ and $x'\in V^*$.


A lot of subsequent computations involve taking OPE of iterated Wick products of vertex operators in
$$
\cW(V):=\cE(V)\otimes\cS(V).
$$
There is a simple tool from physics, known as Wick's theorem, that allows us to compute $A(z)B(w)$ easily where each of $A,B$ has the shape $:a_1\cdots a_p:$ where $a_i$ is one of the generators of $\cW(V)$, or their higher derivatives. For an introduction to Wick's theorem, see \Kac. Here is a typical computation by Wick's theorem
by ``summing over all possible contractions'':
$$\eqalign{
(:a_1(z)a_2(z):)(:a_3(w)a_4(w):)&\sim \bra a_2a_3\ket\bra
a_1a_4\ket+(-1)^{|a_2||a_3|}\bra a_1a_3\ket\bra a_2a_4\ket\cr
&+\bra a_2a_3\ket:a_1(z)a_4(w):+:a_2(z)a_3(w):\bra a_1a_4\ket\cr +
(-1)^{|a_2||a_3|}&\bra a_1a_3\ket:a_2(z)a_4(w):+
(-1)^{|a_2||a_3|}:a_1(z)a_3(w):\bra a_2a_4\ket. }$$ Here the $a_i$
are homogeneous vertex operators with OPE $a_i(z)a_j(w)\sim \bra
a_ia_j\ket$, where the symbol $\bra a_ia_j\ket$ denotes something
of the form $const.~(z-w)^{-p}$ depending on $i,j$. To get the
final answer for the OPE, one formally expands each
$:a_i(z)a_j(w):$ on the right side above in powers of $(z-w)$,
i.e. replacing it formally by $:a_i(w)a_j(w):+:\partial a_i(w)
a_j(w):(z-w)+\cdots$.

Now let $\gg$ be a Lie algebra and $V$ be a finite-dimensional
$\gg$-module via the homomorphism $\rho:\gg\rightarrow End~V$.
Associated to $\rho$ is a $\gg$-invariant bilinear form $B$ on
$\gg$ given by $B(u,v) = Tr\big(\rho(u)\rho(v)\big)$.

\lemma{$\rho:\gg\rightarrow End~V$ induces a vertex algebra
homomorphism $\rho_\cS:O(\gg,-B)\rightarrow\cS(V)$.} \proof Let
$\rho^*:\gg\rightarrow V^*$ be the dual module, let $\bra,\ket$
denote the pairing between $V$ and $V^*$. Choose a basis
$x_1,\dots,x_n$ of $V$, and let $x_1',\dots,x_n'$ be the dual
basis. Put
$$
\Theta_\cS^u(z) = -:\beta^{\rho(u)x_i}(z)\gamma^{x'_i}(z):
$$
(summing over $i$, as usual). We need to show that the following
OPE holds: \eqn\dumb{ \Theta_\cS^u(z)\Theta_\cS^v(w)\sim -B(u,v)
(z-w)^{-2} + \Theta_\cS^{[u,v]}(w)(z-w)^{-1}. } By Wick's
theorem,
$$\eqalign{
\Theta_\cS^u(z)\Theta_\cS^v(w) &=
\left(-:\beta^{\rho(u)x_i}(z)\gamma^{x'_i}(z):\right)
\left(-:\beta^{\rho(v)x_j}(w)\gamma^{x'_j}(w):\right)\cr &=
-\langle\rho(u)x_i,x_j'\rangle \langle
\rho(v)x_j,x'_i\rangle(z-w)^{-2} -\langle \rho(v)x_j,x'_i\rangle
:\beta^{\rho(u)x_i}(w)\gamma^{x_j'}(w):(z-w)^{-1} \cr
&+\langle\rho(u)x_i,x'_j\rangle:\beta^{\rho(v)x_j}(w)\gamma^{x_i'}(w):(z-w)^{-1}
}$$ which yields the right side of \dumb. $\Box$

Likewise we have the fermionic analogues $\cE(V)$ of $\cS(V)$, and
$\Theta_\cE^u$ of $\Theta_\cS^u$ with
$$
\Theta_\cE^u(z)=:b^{\rho(u)x_i}c^{x_i'}:.
$$
 A verbatim computation with gives

\lemma{$\rho:\gg\rightarrow End~V$ induces a vertex algebra homomorphism
$\rho_\cE:O(\gg,B)\rightarrow\cE(V)$.}

Now let's specialize to the case where $V$ is the adjoint module
of $\gg$, where $\gg$ is a finite-dimensional Lie algebra. Then
$\cW(\gg)=\cE(\gg)\otimes\cS(\gg)$  is called the semi-infinite
Weil algebra of $\gg$. This algebra has been studied by numerous
authors (see e.g. \FF\AK\Getz). By the two preceding lemmas, we
have a vertex algebra homomorphism $O(\gg,0)\ra\cW(\gg)$, with
$u(z)\mapsto\Theta_\cE^u(z)\otimes 1+1\otimes\Theta_\cS^u(z)$ for
$u\in\gg$.

\others{Example}{Virasoro elements.}
\thmlab\ExampleIII
\noi
Let $A$ be a vertex algebra. We call a vertex operator $L\in A$ a Virasoro element if
$$
L(z)L(w)\sim{c\over2}(z-w)^{-4}+2L(w)~(z-w)^{-2}+\partial
L(w)~(z-w)^{-1}
$$
where $c$ is a scalar called the central charge of $L$. One often
further requires that $L(w)\circ_1$ acts diagonalizably on $A$ and
that $L(w)\circ_0$ acts by $\partial$. If these two conditions hold, then $A$, equipped with $L$, is called a conformal vertex algebra of central charge $c$. A vertex operator $a\in A$ is said to be primary (with
respect to $L$) of conformal weight $\Delta\in\C$ if
$$
L(z)a(w)\sim~\Delta~ a(w)~(z-w)^{-2}+\partial a(w)~(z-w)^{-1}.
$$
A vertex operator $a\in A$ is said to be quasi-primary of conformal weight $\Delta\in\Z_>$ if
$$
L(z)a(w)\sim~\alpha~(z-w)^{-\Delta-2}+\Delta~
a(w)~(z-w)^{-2}+\partial a(w)~(z-w)^{-1}
$$
for some scalar $\alpha$.

For example if $\gg$ is a finite-dimensional simple Lie algebra
then the vertex algebra $O(\gg,\lambda\kappa)$ has a Virasoro
element given by {\it Sugawara-Sommerfield formula}
$$
L(z)={1\over 2\lambda+1}\sum_i:x^i(z)x^i(z):
$$
where the $x^i$ is an orthornormal basis of $(\gg,\kappa)$. 
 This Virasoro element has central charge ${ 2\lambda~dim~\gg\over 2\lambda+1}$  provided of course that the denominator is nonzero.
(Note that we have chosen a normalization so that we need not explicitly mention the dual Coxeter number of $\gg$.)
More generally if $(\gg,B)$ is any finite dimensional Lie algebra with a non-degenerate
invariant form then $O(\gg,\lambda B)$ admits a Virasoro element $L$ for all but finitely many values of $\lambda$ \LtypeI. The Virasoro element above is characterized by property that for every $x\in\gg$, the vertex operator $x(z)$ is primary of conformal weight 1.



\others{Example}{Topological vertex algebras.}
\thmlab\ExampleV
\noi This notion was introduced in \LZI~Definition 3.4, where we call these objects TVA. It is an abstraction based on examples from physics (see e.g. \DVV\EY\LZI). A topological vertex algebra is a vertex algebra $A$ equipped with four distinguished vertex operators $L,F,J,G$, where $L$ is a Virasoro element with central charge zero, $F$ is an even current which is a conformal weight one quasi-primary (with respect to $L$), $J$ an odd conformal weight one primary with $J(0)^2=0$, and $G$ an odd conformal weight two primary, such that
$$
J(0)G=L,~~~F(0)J=J,~~~F(0)G=-G.
$$
In special cases, further conditions are often imposed, such as that $J(z)$ commutes with itself,
or that $G(z)$ commutes with itself, which we do not require here. Note also that we do not require that the Fourier mode $L(1)$ acts diagonalizably on $A$. There are numerous examples arising from physics. One of the simplest is given by 
$\cW(\C)$,
which has four generators $b,c,\beta,\gamma$. If we put (suppressing $z$) $L=-:b\partial c:+:\beta\partial\gamma:$,
$F=-:bc:$, $J=:c\beta:$, $G=:b\partial\gamma:$, then it is straightforward to check that they give
$\cW(\C)$ the structure of a TVA. The same vertex algebra $\cW(\C)$ supports many TVA structures. One can twist the one above by using the current $F$ to get some other TVA structures on $\cW(\C)$. Another such example is given in \Getz.

\subsec{Differential and graded structures}

The vertex algebras we consider here typically come equipped with a number of graded structures.
The $\Z/2$-graded structure on a vertex algebra often arises from 
a $\Z$-grading we call degree. The semi-infinite exterior algebra $\cE(V)$ is one such example, where the odd generators $b^x,c^{x'}$ are assigned degrees -1 and +1 respectively. Then $\cE(V)$ is a direct sum of subspaces consisting of degree homogeneous elements. Like the $\Z/2$-graded structure in this case, the degree structure on a vertex algebra is additive under the circle products. In general, we say that a vertex algebra $\cA$ is {\it degree graded} if it is $\Z$-graded $\cA=\oplus_{p\in\Z}\cA^p$, and the degree is additive under the circle products. We denote the degree by $deg=deg_\cA$.

In addition to the degree grading, the vertex algebras we consider often come equipped with another $\Z$-grading we call weight. In the example $\cE(V)$, the vertex operator $b^x,c^{x'}$ can be assigned weights 1 and 0 respectively. Then $\cE(V)$ is a direct sum of subspaces consisting of weight homogeneous elements. The weight structure is not additive under the circle products in this case. But rather, we have $wt(a\circ_n b)=wt(a)+wt(b)-n-1$. In general, we say that a vertex algebra $\cA$ {\it weight graded} if it is $\Z$-graded $\cA=\oplus_{n\in\Z}\cA[n]$, and the $n$th circle product has weight $-n-1$. We denote the weight by $wt=wt_\cA$. Note that there can be several different weight structures on the same vertex algebra.

We say that a vertex algebra $\cA$ is {\it degree-weight graded} if it is both degree and weight graded and the gradings are compatible, i.e. $\cA[n]=\oplus_{p\in\Z}\cA^p[n]$, $\cA^p=\oplus_{n\in\Z}\cA^p[n]$, where $\cA^p[n]=\cA^p\cap\cA[n]$.

As a consequence of Lemma \Associator, if a vertex algebra $\cA$ is weight graded and has no negative weight elements, then $\cA[0]$ is a commutative associative algebra with product $\circ_{-1}$ and unit 1. Almost all vertex algebras in this paper have this property.

If a vertex algebra $\cA$ comes equipped with a Virasoro element $L$ where $L\circ_1$ acts diagonalizably on $\cA$ with integer eigenvalues, then the eigenspace decomposition defines a weight grading on $\cA$.

We call a pair $(\cA,\delta)$ a differential vertex algebra if $\cA$ is a vertex algebra equipped with
a linear map which is vertex algebra derivation, i.e. a super-derivation of each circle product, such that $\delta^2=0$. If, furthermore, $\cA$ is degree graded, then we assume that $\delta$ is a degree +1 linear map. If, furthermore, $\cA$ is weight graded, then we assume that $\delta$ is a weight 0 linear map. The categorical notion of homomorphisms and modules of differential vertex algebras are defined in an obvious way.


\subsec{The commutant construction}

This is a way to construct interesting vertex subalgebras of a given vertex algebra, and it is the vertex algebra analogue of the commutant construction in
the theory of associative algebras. 

Let $A$ be a vertex algebra and $S\subset A$ any subset.
The commutant of $S$ in $A$ is the space
$$
Com(S,A)=\{a(z)\in A|\ b(z)\circ_n a(z) = 0,\ \forall b(z)\in S,\
n\geq 0\}.
$$
It is a vertex subalgebra of $A$: this follows from the fact that for any elements $a,b$ in a vertex algebra, we have $[b(z),a(w)]=0$ iff $b(z)\circ_n a(z)=0$ for $n\geq0$, which is an immediate
consequence of Lemma \Commutator. From this, it also follows that if $C$ is the vertex algebra generated by the set $S$, then
$$
Com(C,A)=Com(S,A).
$$
Clearly if $S\subset S'\subset A$, then we have $Com(S',A)\subset
Com(S,A)$.

The commutant subalgebra $Com(C,A)$ has a second interpretation.
It can be thought of as a vertex algebra analogue of the ring of
invariants in a commutative ring with a Lie group or a Lie algebra
action. First we can think of $A$ as a $C$-module via the left
regular action of $A$. Then $Com(C,A)$ is the subalgebra of $A$
annihilated by all $\hat c(n)$, $c\in C$, $n\geq0$. If $C$ is a
homomorphic image of a current algebra $O(\gg,B)$, then
$Com(C,A)=A^{\gg_\geq}$ where the right side is the subspace of
$A$ annihilated by $u(n)$, $u\in\gg$, $n\geq0$. The invariant
theory point of view of the commutant construction is developed in
\Linshaw. In our construction of the chiral equivariant cohomology
later, the commutant subalgebra will play the role of the
classical algebra of invariants in the classical equivariant
cohomology.

\subsec{A vertex algebra for each open set}

{\it Notations.} Here $U,U',V,V'$ will denote open sets in $\R^n$,
$\gamma:\R^n\ra\R$ an arbitrary linear coordinate, and
$\gamma^i:\R^n\ra\R$ the $i$-th standard coordinate. The space of
smooth complex-valued differential forms $\Omega(U)$ can be
thought of as the space of functions on a super manifold. Without
digressing into super geometry, it suffices to think of the linear
coordinates $\gamma^i$ (restricted to $U$) as even variables, and
the coordinate one-forms $c^i:=d\gamma^i$ as odd variables. We can
regard the $\beta^i={\partial\over\partial\gamma^i}$ as even
vector fields acting as derivations, and the
$b^i={\partial\over\partial c^i}$ as odd vector fields acting as
odd derivations, on the function space $\Omega(U)$.

Let $C=\C\{{\partial\over\partial\gamma^i}, {\partial\over\partial
c^i}:1\leq i\leq n\}=C_0\oplus C_1$ denote the $\C$-span of the
constant vector fields. It is an abelian Lie (super) algebra
acting by derivations on the commutative super algebra
$\Omega(U)$. We now apply the current algebra construction to the
semi-direct product Lie algebra
$$
\Lambda(U):=C\triangleright\Omega(U)
$$
equipped with the zero bilinear form 0. Thus we consider
the loop algebra $\Lambda(U)[t,t^{-1}]$ and its module
$$
\cV(U):=N(\Lambda(U),0)
$$
For any $x\in\Lambda(U)$, denote by $x(k)\in End~\cV(U)$ the operator
representing the $xt^k\in\Lambda(U)[t,t^{-1}]$, and form $x(w)=\sum x(k)w^{-k-1}\in QO(\cV(U))$.
Then the $x(w)$ generates the vertex (super) algebra $O(\Lambda(U),0)$ defined in Example \ExampleI. As before, we identify this as a linear space with $\cV(U)$. Recall that we have a linear inclusion $\Lambda(U)\hra\cV(U)$, $x\mapsto x(w)$.
In particular for any function $f\in C^\infty(U)\subset\Omega(U)\subset\Lambda(U)$
and any vector field $\beta\in C_0$, we have
$$
\beta(z)f(w)\sim\beta(f)(w)~(z-w)^{-1}.
$$

The vertex algebra $\cV(U)$ is too large. For example the constant
function 1 gives a vertex operator $1(w)$ which is not equal to
$id$. If $f,g$ are two smooth functions then the vertex operators
$(fg)(w)$ and $:f(w)g(w):$ are not the same. Let $\cI(U)$ be the
two-sided ideal in $\cV(U)$ generated by the vertex operators
$$
{d\over dw} f(w)-{d\over dw}\gamma^i(w)~{\partial f\over\partial\gamma^i}(w)-
{d\over dw} c^i(w)~{\partial f\over\partial c^i}(w),~~~~
(fg)(w)-f(w)g(w),~~~1(w)-id
$$
where $f,g\in\Omega(U)$. (As always, the repeated index $i$ is summed over $i=1,..,n$,
unless said otherwise.) We put
$$
\cQ(U):=\cV(U)/\cI(U).
$$
Note that $\Lambda(U)$ becomes a $\Z_+$-graded Lie super algebra if we declare $C,\Omega(U)$ to have weight 1,0 respectively. This induces
a $\Z$-grading on $\cQ(U)$, and a canonical surjection ${}^2$ $\cQ(U)[0]\ra \Omega(U)$.\footnote{}{We thank B. Song for pointing this out.
If we declare that $\beta^i(z),b^i(z),c^i(z)$ have weights 1,1,0 respectively, $f(z)$ has weight 0 for $f\in C^\infty(U)$, then $\cI(U)$ is homogeneous ideal in $N(\Lambda(U),0)$. Hence $\cQ(U)[0]\ra\Omega(U)$, $g(z)\mapsto g$, is well-defined and surjective.}
On $\cQ(U)$, we also have
$$
b^i(z)c^j(w)\sim\delta_{ij} (z-w)^{-1}.
$$

Let $\Gamma$ be the vertex algebra generated by $\beta(w)$ with $\beta\in C_0$,
and the $f(w)$ with $f\in C^\infty(U)$,  subject to the relations
$$
\beta(z)f(w)\sim\beta(f)(w)~(z-w)^{-1}
$$
and with
$$
{d\over dw} f(w)-{d\over dw}\gamma^i(w)~{\partial f\over\partial\gamma^i}(w),~~~~
(fg)(w)-f(w)g(w),~~~1(w)-id
$$
being set to zero for all $f,g\in C^\infty(U)$. Let $B$ be the
vertex algebra generated by the $b^i(w),c^i(w)$, subject to the
relations $b^i(z)c^j(w)\sim\delta_{ij} (z-w)^{-1}$.
We claim that there is a canonical isomorphism
$$
\cQ(U)\cong \Gamma\otimes B
$$
Since each $f\in\Omega(U)$ can be uniquely written as
$f_I c^I$ where $f_I\in C^\infty(U)$ and $c^I=c^{i_1}c^{i_2}\cdots$, $i_1<i_2<\cdots$,
we can define a map $\Lambda(U)\ra\Gamma\otimes B$, $f_Ic^I\mapsto f_I(z)\otimes c^{i_1}(z)c^{i_2}(z)\cdots$,
$\beta^i\mapsto \beta^i(z)\otimes 1$, $b^i\mapsto 1\otimes b^i(z)$. Since the image of $\Lambda(U)$ satisfies the expected OPE of the current algebra $O(\Lambda(U),0)\equiv N(\Lambda(U),0)$, by the universal property, we have a surjective map $N(\Lambda(U),0)\ra\Gamma\otimes B$.
It is easy to check that this map factors through the ideal of relations $\cI(U)$, hence we have a surjective map $\cQ(U)\ra\Gamma\otimes B$. Likewise, we have maps $\Gamma\ra\cQ(U), B\ra\cQ(U)$ whose images commute in $\cQ(U)$. This yields a map $\Gamma\otimes B\ra\cQ(U)$. Then we verify that this map is the inverse of $\cQ(U)\ra\Gamma\otimes B$ above.

We would like to write down a basis for $\Gamma$. For this, we will construct $\Gamma$ in a
different way. Define a $\Z$-graded Heisenberg Lie algebra $\gg$
by the relations
$$
[\beta_p,\gamma_q]=\delta_{p+q,0}\beta(\gamma)~\One
$$
where the generators $\beta_p$ are linear in $\beta\in C_0$,
and the $\gamma_q$ are linear in $\gamma\in Hom(\R^n,\R)$.
Let $\gg_\geq\subset\gg$ be the subalgebra generated by $\One$ and the $\beta_p,\gamma_p$
with $p\geq0$. We make $C^\infty(U)$ a $\gg_\geq$-module by
$$\eqalign{
&\beta_p\cdot f=\gamma_p\cdot f=0,~~~~p>0\cr
&\beta_0\cdot f=\beta(f),~~~\gamma_0\cdot f=\gamma f,~~~\One\cdot f=f.
}$$
Put
$$
\Gamma':=\gU\gg\otimes_{\gg_\geq} C^\infty(U).
$$
The commutation relations of the operators $\beta_p,\gamma_q$ acting on this $\gg$-module
translate into the equivalent relations
$$\eqalign{
&[\beta(z),\beta'(w)]=[\gamma(z),\gamma'(w)]=0,\cr
&[\beta(z)_+,\gamma(w)]=\beta(\gamma)~(z-w)^{-1},
~~~[\beta(z)_-,\gamma(w)]=\beta(\gamma)~(w-z)^{-1}.
}$$
for $\beta,\beta'\in C_0$ and $\gamma,\gamma'\in Hom(\R^n,\R)$, where $\beta(w)=\sum\beta_n w^{-n-1}$, $\gamma(w)=\sum\gamma_nw^{-n}$. By Lemma \Commutator, it follows that
the operators $\beta(z),\gamma(z)\in QO(\Gamma')$ generate a vertex algebra.

Note that $\gg$ acts on the abelian Lie algebra
$C^\infty(U)[t,t^{-1}]$ by derivations defined by
$$
\beta_p\cdot ft^q=\beta(f)t^{p+q},~~~\gamma_p\cdot ft^q=\One\cdot ft^p=0.
$$
We will extend the $\gg$ action on $\Gamma'$ to an action of the
semi-direct product algebra $\gg\triangleright
C^\infty(U)[t,t^{-1}]$ on $\Gamma'$ as follows. {\it Having this
action is the main point of constructing the module $\Gamma'$.}

By using the existence of a PBW basis of $\gU\gg$, we will first
define a map $f(k): \gU\gg\otimes C^\infty(U)\ra\Gamma'$ for each
$ft^k\in C^\infty(U)[t,t^{-1}]$ inductively, and then show that
$f(k)$ descends to an operator $f(k):\Gamma'\ra\Gamma'$. Given
$f,g\in C^\infty(U)$, we define $f(k)(1\otimes g)=\delta_{k,-1}
(1\otimes fg)$ for $k\geq-1$. For $k<-1$, we put
$$
f(k)(1\otimes g)= {1\over k+1}\sum_{p<0}p~\gamma^i_p~{\partial
f\over\partial\gamma^i}(k-p)(1\otimes g).
$$
Note that this definition is recursive. The formula comes from
solving for the Fourier modes $f(k)$ using the anticipated
relation (from our earlier construction of $\Gamma$)
$$
{d\over dw} f(w)(1\otimes g) ={d\over dw}\gamma^i(w)~{\partial
f\over\partial\gamma^i}(w)(1\otimes g).
$$
The sought-after vertex operator $f(w)=\sum f(k)w^{-k-1}\in
QO(\Gamma')$ will eventually play the role of the $f(w)\in
QO(\Gamma)$ earlier. Now suppose $f(k)(\omega\otimes g)$ is
defined for all $f\in C^\infty(U)$ and all $k\in\Z$. We define
$$\eqalign{
f(k)(\gamma_p\omega\otimes g)&=\gamma_pf(k)(\omega\otimes g)\cr
f(k)(\beta_p\omega\otimes g)&=\beta_p f(k)(\omega\otimes
g)-\beta(f)(k+p)(\omega\otimes g). }$$ This completes the
definition of the $f(k):\gU\gg\otimes C^\infty(U)\ra\Gamma'$. Note
that when $f\in C^\infty(U)$ is the restriction of a given linear
function $\gamma\in Hom(\R^n,\R)$, then $f(k)=\gamma_{k+1}$, i.e.
$f(w)=\gamma(w)$ in this case.

Using the recursive definition above, it is straightforward but tedious to check that $f(k)$ descends to
an operator $f(k):\Gamma'\ra\Gamma'$, and that $f(w)=\sum f(k)w^{-k-1}\in QO(\Gamma')$ satisfies the relations
\eqn\ope{
[\gamma(z),f(w)]=0,~~[\beta(z)_+,f(w)]=\beta(f)(w)~(z-w)^{-1},~~
[\beta(z)_-,f(w)]=\beta(f)(w)~(w-z)^{-1}.
}
To see that $\Gamma'$ has a module structure
over the Lie algebra $\gg\triangleright C^\infty(U)[t,t^{-1}]$, it remains to show that
$[f(z),g(w)]=0$ for all $f,g\in C^\infty(U)$. This follows from the next lemma.

\lemma{For $f,g\in C^\infty(U)$, we have
$$
f(z)g(w)=g(w)f(z),~~~
(fg)(w)=f(w)g(w),~~~
1(z)=id,~~~
{d\over dw} f(w)={d\over dw}\gamma^i(w){\partial f\over\partial\gamma^i}(w).
$$
In particular $\Gamma'$ has a module structure over the Lie
algebra $\gg\triangleright C^\infty(U)[t,t^{-1}]$.}
\proof\thmlab\IdealRelations It follows from \ope~ that the
$\gamma_p$ commute with the commutator $[f(z),g(w)]$ for any
$\gamma\in Hom(\R^n,\R)$. It follows from the recursive definition
of the $f(k)$ that $[f(z),g(w)]1\otimes h=0$ for any $h\in
C^\infty(U)$. For $\beta\in C_0$, \ope~ implies that
$$
[\beta_p,[f(z),g(w)]]=z^p[\beta(f)(z),g(w)]+w^p[f(z),\beta(g)(w)],~~~~[\gamma_p,[f(z),g(w)]]=0.
$$
Using these commutator relations and the existence of a PBW basis of $\gU\gg$, we find by induction that $[f(z),g(w)]$ must be identically zero on $\Gamma'$. This proves the first asserted equation.
The argument for each of the remaining three equations is analogous. $\Box$

\corollary{The circle algebra $G'\subset QO(\Gamma')$ generated by
the operators $\{\beta(w),f(w)| \beta\in C_0,~f\in C^\infty(U)\}$
is a vertex algebra. It is linearly isomorphic to $\Gamma'$.}
\proof By the first equation of the preceding lemma together with
\ope, it follows that $(z-w)[\beta(z),f(w)]=0$ and
$[f(z),g(w)]=0$. This implies that $G'$ is a vertex algebra.
Moreover we have $\beta_k 1_{\Gamma'}=f(k)1_{\Gamma'}=0$ for
$k\geq0$, where $1_{\Gamma'}:=1\otimes
1\in\gU\gg\otimes_{\gg_\geq}C^\infty(U)$. It follows that the
creation map
$$
G'\ra\Gamma',~~~a(w)\mapsto \lim_{w\ra0}a(w)1_{\Gamma'}=a(-1)1_{\Gamma'}
$$
is a well-defined linear map. 

We claim that this is a linear isomorphism. In fact, it follows from the preceding lemma,
that $G'$ is spanned by the vertex operators, each having the shape
$$
:\beta^{I_0}{d\beta\over dw}^{I_1}{d^2\beta\over dw^2}^{I_2}\cdots{d\gamma\over dw}^{J_1}
{d^2\gamma\over dw^2}^{J_2}\cdots f_\alpha(w):\in G'.
$$
Its image under the creation map is a nonzero scalar (given by
products of factorials) times the vector
$$
\beta(-1)^{I_0}\beta(-2)^{I_1}\cdots\gamma(-2)^{J_1}\gamma(-3)^{J_2}\cdots\otimes f_\alpha\in\Gamma'.
$$
Here $\{f_\alpha\}$ is a given basis of $C^\infty(U)$, and
$\beta(k)=\beta_k$, $\gamma(k)=\gamma_{k+1}$;  ${d^k\beta\over
dw^k}^I$ means the usual ${d^k\over dw^k}\beta^{i_1}(w){d^k\over
dw^k}\beta^{i_2}(w)\cdots$ for any given finite list
$I=\{i_1,i_2,...\}$ of indices ranging over $\{1,..,n\}$; likewise
for other multi-index monomials. By the PBW theorem, these vectors
form a basis of $\Gamma'$ indexed by
$(I_0,I_1,..,J_1,J_2,...\alpha)$. This implies that $G'\ra\Gamma'$
is a linear isomorphism. $\Box$

{\it From now on, we identify $G'$ with $\Gamma'$ via this
isomorphism.}

\corollary{$\Gamma'\otimes B$ is canonically a $\Lambda(U)[t,t^{-1}]$-module
such that $xt^p\cdot1_{\Gamma'}\otimes 1_B=0$ for $x\in\Lambda(U)$, $p\geq0$.}
\proof
An element $x\in\Lambda(U)=C\triangleright\Omega(U)$ can be uniquely written as
$x=\beta+b+\sum f_I c^I$ where $\beta\in C_0$, $b=\lambda_i b^i\in C_1$, $f_I\in C^\infty(U)$,
$c^I=c^{i_1}c^{i_2}\cdots,~i_1<i_2<\cdots$, as before. We define the linear operator $x(p)$
representing $xt^p$ on $\Gamma'\otimes B$ to be the $p$-th Fourier mode of the vertex operator
$$
x(w)=\beta(w)+b(w)+\sum f_I(w)c(w)^I
$$
Here it is understood that $b(w)=\lambda_i b^i(w)$ and the $c^i(w)$ act on the factor $B$
while $\beta(w)$ and the $f_I(w)$ act on the factor $\Gamma'$.
Taking a second element $x'=\beta'+b'+\sum f'_{I'}c^{I'}\in\Lambda(U)$, and using \ope~ together with
the first equation of the preceding lemma, we find that
$$
x(z)x'(z')\sim [x,x'](z')(z-z')^{-1}
$$
where $[x,x']$ is the Lie bracket in $\Lambda(U)$. This shows that $xt^p\mapsto x(p)$
defines a $\Lambda(U)[t,t^{-1}]$-module structure on $\Gamma'\otimes B$. Finally,
by construction, the $\beta(w),b(w),f_I(w),c^i(w)$ are vertex operators
whose $p$-th Fourier modes annihilates $1_\Gamma'\otimes 1_B$ for $p\geq0$.
It follows that the Fourier modes $x(p)$ of $x(z)$, which represent the $xt^p$, have the same property. $\Box$

\corollary{We have a vertex algebra isomorphism $\cQ(U)=\Gamma\otimes B\ra
\Gamma'\otimes B=G'\otimes B$
which sends $x(z)$ to $x(z)$ for $x\in\Lambda(U)$.}
\proof
By the universal property of the induced module $\cV(U)$, and by the preceding corollary,
there is a unique $\Lambda(U)[t,t^{-1}]$-module homomorphism sending $1\otimes 1\in\cV(U)$
to $1_{\Gamma'}\otimes 1_B\in\Gamma'\otimes B$. Since we identify $\cV(U)$
with the vertex algebra generated by the $x(z)\in QO(\cV(U))$,
this map sends $x(z)\in QO(\cV(U))$ to $x(z)\in G'\otimes B'$. In particular,
it is a vertex algebra homomorphism. By the preceding lemma, this homomorphism factors through the ideal $\cI(U)\subset\cV(U)$, hence it descends to $\cQ(U)\ra G'\otimes B$.
By construction it is obvious that $\cQ(U)=\Gamma\otimes B$ is spanned by vertex operators of the shape
$$
:\beta^{I_0}{d\beta\over dw}^{I_1}{d^2\beta\over dw^2}^{I_2}\cdots{d\gamma\over dw}^{J_1}
{d^2\gamma\over dw^2}^{J_2}\cdots f_\alpha(w):\otimes
:b^{K_0}{db\over dw}^{K_1}{d^2b\over dw^2}^{K_2}\cdots c^{L_0}{dc\over dw}^{L_1}
{d^2c\over dw^2}^{L_2}\cdots:
$$
But their images form a basis of $G'\otimes B$. It follows that $\cQ(U)\ra G'\otimes B$
is an isomorphism. $\Box$

{\it From now on, we identify $\cQ(U)$ with $G'\otimes B$ via the isomorphism.}
But it will be convenient to use both points of view.

\corollary{There is a canonical map $\Omega(U)\hra\cQ(U)$ such
that $fg\mapsto :f(w)g(w):=f(w)g(w)$.} \proof\thmlab\Inclusion The
map is defined by $f=\sum f_Ic^I\mapsto f(w)=\sum f_I(w)c(w)^I\in
G'\otimes B$. It is clear that this is independent of basis. By
Lemma \IdealRelations, it has the desired multiplicative property.
If $f(w)=0$, then $f_I(w)=0$ for all $I$, since the
$c(w)^I=c^{i_1}(w)c^{i_2}(w)\cdots$ are independent in $B$. In
particular $f_I(-1)1_{\Gamma'}=1\otimes f_I=0$ in $\Gamma'$. It
follows that $f_I=0$ for all $I$. Thus the map $f\mapsto f(w)$ is
injective. $\Box$

\subsec{MSV chiral de Rham complex of a smooth manifold}

Following \MSV, the idea is to first construct a sheaf of vertex algebras on $\R^n$,
and then transfer it onto any given smooth manifold $M$ in a coordinate independent way.

\lemma{Any inclusion of open sets $U{\br\iota\over\subset} U'$
induces a canonical vertex algebra homomorphism
$\cQ(U'){\br\cQ(\iota)\over\ra}\cQ(U)$. Moreover, this defines a
sheaf of vertex algebras on $\R^n$} \proof Given an inclusion
$\iota:U\subset U'$, clearly we have a Lie algebra homomorphism
$\Lambda(U')\ra\Lambda(U)$ induced by restrictions of functions
$\Omega(U')\ra\Omega(U)$. By the functoriality of the current
algebra construction Example \ExampleI, we get a vertex algebra
homomorphism $\cV(U')\ra\cV(U)$. The ideal $\cI(U')$ is mapped
into $\cI(U)$ because $\cI(U')$ is generated by vertex operators
constructed from the circle products $\circ_{-1}$, $\circ_{-2}$
and the vertex operators $f(w)$, $f\in\Omega(U')$. Thus we have a
vertex algebra homomorphism $\cQ(U')\ra\cQ(U)$, which we denote by
$\cQ(\iota)$. A similar argument shows that given inclusions of
open sets
$U_1{\br\iota_1\over\subset}U_2{\br\iota_2\over\subset}U_3$, we
get $\cQ(\iota_2\circ\iota_1)=\cQ(\iota_1)\circ\cQ(\iota_2)$. This
shows that the assignment $U\rightsquigarrow\cQ(U)$ is a presheaf.

To see that $\cQ$ is a sheaf, suppose that $U_i{\br\iota_i\over\subset}U$ form
a covering of $U$ we need to show that the sequence
$$
0\ra\cQ(U)\rightarrow\prod_i\cQ(U_i)\rightrightarrows\prod_{i,j}\cQ(U_i\cap U_j)
$$
is exact. By the PBW theorem applied to $\cQ(U)=\Gamma'\otimes B$,
each element $a\in\cQ(U)$ can be uniquely represented in the shape
$a=\sum_P Pf_P$, where $f_P\in C^\infty(U)$, and $\{P\}$ a basis of the graded polynomial space
$$
\C[\beta^i(k),\gamma^i(k-1),b^i(k),c^i(k)~:~k\leq-1,~1\leq i\leq n].
$$
Given an inclusion $\iota:U\subset U'$,
the restriction map $\cQ(\iota):\cQ(U')\ra\cQ(U)$ sends $\sum P f_P$ to
$\sum P f_P|U$. Now our assertion follows from the exact sequence
$$
0\ra C^\infty(U)\rightarrow\prod_iC^\infty(U_i)\rightrightarrows\prod_{i,j}C^\infty(U_i\cap U_j)
$$
This completes the proof. $\Box$

In order to transfer the sheaf $\cQ$ from $\R^n$ to an arbitrary
smooth manifold, we must be able to compare the vertex algebras
$\cQ(U)$ under diffeomorphisms $U\ra U'$ of open sets. For this it
is convenient to enlarge the category $Open(\R^n)$ of open sets by
allowing any open embedding $U\hra U'$ of open sets to be a
morphism. We shall denote this new category by
$(Open(\R^n),\hra)$. This category is a special example of a
Grothendieck topology. Namely, if $U_i{\br\psi_i\over\hra}U$,
$i=1,2,$ are two morphisms then we declare the fiber product to be
$U_1\times_U U_2:=\psi_1(U_1)\cap\psi_2(U_2)$. We also declare any
collection of morphisms $U_i{\br\psi_i\over\hra} U$ to be a
covering if $\cup_i\psi_i(U_i)=U$. Note that any morphism
$U{\br\psi\over\ra} U'$ can be factorized as a diffeomorphism
followed by an inclusion
$U{\br\varphi\over\ra}\psi(U){\br\iota\over\subset}U'$. But there
may be another open set $W\supset U$ and a diffeomorphism
$W{\br\rho\over\ra}U'$ such that $\psi$ is factorized as $U\subset
W{\br\rho\over\ra}U'$.

\lemma{Any diffeomorphism of open sets $U{\br\varphi\over\ra} U'$
induces a canonical vertex algebra isomorphism
$\cQ(U'){\br\cQ(\varphi)\over\ra}\cQ(U)$. Moreover, given
diffeomorphisms of open sets
$U_1{\br\varphi_1\over\ra}U_2{\br\varphi_2\over\ra}U_3$, we get
$\cQ(\varphi_2\circ\varphi_1)=\cQ(\varphi_1)\circ\cQ(\varphi_2)$.}
\proof\thmlab\FunctorialityI  The diffeomorphism $\varphi$ induces
the pull-back isomorphism $\varphi^*:\Omega(U')\ra\Omega(U)$ on
forms. This induces a vertex algebra isomorphism from the
subalgebra $\bra f(z)|f\in\Omega(U')\ket\subset\cQ(U')$ onto $\bra
f(z)|f\in\Omega(U)\ket\subset\cQ(U)$ by $f(z)\mapsto
\varphi^*(f)(z)$. We would like to extend this to
$\cQ(U')\ra\cQ(U)$. First note that $\varphi^*$ does not extend to
a Lie algebra isomorphism between $\Lambda(U')$ and $\Lambda(U)$
in general. However it extends to a Lie algebra isomorphism
between two larger Lie algebras
$$
\varphi^*:Vect(U')\triangleright\Omega(U')\ra Vect(U)\triangleright\Omega(U)
$$
where $Vect(U')=\Omega(U')\{{\partial\over\partial\gamma^i},{\partial\over\partial c^i}: 1\leq i\leq n\}=\Omega(U')\otimes C$ is the Lie algebra of all smooth super derivations on $\Omega(U')$.
Since the constant vector fields $C$ can be regarded as a subalgebra of $Vect(U')$, $\varphi^*$ maps
$\Lambda(U')=C\triangleright\Omega(U')$ to a Lie subalgebra of $Vect(U)\triangleright\Omega(U)$.
In particular, one can express each of the pull-backs $\varphi^*{\partial\over\partial\gamma^i},
\varphi^*{\partial\over\partial c^i}, \varphi^*c^i,\varphi^*\gamma^i$,
uniquely in terms of ${\partial\over\partial\gamma^i},{\partial\over\partial c^i}, c^i,\gamma^i\in\Omega(U)$, and the coordinates $\varphi^i$ of $\varphi$. In fact,
we have
$$
\varphi^*\gamma^i=\gamma^i\circ\varphi=\varphi^i,~~~~\varphi^*c^i=\varphi^*d\gamma^i=
{\partial\varphi^i\over\partial\gamma^j}c^j
$$
and for any constant vector field $X\in C$, the vector field $\varphi^*X\in Vect(U)$ is determined by the condition that
$$
(\varphi^*X)(\varphi^*f)=[\varphi^*X,\varphi^*f]=\varphi^*[X,f]=\varphi^*(X(f)),~~~f\in\Omega(U').
$$
These pull-back expressions make sense as vertex operators if one
{\it formally} replaces the
${\partial\over\partial\gamma^i},{\partial\over\partial c^i}$,
$c^i, \gamma^i, \varphi^i\in\Omega(U)$, by their vertex operator
counterparts $\beta^i(z),b^i(z),c^i(z)$,
$\gamma^i(z),\varphi^i(z)\in\cQ(U)$. Equivalently, $\varphi^*$
determines an injective linear map $\Phi:\Lambda(U')\ra\cQ(U)$
with the property that $\Phi(f)=(\varphi^*f)(z)$ for
$f\in\Omega(U')$. In particular, by Lemma \IdealRelations,
$\Phi(fg)=\Phi(f)\Phi(g)$ for $f,g\in\Omega(U')$.

A remarkable result of \MSV~says that \eqn\PhiOPE{
\Phi(b^i)(z)~\Phi(c^j)(w)\sim \delta_{ij}(z-w)^{-1},~~~~
\Phi(\beta)(z)~\Phi(f)(w)\sim\Phi(\beta(f))(w)~(z-w)^{-1} } for
$\beta\in C_0,~f\in C^\infty(U')$, and all other OPE are trivial.
Since each element $x\in\Lambda(U')$ can be uniquely expressed as
$x=\beta+b+\sum f_I c^I$, it follows that
$\Phi(x)=\Phi(\beta)+\Phi(b)+\sum \Phi(f_I)\Phi(c^I)$. Taking a
second element $x'\in\Lambda(U')$, it follows from \PhiOPE~ that
$$
\Phi(x)(z)~\Phi(x')(z')\sim\Phi([x,x'])(z')~(z-z')^{-1}.
$$
By the universal property of the current algebra construction of $\cV(U')$, the map $\Phi$ extends to
a vertex algebra homomorphism $\cV(U')\ra\cQ(U)$. Note that
$\Phi$ maps functions $\Omega(U')\subset\Lambda(U')\hra\cV(U')$ to functions
$\Omega(U)\hra\cQ(U)$. Since $\Phi$ preserves the circle products,
it follows that $\Phi$ maps the ideal $\cI(U')\subset\cV(U')$ to
$\cI(U)$ which is zero in $\cQ(U)$. Thus we obtain a homomorphism
$\cQ(U')\ra\cQ(U)$ which we denote by $\cQ(\varphi)$.

Now consider two diffeomorphisms $U_1{\br\varphi_1\over\ra}U_2{\br\varphi_2\over\ra}U_3$.
It is clear that
$$
(\varphi_2\circ\varphi_1)^*:Vect(U_3)\triangleright\Omega(U_3)\ra Vect(U_1)\triangleright\Omega(U_1)
$$
coincides with $\varphi_1^*\circ\varphi_2^*$; this is functoriality of pull-back. However,
because the definition of $\cQ(\varphi_i)$ involves formal substitutions of
non-commuting vertex operators, it is not a priori clear the same should hold for the $\cQ(\varphi_i)$.
But \MSV~ showed (in the case when all three $U_i$ are the same
by direct calculation using generators; but the proof works here), that
$\cQ(\varphi_2\circ\varphi_1)=\cQ(\varphi_1)\circ\cQ(\varphi_2)$ indeed holds;
this result requires a certain anomaly cancellation.  In any case, specializing this
to $U{\br\varphi\over\ra}U'{\br\varphi^{-1}\over\ra}U$, we conclude that
$\cQ(\varphi)$ is an isomorphism. $\Box$

Actually, there is a slight ambiguity in the choice of the isomorphism $\Phi$ above. When the classical objects ${\partial\over\partial\gamma^i},{\partial\over\partial c^i}$, $c^i, \gamma^i, \varphi^i\in\Omega(U)$ are formally replaced by
their vertex operator counterparts $\beta^i(z),b^i(z),c^i(z)$, $\gamma^i(z),\varphi^i(z)\in\cQ(U)$,
care must be taken when deciding in which order the vertex operators should appear because they do not commute in general. This ambiguity occurs at just one place, namely the $\beta^i$. The classical transformation law for the derivation ${\partial\over\partial\gamma^i}$, under coordinate transformations $\tilde\gamma^i=g^i(\gamma)$, $\gamma^i=f(\tilde\gamma)$, is
$$
{\partial\over\partial\tilde\gamma^i}=
{\partial f^j\over\partial{\tilde\gamma}^i}{\partial\over\partial\gamma^i}+
{\partial^2 f^k\over\partial{\tilde\gamma}^i\partial{\tilde\gamma}^j}
{\partial g^j\over\partial\gamma^l}c^l {\partial\over\partial c^k}
$$
There is no ambiguity in the second term on the right side under
the formal replacement. But for the first term, we can have the
vertex operators $:{\partial
f^j\over\partial{\tilde\gamma}^i}\beta^j:$ or $:\beta^j{\partial
f^j\over\partial{\tilde\gamma}^i}:$, which are not equal in
general. Both choices give the right OPE for
$\tilde\gamma^i,\tilde\beta^j$, but only the second choice
guarantees the composition property in the preceding lemma.

\lemma{Consider the following commutative diagram in $(Open(\R^n),\hra)$:
$$\matrix{
V&\br\iota\over\subset&U\cr
\varphi'\da\hskip.2in& &\varphi\da\hskip.2in\cr
V'&\br\iota'\over\subset &U'}
$$
where the vertical arrows are diffeomorphisms.
This induces a commutative diagram under $\cQ$, i.e.
$$
\cQ(\varphi')\circ\cQ(\iota')=\cQ(\iota)\circ\cQ(\varphi).
$$}
\proof\thmlab\FunctorialityII The commutative diagram says that
$\varphi'$ is the restriction of $\varphi$ to $V$. The
$\cQ(\iota')$, $\cQ(\iota)$ are defined by restrictions of
(arbitrary) smooth functions; these operations are clearly
compatible with formal substitutions, hence with $\cQ(\varphi)$.
$\Box$

\definition{For any open embedding $U{\br\psi\over\hra} U'$, we define
$\cQ(\psi):\cQ(U')\ra\cQ(U)$ by $\cQ(\psi)=\cQ(\varphi)\circ\cQ(\iota)$
where $U{\br\varphi\over\ra}\psi(U){\br\iota\over\subset}U'$ is
the factorization of $\psi$ into a diffeomorphism followed by an inclusion.}

\lemma{The assigment
$\cQ:(Open(\R^n),\hra)\rightsquigarrow\cV\cA$, defines a sheaf of
vertex algebras in the Grothendieck topology on $\R^n$.} \proof
The proof has little to do with vertex algebras. Suppose that one
has a sheaf $\cF$ in the ordinary topology $(Open(\R^n),\subset)$
on $\R^n$, and that $\cF$ further assigns to each diffeomorphism
$U{\br\varphi\over\ra}U'$ an isomorphism
$\cF(\varphi):\cF(U')\ra\cF(U)$ in a way that Lemmas
\FunctorialityI~ and \FunctorialityII~ hold for $\cF$. Then $\cF$
defines a sheaf in the Grothendieck topology $(Open(\R^n),\hra)$.
Namely, one keeps the same assignment of objects $\cF(U)$ for open
sets $U$, and assign to each open embedding $\psi$ a morphism
$\cF(\psi)$ as in the preceding definition. Then one finds, by
straightforward checking, that

\bu $\cF:(Open(\R^n),\hra)\rightsquigarrow\cV\cA$ is a functor;
\bu if $U_i{\br\psi_i\over\hra}U$ is a covering, then
the sequence
$$
0\ra \cF(U)\ra\prod_i\cF(U_i)\rightrightarrows\prod_{i,j}\cF(U_i\times_U U_j)
$$
is exact.

Now apply this to $\cQ$. $\Box$

\lemma{Given any sheaf of vertex algebras $\cF:(Open(\R^n),\hra)\rightsquigarrow\cV\cA$ in the
Grothendieck topology on $\R^n$, we can attach, to every smooth manifold $M^n$,
a sheaf of vertex algebras $\cF_M$ in the ordinary topology of $M$.}
\proof
This again has little to do with vertex algebras.

Let $\cB$ be the set of all coordinate open subsets of $M$, i.e.
$O\in\cB$ iff there is a chart $\psi:O\ra\R^n$. Let $\cC_O$ be the
set of such charts, and let $\cG_O$ be the groupoid consisting of
objects $\cF(\psi(O))$, $\psi\in\cC_O$, and morphisms
$$
\cF(\psi(O)){\br g_{\psi\varphi}:=\cF(\psi\circ\varphi^{-1})\over\longrightarrow}\cF(\varphi(O)),~~~
\psi,\varphi\in\cC_O.
$$
By definition this groupoid has just one morphism
$g_{\psi\varphi}$ between any pair of objects. In particular,
these morphisms satisfy the ``cocycle condition'' that
$$
g_{\psi\varphi}g_{\rho\psi}=g_{\rho\varphi}.
$$
We define the ``average'' of $\cG_O$ by
$$
\bar\cG_O:=\{(v_\varphi)\in\prod_{\psi\in\cC_O}\cF(\psi(O))|v_\varphi=g_{\psi\varphi}v_\psi~
\forall\psi,\varphi\}.
$$
Note that each projection $\pi_\varphi:\prod_\psi\cF(\psi(O))\ra\cF(\varphi(O))$,
restricts to an isomorphism $\pi_\varphi:\bar\cG_O\ra\cF(\varphi(O))$. Hence
each tuple $(v_\varphi)\in\bar\cG_O$ is determined by any one of
its entries $v_\varphi$. Conversely, for any $\varphi\in\cC_O$
and any $x\in\cF(\varphi(O))$, there exists a unique tuple in $\bar\cG_O$,
which we denote by $\bar x$, such that $\pi_\varphi(\bar x)=x$.

Since a sheaf on $M$ is determined by its values
on a collection of open sets forming a basis of $M$,
it suffices to define $\cF_M(O)$ for $O\in\cB$. We define
$$
\cF_M(O)=\bar\cG_O~~~O\in\cB,
$$
and for $O\subset P$, we define the restriction map
$$
res_{PO}:\bar\cG_P\ra\bar\cG_O,~~~v\mapsto \ol{\cF(\iota)\pi_\psi(v)}.
$$
Here $\psi:P\ra\R^n$ is any given chart and $\iota:\psi(O)\subset\psi(P)$ is the inclusion.
It is straightforward to check that

\bu $res_{PO}$ is well-defined, i.e. independent
of the choice of $\psi$;
\bu  $\cF_M$ defines a sheaf on $M$. In other words, we have that
\bu if $O\subset P\subset Q\in\cB$, then
$$
res_{PO}~res_{QP}=res_{QO};
$$
\bu if $O_i\subset O\in\cB$ form a covering by open sets, then the sequence
$$
0\ra\cF_M(O)\ra\prod_i\cF_M(O_i)\rightrightarrows\prod_{i,j}\cF_M(O_i\cap O_j)
$$
is exact. $\Box$

Applying this to the sheaf $\cQ$, we obtain a sheaf $\cQ_M$ of vertex (super) algebras
for every smooth manifold $M$. This is the chiral de Rham sheaf.
Since $\Omega(U)\hra\cQ(U)$ by Corollary \Inclusion,
it follows that $\cQ_M$ contains the de Rham sheaf $\Omega_M$.
To summarize:

\theorem{\MSV~For every smooth manifold $M$, we have a sheaf
$\cQ_M$ of vertex algebras which contains $\Omega_M$ as a
subsheaf of vector spaces.}

When $M$ is fixed and no confusion arises, we shall denote the
chiral de Rham sheaf by $\cQ$ without the subscript. It was shown
further in \MSV~that for any $M$, the vertex algebra $\cQ(M)$ of
global sections contains a Virasoro element with central charge 0
given in local coordinates by:
$$
\omega_\cQ(z)=\omega_{bos}(z)+\omega_{fer}(z),~~~~
\omega_{bos}(z)=-:b^{i}(z)\partial c^{i}(z):,~~~
\omega_{fer}(z)=:\beta^{i}(z)\partial\gamma^{i}(z):,
$$ (cf. \FMS).
$\cQ(M)$ contains another vertex operator $g(z)$, which will also
be useful to us. It is defined locally by
$$g(z)=:b^i(z)\partial\gamma^i(z):.$$
When $M$ is Calabi-Yau, $\cQ(M)$ contains a topological vertex
algebra (TVA), as in Example \ExampleV, where $\omega_\cQ(z),
g(z)$ play the roles of $L(z),G(z)$, respectively. The
differential $d_\cQ$ of this TVA still makes sense when $M$ is not
Calabi-Yau, and the equation $[d_\cQ,g(z)]=\omega_\cQ$ still
holds. $d_\cQ$ is given by the zeroth Fourier mode of the vertex
operator
$$
:\beta^i(z)c^i(z):.
$$
Moreover $(\cQ,d_\cQ)$ has the structure of a complex of sheaves
containing the de Rham complex of sheaves $(\Omega,d_{dR})$. In
particular $(\Omega(M),d_{dR})$ is a subcomplex of
$(\cQ(M),d_\cQ)$. $\cQ(M)$ is a $\Z_{\geq}$-graded module over the
Virasoro algebra, where the grading is given by the eigenvalues of
the Fourier mode $\omega_\cQ(1)$ of the Virasoro element. Moreover
the eigenspace of zero eigenvalue is $\Omega(M)$. Since
$[d_\cQ,g(1)]=\omega_\cQ(1)$, it follows that the Fourier mode
$g(1)$ is a contracting homotopy for $d_\cQ$ in every eigenspace
of nonzero eigenvalue. This means that the chiral de Rham
cohomology, i.e. the cohomology of $(\cQ(M),d_\cQ)$, coincides
with the classical de Rham cohomology. 

\newsec{From Vector Fields on $M$ to Global Sections of $\cQ_M$}

Recall that to any given Lie algebra $\gg$, we can attach a Lie (super) algebra as follows.
Let $\gg_{-1}$ be the adjoint module of $\gg$, but declared to be an odd vector space.
We can then form the semi-direct product Lie algebra $\gs\gg:=\gg\triangleright\gg_{-1}$.
Define a linear map $d:\gs\gg\ra\gs\gg$, $(\xi,\eta)\mapsto(\eta,0)$, which is a square-zero odd derivation.
The result is an example of a {\it differential (graded) Lie algebra}
$(\gs\gg,d)$, i.e. a Lie algebra equipped with a square-zero derivation.

Assign to $\gs\gg$ the zero bilinear form, and consider the
current algebra $O(\gs\gg,0)$ defined in Example \ExampleI. Then
the Lie algebra derivation $d:\gs\gg\ra\gs\gg$ induces a vertex
algebra derivation
$$
{\bf d}:O(\gs\gg,0)\ra O(\gs\gg,0)
$$
such that $(\xi,\eta)(z)\mapsto(\eta,0)(z)$. Note that $(\gs\gg,d)$ is a $\Z$-graded Lie algebra where $\gg,\gg_{-1},d$ have degrees 0,-1,1 respectively. This makes $O(\gs\gg){\br def\over=}(O(\gs\gg,0),{\bf d})$ a degree graded differential vertex algebra where $(\xi,0)(z),(0,\xi)(z)$ have degrees $0,-1$ respectively. Note also that $O(\gs\gg)$ is also weight graded where $(\xi,\eta)(z)$ has weight 1.

\definition{An $O(\gs\gg)$-algebra is a degree-weight graded differential vertex algebra $(\cA^*,\delta)$ equipped with a homomorphism $\Phi_\cA:O(\gs\gg)\ra(\cA,\delta)$.
We shall often denote the $O(\gs\gg)$-structure on $\cA$ simply by the map $\gs\gg\ra\cA$,
$(\xi,\eta)\mapsto L_\xi(z)+\iota_\eta(z)$.
An $O(\gs\gg)$-algebra homomorphism is a differential vertex algebra homomorphism
$f:(\cA,\delta)\ra (\cA',\delta')$ such that $f\circ\Phi_\cA=\Phi_{\cA'}$. Likewise we have the categorical notions of $O(\gs\gg)$-modules and $O(\gs\gg)$-module homomorphisms.}
\thmlab\OsgDef

\subsec{From vector fields to $O(\gs\gX)$-algebras}

Consider $\gX=\gX(M)$, the Lie algebra of vector fields on $M$.
For $X\in\gX$ let
$$
L_X:\Omega^k(M)\ra\Omega^k(M),~~~\iota_X:\Omega^k(M)\ra\Omega^{k-1}(M)
$$
respectively be the Lie derivative and the interior multiplication
by $X$. They have the familiar (super) commutators:
$$
[L_X,L_Y]=L_{[X,Y]},~~~[L_X,\iota_Y]=\iota_{[X,Y]},~~~[\iota_X,\iota_Y]=0.
$$
Thus the map $\phi:\gs\gX=\gX\triangleright\gX_{-1}\ra Der~\Omega(M)$, $(X,Y)\mapsto L_X+\iota_Y$ defines an injective Lie algebra homomorphism. Since $\phi\circ d=d_{dR}\circ\phi$, it follows that $(\gs\gX,d)\cong(\phi(\gs\gX),d_{dR})$ as differential Lie algebras.
We shall identify these two algebras.

As before, we can consider the corresponding differential vertex
algebra $O(\gs\gX)=(O(\gs\gX,0),{\bf d})$. Since the invariant
bilinear form we have chosen for $\gs\gX$ is zero, the diagonal
map $\gs\gX\ra\gs\gX\oplus\gs\gX$ induces a Lie algebra
homomorphism of the corresponding loop algebras, and ultimately a
differential vertex algebra homomorphism $O(\gs\gX)\ra
O(\gs\gX)\otimes O(\gs\gX)$. This makes the tensor product of any
two $O(\gs\gX)$-algebras also an $O(\gs\gX)$-algebra; likewise for
modules.

Let $X\in\gX$. Write $X=f_i{\partial\over\partial\gamma^i}$ in some local coordinates
$\psi:O\ra\R^n$. Then $\iota_X=f_i{\partial\over\partial c^i}$ as a derivation on $\Omega(M)$.

\lemma{For $X\in\gX(M)$, there exists a global section
$\iota_X(z)$ of the chiral de Rham sheaf such that on a local
coordinate open set we have $\iota_X(z)=:f_i(z)b^i(z):$.} \proof
Note that $:f_i(z)b^i(z):=f_i(z)b^i(z)$ because the
$f_i(z),b^i(z)$ commute. It is enough to show that on two
overlapping coordinate open sets $O,\tilde O\subset M$, the two local
expressions for $\iota_X(z)$ agree on $O\cap\tilde O$. Denote the two
local expressions for $X$ by $f_i{\partial\over\partial\gamma^i}$
and $\tilde f_i{\partial\over\partial{\tilde\gamma}^i}$ respectively. Since
$X$ and $\iota_X$ are both globally defined on $M$, on the overlap
we have the relations
$\tilde f_i=f_j{\partial{\tilde\gamma}^i\over\partial{\gamma}^j}$,
${\partial\over{\partial\tilde c}^i}={\partial\over{\partial c}^j}
{\partial{\gamma}^j\over\partial{\tilde\gamma}^i}$. This means that
$f_i(z),b^i(z)\in\cQ(O)$ and $\tilde f_i(z),{\tilde b}^i(z)\in\cQ(\tilde O)$, when
restricted to $O\cap \tilde O$, are related by
$\tilde f_i(z)=f_j(z){\partial{\tilde \gamma}^i\over\partial{\gamma}^j}(z)$,
${\tilde b}^i(z)=b^j(z)
{\partial{\gamma}^j\over\partial{\tilde\gamma}^i}(z)$. It follows that
$$
\tilde f_i(z){\tilde b}^i(z)=f_j(z)b^j(z).
$$
Here we have used the identity (a relation coming from the ideal $\cI(U)$)
${\partial\gamma^j\over\partial{\tilde\gamma}^i}(z)
{\partial{\tilde\gamma}^i\over\partial{\gamma}^k}(z)
=({\partial\gamma^j\over\partial{\tilde\gamma}^i}
{\partial{\tilde\gamma}^i\over\partial{\gamma}^k})(z)
=\delta_{j,k}id$. This completes the proof. $\Box$

\remark{Even though $X=f_i{\partial\over\partial\gamma^i}$ is
globally defined as a vector field, the formal substitution
$:f_i(z)\beta^i(z):$ does not give a well-defined global section
of $\cQ(M)$. There are two reasons. The first is that as a
derivation on $\Omega(O)$, ${\partial\over\partial\gamma^i}$ does
not transform like a vector. The second is that
$f_i(z),\beta^i(z)$ do not commute as vertex operators in
$\cQ(O)$. Both of these make $:f_i(z)\beta^i(z):$ transform in a
complicated way and fail to be globally well-defined.}

Recall that the vertex algebra $\cQ(M)$ of global sections of the chiral de Rham sheaf
has a well-defined differential
$$
d_\cQ:\cQ(M)\ra\cQ(M),~~~ a(z)\mapsto (\beta^i(z)c^i(z))\circ_0 a(z).
$$
Note that the last expression is also the commutator of the zeroth
Fourier mode of $\beta^i(z)c^i(z)$ with $a(z)$. Thus we obtain a
global section
$$
L_X(z){\br def\over=}d_\cQ~\iota_X(z)
$$
in $\cQ(M)$. In local coordinates, $L_X(z)$ is given by
$$
 L_X(z) =:\beta^i(z)f^i(z): + :{\partial f^j\over\partial\gamma^i}(z)~c^i(z)~ b^j(z):.
 $$

\lemma{Both $\iota_X(z), L_X(z)$ are primary vertex operators of conformal
weight 1.} \proof\thmlab\PrimaryLiota This is a local calculation.
Recall that the Virasoro element $\omega_\cQ(z)\in\cQ(M)$ is
characterized by the fact that locally the vertex operators
$b(z),c(z),\beta(z),\gamma(z)$ are primary of conformal weights 1,0,1,0
respectively \FMS.
In particular, $f(z)$ is primary of conformal weight 0 for any $f\in C^\infty(M)$. That makes
$\iota_X(z)=:f_i(z)b^i(z):$ primary of conformal weight 1. Since $\omega_\cQ(z)$ is $d_\cQ$-exact
it commutes with $d_\cQ$.
It follows that $L_X(z)=d_\cQ~\iota_X(z)$
is also primary of conformal weight 1. $\Box$

\theorem{The differential vertex algebra $(\cQ(M),d_\cQ)$ is an $O(\gs\gX)$-algebra.}
\proof
Using the local formulas for $L_X(z),\iota_Y(w)\in\cQ(M)$, we get easily
$$
L_X(z)\iota_Y(w)\sim\iota_{[X,Y]}(w)~(z-w)^{-1}.
$$
Taking commutators on both sides with the zeroth Fourier mode
$D(0)$ of $D(z)=\beta^i(z)c^i(z)$, and recalling that
$L_Y(w)=d_\cQ~\iota_Y(w)=[D(0),\iota_Y(w)]$, we get
$$
L_X(z)L_Y(w)\sim L_{[X,Y]}(w)~(z-w)^{-1}.
$$
We also have $\iota_X(z)\iota_Y(w)\sim0$.
By the universal property of a current algebra (Example
\ExampleI), there is a unique vertex algebra homomorphism
$\Phi_\cQ:O(\gs\gX,0)\ra\cQ(M)$ such that, for any
$(X,Y)\in\gs\gX$, we have $(X,Y)(z)\mapsto L_X(z)+\iota_Y(z)$. By
definition, the differential for $O(\gs\gX)$ is given by ${\bf
d}:(X,Y)(z)\mapsto(Y,0)(z)$, while the differential for $\cQ(M)$
is $d_\cQ:L_X(z)+\iota_Y(z)\mapsto L_Y(z)$. This shows that
$\Phi_\cQ$ intertwines ${\bf d}$ and $d_\cQ$, hence we have a
differential vertex algebra homomorphism
$\Phi_\cQ:O(\gs\gX)\ra(\cQ(M),d_\cQ)$. $\Box$

The formula for $L_X(z)$ and the statement that $\cQ(M)$ is a
module over the current algebra $O(\gX,0)$ via $X\mapsto L_X(z)$,
also appear in \MalikovSchectman. \remark{Consider any given Lie subalgebra
$\gg\subset\gX$. Then we have a differential Lie subalgebra
$\gs\gg\subset\gs\gX$. By Example \ExampleI, this induces an
inclusion of differential vertex algebras $O(\gs\gg)\subset
O(\gs\gX)$. This makes every $O(\gs\gX)$-algebra canonically an
$O(\gs\gg)$-algebra, and likewise for modules.}

\subsec{From group actions to global sections}

Let $G$ be a Lie group with Lie algebra $\gg$, and $M$ a smooth
$G$-manifold, i.e. $M$ is equipped with an effective $G$-action.
Then the group homomorphism $G\ra Diff(M)$ induces an injective
Lie algebra homomorphism given by
$$
\gg\ra\gX(M),~~~\xi\mapsto X_\xi,~~~~X_\xi(x)={d\over dt}e^{t\xi}\cdot x|_{t=0}.
$$
Thus $\gg$ can be viewed as a Lie subalgebra of $\gX=\gX(M)$. We
shall denote $L_{X_{\xi}}$ simply by $L_\xi$, and likewise for
$\iota_\xi$. Now it follows immediately from the preceding remark
that we have

\theorem{Let $G$ be Lie group with Lie algebra $\gg$, and $M$
be a $G$-manifold. Then $(\cQ(M),d_\cQ)$ is canonically
an $O(\gs\gg)$-algebra.}

In particular, each $\xi\in\gg$ gives rise to two vertex operators
$L_\xi(z),\iota_\xi(z)\in\cQ(M)$.

\newsec{Classical Equivariant Cohomology Theory}

In this section we summarize the theory of classical equivariant
cohomology from the de Rham theoretic point of view. All material in this section is taken from the book of Guillemin-Sternberg \GS. This summary will be used as a guide for
formulating the vertex algebra analogue of the theory.

Let $G$ be a compact Lie group with Lie algebra $\gg$, and $M$ a
topological space equipped with an action of $G$ by
homeomorphisms. The equivariant cohomology of $M$, denoted by
$H^*_G(M)$, is defined to be the ordinary cohomology of the
quotient $(M\times \cE)/G$ where $\cE$ is any contractible
topological space on which $G$ acts freely. It is well-known that
this is independent of the choice of $\cE$. Furthermore, if the
$G$ action on $M$ is free then $H^*_G(M)=H^*(M/G)$. If the action
is not free, the quotient space $M/G$ may be pathological, and
$H^*_G(M)$ is the appropriate substitute for $H^*(M/G)$. From now
on, we let $M$ be a finite-dimensional manifold on which $G$ acts
by diffeomorphisms. Then there is a de Rham model for $H^*_G(M)$
together with an equivariant version of the de Rham theorem which
asserts the equivalence between this model and the topological
definition.


The $G$-action on $M$ induces a group homomorphism $\rho:G\ra
Aut~\Omega(M)$ and a differential Lie algebra homomorphism
$\gs\gg\ra Der~\Omega(M)$, $(\xi,\eta)\mapsto L_\xi+\iota_\eta$,
making $\Omega(M)$ a $(\gs\gg,d)$-module such that the following
identities hold
\eqn\Compatibility{\eqalign{
[d,\iota_\xi]&=L_\xi\cr \frac{d}{dt}\rho(e^{t\xi})|_{t=0} &=
L_{\xi}\cr \rho(a) L_{\xi}\rho(a^{-1}) &= L_{Ad_a(\xi)}\cr \rho(a)
\iota_{\xi}\rho(a^{-1}) &= \iota_{Ad_a(\xi)}\cr \rho(a) d
\rho(a^{-1}) &= d.
}}
Here $d=d_{dR}$. Note that $\Omega(M)$ is $\Z$-graded by the form degree and that operators $L_\xi,\iota_\xi,d$ have degrees $0,-1,1$ respectively.
In the terminology of \GS, ($\Omega(M),d_{dR})$ is an example of

\definition{A $G^*$-algebra is a $\Z$-graded differential commutative superalgebra $(A,d_A)$ on which $G$ acts by automorphisms $\rho:G\ra Aut(A)$ and $\gs\gg$ acts by (super) derivations
$\gs\gg\ra Der~A$, $(\xi,\eta)\mapsto L_\xi+\iota_\eta$, such that
\Compatibility~hold with $d=d_A$. A $G^*$-algebra morphism is an
algebra homomorphism which preserves the above structures in an
obvious way.}

To make sense of \Compatibility~ in this generality, one can either restrict the $G$-action to $G$-finite vectors, or give an appropriate notion of differentiation along a curve in $A$.

We often denote $(A,d_A)$ simply by $A$.

Now suppose $M$ is a principal $G$-bundle. Then the vector fields
$X_\xi$, $\xi\in\gg$, generate a $G$-invariant trivial subbundle
$V\subset TM$. Choose a $G$-invariant splitting $TM=V\oplus H$; we
get $T^*M=V^*\oplus H^*$. Corresponding to the choice of $H$ is a
canonical map \eqn\Orthogonality{\eqalign{
&\gg^*\hra\Omega^1(M),~~~\xi'\mapsto\theta^{\xi'}\cr
&\theta^{\xi'}|H=0,~~~
\iota_\xi\theta^{\xi'}=\theta^{\xi'}(X_\xi)=\bra\xi',\xi\ket. }}
Using the $G$-invariance of the splitting $T^*M=V^*\oplus H^*$, it
is easy to see that:
$$
L_\xi\theta^{\xi'}=\theta^{ad^*(\xi)\xi'}.
$$
The $\theta^{\xi'}$ are called connection one-forms, and the
two-forms
$\mu^{\xi'}=d\theta^{\xi'}+\half\theta^{ad^*(\xi_i)\xi'}\theta^{\xi_i'}$
are called the curvature forms of the $\theta^{\xi'}$, where
$ad^*:\gg\ra End~\gg^*$ is the coadjoint module.

\definition{A $G^*$-algebra $A$ is said to be of type $C$ if there is an inclusion
$\gg^*\hra A_1$, $\xi'\mapsto\theta^{\xi'}$ such that
$\iota_\xi\theta^{\xi'}=\bra\xi',\xi\ket$ and the image $C\subset
A_1$ is $G$-invariant.}

Every one-form $\omega\in\Omega^1(M)$ such that
$\iota_\xi\omega=0$ for all $\xi$ can be thought of as a section
of $H^*$, hence it is called a horizontal one-form. Likewise, any
form $\omega\in\Omega(M)$ satisfying the same condition is called
horizontal.

\definition{Let $A$ be a $G^*$-algebra. An element $\omega\in A$ is said to be horizontal
if $\iota_\xi\omega=0$ for all $\xi\in\gg$. It is called basic if
it is horizontal and $\rho(a)\omega=\omega$ for all $a\in G$. We denote
by $A_{hor}$ and $A_{bas}$, respectively, the subalgebras of
horizontal and basic elements in $A$.}

If $G$ is connected then the condition that $\rho(a)\omega=\omega$ can
be replaced by the equivalent condition that $L_\xi\omega=0$. In
particular, $A_{bas}$ is the $\gs\gg$-invariant subalgebra of $A$
in this case.

It is easy to see that $d(A_{bas})\subset A_{bas}$.
Thus $A_{bas}$ is a subcomplex of $A$ and its cohomology $H^*_{bas}(A)$
is well-defined. Moreover, if $\phi:A\rightarrow B$ is a morphism
of $G^*$-algebras, then $\phi(A_{bas})\subset B_{bas}$, and so $\phi$ induces a map
$$
\phi_{bas}:H^*_{bas}(A)\rightarrow H^*_{bas}(B).
$$
Guillemin-Sternberg define the equivariant cohomology of a
$G^*$-algebra $A$, in a way that is analogous to the topological
definition.

\definition{Let $E$ be any acyclic $G^*$-algebra of type $C$. For any
$G^*$-algebra $A$, its equivariant cohomology ring
$H^*_G(A)$ is defined to be $H^*_{bas}(A\otimes E)$. The usual rule of graded tensor product of two superalgebras applies here.}

Thus in the category of $G^*$-algebras, an acyclic $G^*$-algebra $E$ plays the role of a contractible space $\cE$ with free $G$-action.
This definition is shown to be independent of the choice of $E$. Moreover, the equivariant de Rham theorem asserts that for any $G$-manifold $M$,
$$
H^*_G(M) =H^*_G(\Omega(M)),
$$
where the right side is defined by taking $A=\Omega(M)$.

\subsec{Weil model for $H^*_G(A)$}

There is a natural choice for the acyclic $G^*$-algebra of type
$C$, namely, the Koszul algebra
$$
W(\gg)=\Lambda(\gg^*)\otimes S(\gg^*).
$$
This algebra is $\Z_\geq$-graded where the generators
$c^{\xi'}=\xi'\otimes 1, ~\xi'\in\gg^*,$ of the exterior algebra
$\Lambda(\gg^*)$, and the generators $z^{\xi'}=1\otimes\xi'$ of
the symmetric algebra $S(\gg^*)$, are assigned degrees 1 and 2
respectively. The $G$ action $\rho:G\ra Aut~W(\gg)$ is the action
induced by the coadjoint action on $\gg^*$, and the
$(\gs\gg,d)$-structure on $W(\gg)$ is defined on generators by the
formulas
$$\eqalign{
L_\xi c^{\xi'}&=c^{ad^*(\xi)\xi'},~~~L_\xi z^{\xi'}=z^{ad^*(\xi)\xi'}\cr
\iota_\xi c^{\xi'}&=\bra\xi',\xi\ket,~~~\iota_\xi z^{\xi'}=c^{ad^*(\xi)\xi'}\cr
d_Wc^{\xi'}&=z^{\xi'},  ~~~d_Wz^{\xi'}=0.
}$$
We have an inclusion $\gg^*\hra W(\gg)$, $\xi'\mapsto c^{\xi'}$. This defines a $G^*$-algebra structure of type $C$ on $W(\gg)$.
This $G^*$-algebra is acyclic because $d_W$ has a contracting homotopy $Q$
given by $Qz^{\xi'}=c^{\xi'}$, $Qc^{\xi'}=0$. In fact, we have
$$
[d_W,Q]=c^{\xi_i'}{\partial\over\partial c^{\xi_i'}}+
z^{\xi_i'}{\partial\over\partial z^{\xi_i'}}
$$
operating on the polynomial super algebra
$W(\gg)=\C[c^{\xi_i'},z^{\xi_i'}|i=1,..,dim~\gg]$, clearly
diagonalizably. This implies that the only $d_W$-cohomology occurs
in degree 0 and is one-dimensional. Here the $\xi_i'$ form the
dual basis of a given basis $\xi_i$ of $\gg$. Hence
$H^*_{bas}(A\otimes W(\gg))$ provides a model for the equivariant
cohomology of $A$, called the Weil model of $H_G^*(A)$.

There is very useful (and crucial for us) change of variables we can perform on $W(\gg)$.
Note that the $c^{\xi'}$ play the role of connection one-forms, and the corresponding elements playing the role of the curvature two-forms are
$$
\gamma^{\xi'}=z^{\xi'}+\half c^{ad^*(\xi_i)\xi'}c^{\xi_i'}.
$$
Note that they are homogeneous of degree 2 in $W(\gg)$. We can
view $W(\gg)$ as an algebra generated by the
$c^{\xi'},\gamma^{\xi'}$. The defining relations of the
$G^*$-algebra structure now become
\eqn\Wdef{\eqalign{ L_\xi
c^{\xi'}&=c^{ad^*(\xi)\xi'},~~~L_\xi\gamma^{\xi'}=\gamma^{ad^*(\xi)\xi'}\cr
\iota_\xi
c^{\xi'}&=\bra\xi',\xi\ket,~~~\iota_\xi\gamma^{\xi'}=0\cr
d_Wc^{\xi'}&=-\half c^{ad^*(\xi_i)\xi'}c^{\xi_i'}+\gamma^{\xi'},
~~~d_W\gamma^{\xi'}=\gamma^{ad^*(\xi_i)\xi'}c^{\xi_i'} }
}
The differential $d_W$ can be written as a sum $d_{CE}+d_K$ where
$d_{CE}$ is the Chevalley-Eilenberg differential of the Lie
algebra cohomology complex of $\gg$ with coefficients in the
module $S(\gg^*)$, and $d_K$ has the shape of a Koszul
differential:
\eqn\CEK{
d_{CE}=-c^{\xi_i'}\gamma^{\xi_j'}\beta^{[\xi_i,\xi_j]}-\half
c^{\xi_i'}c^{\xi_j'}b^{[\xi_i,\xi_j]},
~~~d_K=\gamma^{\xi_i'}b^{\xi_i}.
}
where the $b^\xi$ is an odd derivation on $W(\gg)$ defined by $b^\xi
c^{\xi'}=\bra\xi',\xi\ket$, $b^\xi\gamma^{\xi'}=0$, and the
$\beta^\xi$ is an even derivation on $W(\gg)$ defined by
$\beta^\xi c^{\xi'}=0$, $\beta^\xi\gamma^{\xi'}=\bra\xi',\xi\ket$.
The contracting homotopy $Q$ of $d_W$ in the new variables is
given by $Qc^{\xi'}=0$, $Q\gamma^{\xi'}=c^{\xi'}$.

The $(\gs\gg,d)$-module structure on $W(\gg)$ can be given in a more instructive way by rewriting
\Wdef-\CEK~ as follows. Introduce the Clifford-Weyl algebra
$$
Clifford(\gg)\otimes Weyl(\gg)
$$
to be the associative $\C$-superalgebra with odd generators $c^{\xi'},b^\xi$ and even generators $\gamma^{\xi'},\beta^\xi$, linear in $\xi'\in\gg^*$, $\xi\in\gg$, subject to the commutation relations
$$
[b^\xi,c^{\xi'}]=\bra\xi',\xi\ket=[\beta^\xi,\gamma^{\xi'}].
$$
Note that $d_{CE}+d_K$ given by the formulas \CEK~ can be thought of as an element in this algebra.
Moreover $W(\gg)$ becomes a module over this algebra where the $c^{\xi'},\gamma^{\xi'}$ act by left multiplications, and the $b^\xi,\beta^\xi$ act by derivations, as defined above.
There is a canonical Lie algebra homomorphism (with commutator as  the Lie bracket in the target)
\eqn\dumb{
(\gs\gg,d)\hra Clifford(\gg)\otimes Weyl(\gg)
}
defined by
$$\eqalign{
(\xi,\eta)&\mapsto\Theta_W^\xi+b^\eta,~~~~d\mapsto d_{CE}+d_K\cr
\Theta_W^\xi&:=\Theta_\Lambda^\xi+\Theta_S^\xi,~~~~
\Theta_\Lambda^\xi=-c^{\xi_i'}b^{[\xi,\xi_i]},~~~
\Theta_S^\xi=-\gamma^{\xi_i'}\beta^{[\xi,\xi_i]}.
}$$
Then, here is the main observation: the $(\gs\gg,d)$-module structure $(\gs\gg,d)\ra End~W(\gg)$ factors through the Clifford-Weyl algebra via the map \dumb. In other words, the operators $L_\xi,\iota_\xi$ defining the $G^*$-algebra structure on $W(\gg)$ can be explicitly represented by $\Theta_W^\xi,b^\xi$, regarded as operators on $W(\gg)$. The vertex algebra analogue of this structure will be crucial later.


From \Wdef, we find that $W(\gg)_{bas}=S(\gg^*)^G$, the space of $G$-invariants in $S(\gg^*)$.
From \CEK, we find that $d_W=0$ on $W(\gg)_{bas}$. It follows that
$$
H^*_{bas}(W(g))=H_G^*(\C)=S(\gg^*)^G.
$$

Since $W(\gg)$ is freely generated as an algebra by the variables
$c^{\xi_i'}$ and $d_Wc^{\xi_i'}$, $W(\gg)$ is easily seen to be an initial object in
the category of $G^*$-algebras of type $C$.

\theorem{Let $A$ be any $G^*$-algebra of type $C$. Then there
exists a $G^*$-algebra morphism $\rho:W(\gg)\rightarrow A$.
Furthermore, any two such morphisms are chain homotopic and hence
induce the same map from $S(\gg^*)^G\rightarrow H^*_{bas}(A)$.}

In particular, the notion of a $G^*$-algebra of type $C$ is equivalent to the
notion of a $G^*$-algebra $A$ which admits a $G^*$-algebra
morphism $\rho:W(\gg)\rightarrow A$. From now on, we will refer to
such algebras as $W(\gg)$-algebras (not to be confused with the
term $W^*$-algebra, which has a different meaning in \GS).

The preceding theorem shows that associated to any
$W(\gg)$-algebra $B$ is a canonical map $\kappa_G:
S(\gg^*)^G\rightarrow H^*_{bas}(B)$. Since for any $G^*$-algebra
$A$, $A\otimes W(\gg)$ is a $W(\gg)$-algebra, we have
$$
\kappa_G:
S(g^*)^G\rightarrow H^*_{bas}(A\otimes W(\gg))\cong H^*_G(A).
$$
This is called the Chern-Weil map. Consider the case where $A=\Omega(M)$ and $M=pt$. Using the Weil model, we see that
$$
H^*_G(\C)=H^*_{bas}(\Omega(pt)\otimes W(\gg)) =
H^*_{bas}(W(\gg)) = S(\gg^*)^G.
$$
Topologically, the map $\kappa_G: S(\gg^*)^G\rightarrow
H^*_G(M)$ is induced by $M\rightarrow pt$.

\subsec{Cartan model for $H^*_G(A)$}

Let $A$ be a $W(\gg)$-algebra and $B$ a $G^*$-algebra. Define
$$
\phi = c^{\xi_i'}\otimes\iota_{\xi_i}\in End(A\otimes B).
$$
Since $\phi$ is a derivation and $\phi^{n+1}=0$,
$n=dim~G$, it follows that $\Phi = exp(\phi)$ is a well-defined
automorphism of the commutative (super) algebra $A\otimes B$,
known as the Mathai-Quillen isomorphism.

\theorem{The Mathai-Quillen isomorphism satisfies
$$\eqalign{
\Phi (L_\xi\otimes 1 + 1\otimes L_\xi) \Phi^{-1} &= L_\xi\otimes 1
+1\otimes L_\xi\cr \Phi (\iota_\xi\otimes 1 + 1\otimes \iota_\xi)
\Phi^{-1} &= \iota_\xi\otimes 1 \cr \Phi d \Phi^{-1} &=
d-\gamma^{\xi_i'}\otimes \iota_{\xi_i} + c^{\xi_i'}\otimes
L_{\xi_i} }$$ where $d=d_A\otimes 1+1\otimes d_B$.} The second
relation shows that $\Phi((A\otimes B)_{hor})=A_{hor}\otimes B$.
Let's specialize to $A=W(\gg)$. Since
$d_W|W(\gg)_{hor}=c^{\xi_i'}L_{\xi_i}$ it follows that
$$
\Phi d \Phi^{-1}|W(\gg)_{hor}\otimes B = (c^{\xi_i'}\otimes
1)(L_{\xi_i}\otimes 1+1\otimes L_{\xi_i}) + 1\otimes d_B
-\gamma^{\xi_i'}\otimes \iota_{\xi_i}.
$$
Since $\Phi$ is $G$-equivariant, we have $\Phi((W(\gg)\otimes
B)_{bas})= (S(\gg^*)\otimes B)^G=:C_G(B)$. On $C_G(B)$, the
operator $L_a\otimes 1+1\otimes L_a$ is zero, so
$$
\Phi d \Phi^{-1}|C_G(B) = 1\otimes d_B - \gamma^{\xi_i'}\otimes
\iota_{\xi_i}=:d_G.
$$
In particular, $H_G^*(B)\cong H^*(C_G(B),d_G)$. The right  side
is called the Cartan model for $H_G^*(B)$.

\newsec{Chiral Equivariant Cohomology Theory}

For simplicity, we shall assume throughout that $G$ is a connected
compact Lie group. Under this assumption, the appropriate analogue
of a $G^*$-algebra is the notion of an $O(\gs\gg)$-algebra given
by Definition \OsgDef. Our theory can be easily modified to allow
disconnected $G$ by further requiring that an $O(\gs\gg)$-algebra
comes equipped with a compatible $G$-action, as in the classical
setting.

Using the classical Weil model of $H_G^*(A)$ as a guide, we will
define the chiral equivariant cohomology $\H_G^*(\cA)$ of an
arbitrary $O(\gs\gg)$-algebra $\cA$. We have seen that the chiral
de Rham complex $\cQ(M)$ of a $G$-manifold $M$ is an example of an
$O(\gs\gg)$-algebra. For $\cA=\cQ(M)$, the chiral equivariant cohomology of
$\H^*_G(\cQ(M))$ is a vertex algebra analogue of the classical equivariant cohomology of the $G^*$-algebra $\Omega(M)$. We will see that there is a
canonical inclusion $H_G^*(M)\hra\H^*_G(\cQ(M))$. 

Recall that an $O(\gs\gg)$-algebra is a $\Z$-graded differential vertex
algebra $(\cA^*,d_\cA)$ equipped with a differential vertex algebra
homomorphism $\Phi_\cA:O(\gs\gg)=(O(\gs\gg,0),{\bf
d})\ra(\cA,d_\cA)$, $(\xi,\eta)(z)\mapsto L_\xi(z)+\iota_\eta(z)$.

\definition{Let $\cI\subset O(\gs\gg,0)$ be the vertex subalgebra generated by the odd currents $(0,\xi)(z)$, $\xi\in\gg$. Let $\cA$ be a given $O(\gs\gg)$-algebra. We define the horizontal and basic subalgebras of $\cA$ to be respectively
$$
\cA_{hor}=Com(\Phi_\cA\cI,\cA),~~~~\cA_{bas}=Com(\Phi_\cA
O(\gs\gg,0),\cA).
$$
Thus $\cA_{hor}$ consists of $a(z)\in\cA$ which strictly commute with the elements $\iota_\xi(z)\in\cA$, and $\cA_{bas}$ consists of $a(z)\in\cA_{hor}$ which strictly commute with the elements $L_\xi(z)\in\cA$.}

Since $d_\cA~\iota_\xi(z)=L_\xi(z)$, and $d_\cA$ is a square-zero derivation
of all the circle products, it follows that $(\cA_{bas},d_\cA)$ is vertex algebra with a compatible
structure of a cochain complex. Its cohomology
$\H_{bas}^*(\cA)$ is therefore a vertex algebra, which we will
call the chiral basic cohomology of $\cA$.  An $O(\gs\gg)$-algebra
homomorphism $\phi:\cA\ra\cB$ sends $\cA_{bas}$ to $\cB_{bas}$, so
it induces a vertex algebra homomorphism
$$
\phi_{bas}:\H_{bas}^*(\cA)\ra \H_{bas}^*(\cB).
$$

Here is a small but important surprise: to get an induced
homomorphism on basic cohomology, one needs less than an
$O(\gs\gg)$-algebra homomorphism. That is because the notion of
basic subalgebra uses only half the $O(\gs\gg)$-algebra
structures.

\definition{Let $\cA,\cB$ be $O(\gs\gg)$-algebras.  A differential vertex algebra homomorphism
$\phi:\cA\ra\cB$ is said to be basic if
$\phi(\cA_{bas})\subset\cB_{bas}$. In particular, a basic
homomorphism induces a vertex algebra homomorphism $\phi_{bas}$ on
basic cohomology.}

\others{Example}{}
\thmlab\ExampleVI
\noi Let $\cA,\cB$ be $O(\gs\gg)$-algebras and consider the map
$$
\phi:\cA\ra\cA\otimes\cB,~~~a\mapsto a\otimes 1.
$$
It is obviously a vertex algebra homomorphism, and it respects the differentials because
$(1\otimes d_\cB)(a\otimes 1)=0$. Now note that $\cA_{bas}\otimes 1\subset(\cA\otimes\cB)_{bas}$,
again because $1\otimes L_\xi(n)(a\otimes 1)=1\otimes\iota_\xi(n)(a\otimes 1)=0$ for $n\geq0$. Thus $\phi$ is a basic homomorphism. But $\phi$ will not be an $O(\gs\gg)$-algebra homomorphism unless
$O(\gs\gg)$ acts trivially on $\cB$, in which case $1\otimes L_\xi(z)=1\otimes\iota_\xi(z)=0$.

\subsec{Semi-infinite Weil algebra}

We saw that a crucial ingredient in the Weil model of the classical equivariant cohomology is the Koszul algebra $W(\gg)$, and that via the Clifford-Weyl algebra, one can write down the $G^*$-algebra structure on $W(\gg)$ very explicitly. It turns out that in the vertex algebra setting, there is a natural algebra that unifies the Koszul and the Clifford-Weyl algebras into a single object. This is the semi-infinite Weil algebra.
$$
\cW=\cW(\gg)=\cE(\gg)\otimes\cS(\gg)
$$
which we introduced in Example \ExampleIV.

If we declare that the vertex operators
$b^\xi,c^{\xi'},\beta^\xi,\gamma^{\xi'},$ have the respective
degrees $-1,1,-1,1$, then this defines a $\Z$-valued bigrading
$$
\cW=\oplus_{p,q}\cW_{p,q}
$$
where $\cW_{p,q}=\cE^p\otimes\cS^q$ is the degree $(p,q)$ subspace. This bigrading turns out to come from the following two vertex operators in $\cE(\gg),\cS(\gg)$:
$$
j_{bc}(z)=-:b^{\xi_i}(z)c^{\xi_i'}(z):,~~~ j_{\beta\gamma}(z)=
:\beta^{\xi_i}(z)\gamma^{\xi_i'}(z):.
$$
Their zeroth Fourier modes are diagonalizable operators on $\cE(\gg),\cS(\gg)$ respectively with integer eigenvalues.  The eigenspaces $\cE^p,\cS^q$ are also called the subspaces of $bc$-number $p$ and $\beta\gamma$-number $q$. We put
$$
\cW^n=\oplus_{n=p+2q}\cW_{p,q}.
$$
There is yet another bigrading on $\cW$ which is compatible with the one above. If we declare that the vertex operators $\beta^\xi,\gamma^{\xi'},b^\xi,c^{\xi'}$, have the respective weights 1,0,1,0, then this defines a $\Z_\geq$-valued bigrading
$$
\cW=\oplus_{m,n\geq0}\cW[m,n]
$$
where $\cW[m,n]=\cE[m]\otimes\cS[n]$ is the weight $(m,n)$ subspace. This bigrading turns out to come from Virasoro elements in the vertex algebras $\cE(\gg),\cS(\gg)$. Put
$$
\omega_\cW(z)=\omega_\cE(z)+\omega_\cS(z),~~~~
\omega_\cE(z)=-:b^{\xi_i}(z)\partial c^{\xi_i'}(z):,~~~
\omega_\cS(z)=:\beta^{\xi_i}(z)\partial\gamma^{\xi_i'}(z):.
$$
An OPE calculation by Wick's theorem yields

\lemma{The vertex operators $\omega_\cE(z), \omega_\cS(z)$ are Virasoro elements of central charges $\mp 2 dim~\gg$ in the respective vertex algebras $\cE(\gg),\cS(\gg)$. Moreover, $\omega_\cE(z)$ is the unique Virasoro element such that $b^\xi(z),c^{\xi'}(z)$ are primary of conformal weight 1,0 respectively. Likewise $\omega_\cS(z)$ has a similar characterization in $\cS(\gg)$. The $\cE[m],\cS[n]$ are the respective eigenspaces of $\omega_\cE(1),\omega_\cS(1)$, of eigenvalues $m,n$. Moreover, $\omega_\cE(0)$, $\omega_\cS(0)$, act respectively on $\cE(\gg),\cS(\gg)$, as the derivation $\partial$.}
\thmlab\FMSVirasoro

Note that the subspace $\cE[0]$ consists of the vertex operators
which are polynomial in the (anti-commuting) vertex operators
$c^{\xi'}(z)$, and is canonically isomorphic to the
classical exterior algebra $\Lambda(\gg^*)$. Likewise $\cS[0]$ is
canonically isomorphic to the classical symmetric algebra
$S(\gg^*)$. It follows that $\cW[0,0]$ is nothing but a copy of
the classical Koszul algebra $W(\gg)$.

Define the vertex operators (suppressing the variable $z$):
\eqn\ThetaW{
\Theta_\cW^\xi=\Theta_\cE^\xi+\Theta_\cS^\xi,~~~~
\Theta_\cE^\xi=:b^{[\xi,\xi_i]}c^{\xi_i'}:,~~~~
\Theta_\cS^\xi=-:\beta^{[\xi,\xi_i]}\gamma^{\xi_i'}:.
}
\eqn\DJK{
D=J+K,~~~J=-:c^{\xi_i'}\gamma^{\xi_j'}\beta^{[\xi_i,\xi_j]}:
-\half :c^{\xi_i'}c^{\xi_j'}b^{[\xi_i,\xi_j]}:,~~~K=:\gamma^{\xi_i'}b^{\xi_i}:
}
The Fourier mode $J(0)$ is called the semi-infinite differential. The next four lemmas follows easily from direct computations by Wick's theorem.

\lemma{$J(0)^2=K(0)^2=D(0)^2=[K(0),J(0)]=0$.}
\thmlab\SquareZero

\corollary{The complex $(\cW^*,D(0))$ has the structure of a
double complex $(\cW^{*,*},J(0),K(0))$ where $\cW^{p,q}$ is
defined to be $\cW_{p-q,q}$. Thus we have
$$
D(0):\cW^n\ra\cW^{n+1},~~~
J(0):\cW^{p,q}\ra\cW^{p+1,q},~~~K(0):\cW^{p,q}\ra\cW^{p,q+1}.
$$}

\lemma{$D(0)b^\xi(z)=\Theta_\cW^\xi(z)$.}
\thmlab\VertexCartan

\lemma{The vertex operators $\Theta_\cE^\xi(z),\Theta_\cS^\xi(z)$ are characterized in their respective algebras $\cE(\gg),\cS(\gg)$ by the properties that they are the only weight one elements such that
$$\eqalign{
&\Theta_\cE^\xi(z)b^\eta(w)\sim b^{[\xi,\eta]}(w)~(z-w)^{-1},~~~
\Theta_\cE^\xi(z)c^{\eta'}(w)\sim c^{ad^*(\xi)\eta'}(w)~(z-w)^{-1}.\cr
&\Theta_\cS^\xi(z)\beta^\eta(w)\sim \beta^{[\xi,\eta]}(w)~(z-w)^{-1},~~~
\Theta_\cS^\xi(z)\gamma^{\eta'}(w)\sim \gamma^{ad^*(\xi)\eta'}(w)~(z-w)^{-1}.\cr
}$$}
\thmlab\AdjointCoadjoint

\lemma{The $\Theta_\cE^\xi$ are primary of conformal weight 1 with respect
to $\omega_\cE$. Likewise for the $\Theta_\cS^\xi$ with respect to
$\omega_\cS$.}
\thmlab\ThetaPrimary

\lemma{There is a vertex algebra homomorphism
$O(\gg,\kappa)\ra\cE(\gg)$ such that
$\xi(z)\mapsto\Theta_\cE^\xi(z)$. Likewise we have
$O(\gg,-\kappa)\ra\cS(\gg)$. Here $\kappa(\xi,\eta)=Tr(ad(\xi)ad(\eta))$, is the Killing form of $\gg$.
}
\proof\thmlab\OgModule We have
$$
\Theta_\cE^\xi(z)\Theta_\cE^\eta(w)\sim
\kappa(\xi,\eta)~(z-w)^{-2}+\Theta_\cE^{[\xi,\eta]}(w)~(z-w)^{-1}.
$$
By the universal property of $O(\gg,\kappa)$ given in Example \ExampleI, we
get the first desired homomorphism. The case for $\cS(\gg)$ is  analogous. $\Box$

Combining Lemmas \SquareZero, \VertexCartan, and \OgModule, we get

\theorem{$O(\gs\gg)\ra \cW(\gg)$, $(\xi,\eta)(z)\mapsto\Theta_\cW^\xi(z)+b^\eta(z)$, with ${\bf d}\mapsto D(0)$, defines an $O(\gs\gg)$-algebra structure on $\cW(\gg)$.}

We also have the vertex algebra analogues of the relations \Wdef.

\lemma{
$$\eqalign{
\Theta^\xi_\cW(z) c^{\xi'}(w)&\sim c^{ad^*(\xi)\xi'}(w)(z-w)^{-1},~~~\Theta^\xi_\cW(z)\gamma^{\xi'}(w)\sim\gamma^{ad^*(\xi)\xi'}(w)(z-w)^{-1}\cr
b^\xi(z) c^{\xi'}(w)&\sim\bra\xi',\xi\ket(z-w)^{-1},~~~b^\xi(z)\gamma^{\xi'}(w)\sim0\cr
D(0)c^{\xi'}&=-\half :c^{ad^*(\xi_i)\xi'}c^{\xi_i'}:+\gamma^{\xi'},  ~~~D(0)\gamma^{\xi'}=:\gamma^{ad^*(\xi_i)\xi'}c^{\xi_i'}:
}$$}
\thmlab\WRelations

\lemma{\FF\AK ~$(\cW^*,D(0))$ is acyclic.}
\proof\thmlab\ChiralAcyclic Put $h=:\beta^{\xi_i}\partial
c^{\xi_i'}:$. Then we find that $J(0)h=0,~~~K(0)
h=\omega_\cW$. It follows that, $D(0)h(z)=\omega_\cW(z)$, implying that $[D(0),h(1)]=\omega_\cW(1)$. Since
$\omega_\cW(1)$ is diagonalizable with nonnegative eigenvalues, it
follows that the cohomology of $(\cW^*,D(0))$ is the same as the
cohomology of the subcomplex  $(\cW^*[0,0],D(0))$. Recall that $\cW^*[0,0]$ is canonically isomorphic to the classical Weil algebra $W$. From
the formulas for the vertex operators $J,K$, we see that
$J(0),K(0)$ restricted to $\cW[0,0]$ coincide with their classical
counterparts $d_{CE},d_K$ respectively under the isomorphism.
Hence $d_W$ coincides with $D(0)$ restricted to $\cW[0,0]$. Thus $(\cW[0,0],D(0))$
and $(W,d_W)$ are isomorphic as complexes. But the latter is acyclic. $\Box$

\lemma{$J=:(\Theta_\cS^{\xi_i}+\half\Theta_\cE^{\xi_i})c^{\xi_i'}:$.}
\proof\thmlab\ThetaC By definition, we have
$J=:c^{\xi_i'}\Theta_\cS^{\xi_i}:+\half:c^{\xi_i'}\Theta_\cE^{\xi_i}:$.
Since the $\Theta_\cS^\xi$ commute with the $c^{\xi'}$, it
suffices to show that
$$
:c^{\xi_i'}\Theta_\cE^{\xi_i}:=:\Theta_\cE^{\xi_i}c^{\xi_i'}:
$$
By Lemma \Associator, we have
$$
:\Theta_\cE^{\xi_i}c^{\xi_i'}:=-:(c^{\xi_j'}b^{[\xi_i,\xi_j]})c^{\xi_i'}:
=-:c^{\xi_j'}b^{[\xi_i,\xi_j]}c^{\xi_i'}: - \partial c^{\xi_j'}\bra[\xi_i,\xi_j],\xi_i'\ket.
$$
The second term on the right vanishes because $ad^*(\xi_i)\xi_i'=0$, while the first term coincides with
$:c^{\xi_i'}\Theta_\cE^{\xi_i}:$. $\Box$

We now define the vertex algebra analogue of a $G^*$-algebra of type $C$.

\definition{A $\cW(\gg)$-algebra is a differential vertex algebra $(\cA,d_\cA)$ equipped with a differential vertex algebra homomorphism $\rho_\cA:\cW(\gg)\ra\cA$. We define a homomorphism of $\cW(\gg)$-algebras (and modules) in an obvious way.}

\lemma{$\cW_{hor}=\bra b\ket\otimes\cS(\gg)$, $\cW_{bas}=\cW_{hor}^{\gg_\geq}$.}
\proof
Here $\bra b\ket$ is the vertex algebra generated by the $b^{\xi}\in\cE(\gg)$,
and $(\cdots)^{\gg_\geq}$ is the subspace of $(\cdots)$ annihilated by
the $\Theta_\cW^\xi(n)$, $n\geq0$, $\xi\in\gg$. The first equality follows immediately from the fact that a vertex operator in $\cE(\gg)$ commutes with the $b^\xi$ iff it is in $\bra b\ket$. The second equality follows from the definition of the basic subalgebra and the fact that $a\in\cW_{hor}$ commutes with the $\Theta_\cW^\xi$
iff it is annihilated by the $\Theta_\cW^\xi(n)$. $\Box$

Here is another small surprise: unlike in the classical case where $d_K|W_{hor}=0$ and $d_{bas}=d_W|W_{bas}=0$, neither $K(0)|\cW_{hor}$ nor $D(0)|\cW_{bas}$ is zero in general.

Clearly $\omega_\cW\notin\cW_{hor}$. Since the vertex operators
$\Theta_\cW^\xi$ are primary of conformal weight 1 with respect to the
Virasoro element $\omega_\cW$, they do not commute with
$\omega_\cW$ unless $\gg$ is abelian, in which case the
$\Theta_\cW^\xi$ are identically zero. However, since the Fourier
mode $\omega_\cW(1)$ acts diagonalizably on $\cW$ and the vertex
operators $b^\xi,\Theta_\cW^\xi$ have weight 1, it follows that
$\omega_\cW(1)$ also acts diagonalizably on the basic vertex
subalgebra $\cW_{bas}$. Again, the subspace of zero eigenvalue is
canonically isomorphic to the classical basic subalgebra
$W_{bas}$. Since $[D(0),h(1)]=\omega_\cW(1)$ and $D(0)^2=0$, it
follows that $\omega_\cW(1)$ commutes with $D(0)$, hence its
action on $\cW_{bas}$ descends to $H_{bas}^*(\cW)$.  Here
$h=\beta^{\xi_i}\partial c^{\xi_i'}$. Note that
$h\notin\cW_{bas}$, so we cannot conclude that $\omega_\cW(1)$
acts by zero on cohomology. In fact, we will see that
$\omega_\cW(1)$ does not act by zero on cohomology.

\lemma{$H_{bas}^*(W(\gg))$ is canonically isomorphic to the
eigenspace of zero eigenvalue of $\omega_\cW(1)$ in
$\H_{bas}^*(\cW(\gg))$.}
\proof
Recall that $(W,d_W)$ is isomorphic to $(\cW[0,0],D(0))$. Restricted to
$\cW[0,0]$, the basic subalgebra condition reduces to
$b^\xi(0)a=\Theta^\xi_\cW(0)a=0$, $a\in\cW[0,0]$. It is easily
seen that this coincides with the basic subalgebra condition on
$W$ under the isomorphism above. This shows that $(\cW_{bas}[0],D(0))$ is isomorphic to $(W_{bas},d_W)$, hence $\H^*_{bas}(\cW(\gg))[0]\cong H^*_{bas}(W)$.  $\Box$

\subsec{Weil model for $\H^*_G(\cA)$}

\definition{For a given $O(\gs\gg)$-algebra $\cA$, we define its chiral $G$-equivariant cohomology to be
$$
\H^*_G(\cA)=\H^*_{bas}(\cA\otimes\cW(\gg)).
$$}

For $\cA=\C$, $\H^*_G(\C)=\H^*_{bas}(\cW(\gg))$, a vertex algebra which is already interesting and difficult to compute. Consider the map
$$
\cW(\gg)\ra\cA\otimes\cW(\gg),~~~a\mapsto 1\otimes a.
$$
In Example \ExampleVI, we saw that this is a basic homomorphism but not an $O(\gs\gg)$-algebra homomorphism in general. This induces a vertex algebra homomorphism
$$
\kappa_G:\H_G^*(\C)=\H_{bas}^*(\cW(\gg))\ra\H_{bas}^*(\cA\otimes\cW(\gg))=\H_G^*(\cA).
$$
This is our vertex algebra analogue of the Chern-Weil map.

Recall that given a manifold $M$, the vertex algebra $\cQ(M)$ has a Virasoro element (in local form)
$$
\omega_\cQ=:\beta^i\partial\gamma^i:-:b^i\partial c^i:
$$
which is $d_\cQ$-exact in $(\cQ(M),d_\cQ)$, and $\omega_\cQ(1)$ acts diagonalizably with eigenvalues in $\Z_\geq$. Since $\omega_\cW$ is $D(0)$-closed in $(\cW(\gg),D(0))$, it follows that
$$
\omega_{\cQ\otimes\cW}=\omega_\cQ\otimes 1+1\otimes\omega_\cW
$$
is also $d_{\cQ\otimes\cW}$-closed. In particular, $\omega_{\cQ\otimes\cW}(1)$ commutes with $d_{\cQ\otimes\cW}$. Again, since the vertex operators $\iota_\xi\otimes 1+1\otimes b^\xi$, $L_\xi\otimes 1+1\otimes\Theta^\xi_\cW$ have weight 1, it follows that  $\omega_{\cQ\otimes\cW}(1)$ acts diagonalizably on the basic subalgebra $(\cQ(M)\otimes\cW(\gg))_{bas}$ and on the basic cohomology $\H_G^*(\cQ(M))$.
Note, however, that the vertex operator $\omega_{\cQ\otimes\cW}$ is not a basic element in general.

More generally, consider an $O(\gs\gg)$-algebra $(\cA,d_\cA)$ with $O(\gs\gg)$-structure $(\xi,\eta)\mapsto L_\xi+\iota_\eta$. Suppose that $\cA$ has no negative weight elements. Then $\cA[0]$ is canonically a commutative associative algebra with product $\circ_{-1}$. Moreover, the operators $d_\cA,L_\xi\circ_0,\iota_\eta\circ_0$ on $\cA[0]$ define a $G^*$-structure on $\cA[0]$. If we assume, furthermore, that $\cA$ has a $d_\cA$-closed Virasoro element $\omega_\cA$ such that $L_\xi,\iota_\eta$ are primary of conformal weight 1 with respective to $\omega_\cA$, then
$\omega_{\cA\otimes\cW}=\omega_\cA\otimes1+1\otimes\omega_\cW$ is $(d_\cA+D(0))$-closed
Virasoro element in $\cA\otimes\cW(\gg)$.

\lemma{Let $\cA$ be a vertex algebra and $\omega\in\cA$ be a Virasoro element.
If $a\in\cA$ is primary of conformal weight one with respect to $\omega$, then
$\omega(m)$, $m\geq0$, preserves the subalgebra $Com(\bra a\ket,\cA)$.}
\proof
By Lemma \Commutator, the OPE $\omega(z)a(w)\sim a(w)(z-w)^{-1}+\partial a(w)(z-w)^{-2}$ translates into
$$
[\omega(m),a(n)]=-n a(m+n-1).
$$
Recall that $b\in Com(\bra a\ket,\cA)$ iff $a(n)b=0$ for all $n\geq0$.
Thus for such an element $b$, we have
$$
a(n)\omega(m)b=na(m+n-1)b=0
$$
for all $n,m\geq0$. Thus $\omega(m)b\in Com(\bra a\ket,\cA)$ for all $m\geq 0$. $\Box$

\theorem{Let $(\cA,d_\cA)$ be an $O(\gs\gg)$-algebra with no negative weight elements. Then the chiral equivariant cohomology $\H_G^*(\cA)$ is a degree-weight graded vertex algebra with
$\Z_\geq$-valued weights such that $\H^*_G(\cA)[0]=H^*_G(\cA[0])$.
If $\cA$ has a $d_\cA$-closed Virasoro structure $\omega_\cA$ such that the $O(\gs\gg)$-structure on $\cA$ are given by primary operators $L_\xi,\iota_\eta$ of conformal weight 1, then 
the operators $\omega_{\cA\otimes\cW}(m)$ induce an action on $\H_G^*(\cA)$ for all $m\geq0$.}
\proof\thmlab\ContainClassical 
Obviously $\cA\otimes\cW(\gg)$ has no negative weight elements, so that the same holds for its chiral basic cohomology. The weight zero subspace of $\cA\otimes\cW(\gg)$ is the tensor product of weight zero spaces $\cA[0]\otimes\cW[0,0]$. We saw that $(\cW[0,0],D(0))=(W,d_W)$ is the classical Weil algebra, and that $(\cA[0],d_\cA)$ is canonically a $G^*$-algebra.
It is clear that $(\cA\otimes\cW(\gg))_{bas}[0]$ coincides with the classical  basic complex $(\cA[0]\otimes W)_{bas}$. This yields the first assertion.

The basic subalgebra $(\cA\otimes\cW(\gg))_{bas}$ consists of elements annihilated by
the $n\geq0$ Fourier modes of the vertex operators $\iota_\xi\otimes 1+1\otimes b^\xi$ and the $L_\xi\otimes 1 +1\otimes\Theta_\cW^\xi$, each of which is primary of conformal weight one with respect to the Virasoro element $\omega_{\cA\otimes\cW}$. By the preceding lemma, $\omega_{\cA\otimes\cW}(m)$, for all $m\geq0$, acts on the basic subalgebra. Since $\omega_\cW$ is
$D(0)$-exact, and $\omega_\cA$ is assumed  $d_\cA$-closed, it follows that
$d_\cA+D(0)$ commutes with $\omega_{\cA\otimes\cW}$. Hence the action of its $m\geq0$ Fourier modes descends to an action on the basic cohomology $\H_G^*(\cA)$. $\Box$

\remark{From this, it is clear that our Chern-Weil map $\kappa_G$ restricts to the classical Chern-Weil map on the weight zero subspaces.}

%

\subsec{Cartan model for $\H_G^*(\cA)$}

We introduce the vertex algebra analogues of the Mathai-Quillen
isomorphism and the Cartan model. We will show that the Cartan
model is equivalent to the Weil model in the vertex algebra
setting.

Given a vector space $V$, we call a linear map $\phi\in End(V)$ pronilpotent if the restriction of $\phi$ to any finite dimensional subspace of $V$ is nilpotent. In this case,
$$
e^\phi=1+\phi+{1\over2!}\phi^2+{1\over3!}\phi^3+\cdots
$$
is a well-defined automorphism of $V$. Let $\cA$ be a vertex
(super) algebra and $a\in\cA$ a homogeneous vertex operator such
that the Fourier mode $\hat a(0)$ is pronilpotent. Since $\hat
a(0)$ is a derivation of all circle products, it follows that
$e^{\hat a(0)}:\cA\ra\cA$ is an automorphism of the vertex
algebra. As a reminder, we shall write $a(0)$ instead
of $\hat a(0)$ for notational simplicity.

Let $(\cA,d_\cA)$ be a $\cW(\gg)$-algebra and $(\cB,d_\cB)$ an $O(\gs\gg)$-algebra. Consider the vertex operator
$$
\phi(z)=c^{\xi_i'}(z)\otimes\iota_{\xi_i}(z)\in\cA\otimes\cB.
$$
We claim that the zeroth Fourier mode $\phi(0)$ is pronilpotent, as an operator on $\cA\otimes\cB$. It suffices to show that for any given homogeneous element $u\otimes v\in\cA\otimes\cB$, we have $\phi(0)^k(u\otimes v)=0$ for $k>>0$.
First, note that $c^{\xi'}(p)u=0=\iota_\xi(p)v$ for $p>>0$. In particular, there is an integer $N>0$ such that $\phi(0)(u\otimes v)=\sum_{|p|<N} (-1)^{|u|}c^{\xi_i'}(-p-1)u\otimes\iota_{\xi_i}(p)v$. Second, note that the modes $c^{\xi_i'}(p),\iota_{\xi_i}(q)$ are all pairwise anti-commuting, and in particular, each one is square-zero. This shows that $\phi(0)^k(u\otimes v)=0$ for $k>2N~dim~\gg$.

By analogy with the classical case, we call
$$\Phi=e^{\phi(0)}\equiv e^{\hat\phi(0)}$$ the Mathai-Quillen isomorphism of
$\cA\otimes\cB$.

\theorem{The Mathai-Quillen isomorphism satisfies
$$\eqalign{
\Phi (L_\xi\otimes 1 + 1\otimes L_\xi) &= L_\xi\otimes 1 +1\otimes
L_\xi\cr \Phi (\iota_\xi\otimes 1 + 1\otimes \iota_\xi) &=
\iota_\xi\otimes 1 \cr \Phi d \Phi^{-1} &=
d+(-\gamma^{\xi_i'}\otimes \iota_{\xi_i} + c^{\xi_i'}\otimes
L_{\xi_i})(0) }$$ where $d=d_\cA\otimes 1+1\otimes d_\cB$.} \proof
The vertex operators  $L_\xi\otimes 1+1\otimes L_\xi$ commute with
the vertex operator $\phi$ because their OPE with $\phi$ is
regular, since the $c^{\xi_i'}$ transform in the coadjoint module
of $\gg$ while the $\iota_{\xi_i}$ transform in the adjoint
module. It follows that
$$
\phi(0)(L_\xi\otimes 1+1\otimes L_\xi)=\phi\circ_0(L_\xi\otimes 1+1\otimes L_\xi)=0.
$$
This proves the first equality.

Next a simple OPE calculation yields
$$\eqalign{
\phi\circ_0(\iota_\xi\otimes 1) &= -1\otimes \iota_\xi\cr
\phi\circ_0(1\otimes \iota_\xi) &= 0\cr
}$$
This gives the second equality.

Finally, as operators on $\cA\otimes\cB$:
$$
[\phi(0),d]=-(d\phi)(0),~~~[\phi(0),[\phi(0),d]]=-(\phi\circ_0d\phi)(0).
$$
Applying Lemma \WRelations~and using that $(\cA,d_\cA)$ is a $\cW(\gg)$-algebra, we get
$$\eqalign{
d\phi&=d_\cA c^{\xi_j'}\otimes\iota_{\xi_j}-c^{\xi_j'}\otimes d_\cB\iota_{\xi_j}\cr
&=(-\half :c^{ad^*(\xi_i)\xi_j'}c^{\xi_i'}:+\gamma^{\xi'})\otimes\iota_{\xi_j}-c^{\xi_j'}\otimes L_{\xi_j}\cr
}$$
Since $\phi$ has regular OPE with the first term, we find that
$$
\phi\circ_0 d\phi=\phi\circ_0(c^{\xi_j'}\otimes L_{\xi_j})
=c^{\xi_i'}c^{\xi_j'}\otimes\iota_{[\xi_i,\xi_j]}.
$$
The last expression has regular OPE with $\phi$, hence applying $[\phi(0),-]$ to $d$
more than twice yields zero. It follows that
$$
\Phi d\Phi^{-1}=d-(d\phi)(0)-\half \phi\circ_0 d\phi
$$
yields the third equality. $\Box$

It follows that $\Phi(\cA\otimes\cB)_{hor}=\cA_{hor}\otimes\cB$.
Specializing to the case $\cA=\cW(\gg)$, we have $\cA_{hor}=\bra
b\ket\otimes\cS(\gg)$. Since $\Phi$ is an $O(\gg,0)$-module
homomorphism by the preceding theorem, it follows that
$$
\Phi(\cW\otimes\cB)_{bas}=(\bra b\ket\otimes\cS(\gg)\otimes\cB)^{\gg_\geq}=:C_G(\cB).
$$
Put
$$
d_G=\Phi d\Phi^{-1}|C_G(\cB).
$$
For any $O(\gs\gg)$-algebra $\cB$, the cohomology of the differential vertex algebra $(C_G(\cB),d_G)$
will be called the Cartan model for $\H_G^*(\cB)$. We have the following vertex algebra analogue of Cartan's fundamental theorem:

\theorem{The Mathai-Quillen isomorphism induces the differential vertex algebra homomorphism
$$
((\cW(\gg)\otimes\cB)_{bas},d_{\cW\otimes\cB})\ra(C_G(\cB),d_G).
$$
Hence
$$
\H_G^*(\cB)=H^*((\cW(\gg)\otimes\cB)_{bas},d_{\cW\otimes\cB})\cong H^*(C_G(\cB),d_G).
$$}

\newsec{Abelian Case}

We now specialize to the case when $G=T$ is an $n$ dimensional
torus. For the trivial $O(\gs\gg)$-algebra $\C$, we give a complete description of the chiral
$T$-equivariant cohomology $\H^*_G(\C)$. Recall that $\C=\cQ(pt)$, so that $\H^*_G(\C)$ is a vertex algebra analogue of $H^*_G(\C)$, the classical equivariant cohomology of a point. For general $O(\gs\gg)$-algebras $\cA$, we derive the analogue of a well-known spectral sequence that
computes the classical $T$-equivariant cohomology.

\subsec{The case $\cA=\C$}


Our first task is to compute
$$
\H^*_T(\C)=\H^*_{bas}(\cW(\gt)).
$$
Since $T$ is abelian, both the adjoint and the coadjoint modules
are trivial. In follows that all the vertex operators
$\Theta^\xi_\cW$ in \ThetaW~are identically zero. Likewise the $J$
\DJK~ is also zero. The differential $D(0)$ on $\cW$ is just
$$
K(0)=\sum_{n\in\Z}\gamma^{\xi_i'}(n-1)~b^{\xi_i}(-n)
$$
Let $\bra\gamma\ket$ be the abelian vertex algebra generated by the $\gamma^{\xi'}$, $\xi'\in\gt^*$, in $\cW(\gt)$.

\theorem{The inclusion $\bra\gamma\ket\subset\cW(\gt)$ induces a canonical isomorphism
$$
\H_T^*(\C)=\H_{bas}^*(\cW(\gt))\cong\bra\gamma\ket.
$$}
\proof\thmlab\AbelianCase Since the $\Theta^\xi_\cW$ are zero, it
follows that $\cW_{bas}=\cW_{hor}=\bra b\ket\otimes\cS(\gt)=\bra
b,\beta,\gamma\ket$, the vertex subalgebra of $\cW(\gg)$ generated
by the vertex operators $b^\xi,\beta^\xi,\gamma^{\xi'}$. Since the
$b^{\xi}(n)$, $n\geq0$, act by zero on this space, the basic
differential becomes
$$
d_{bas}=D(0)|\cW_{bas}=\sum_{n>0}\gamma^{\xi_i'}(n-1)~b^{\xi_i}(-n).
$$
We claim that the odd operator
$$
F=\sum_{n>0}\beta^{\xi_i}(-n)c^{\xi_i'}(n-1)
$$
is the a homotopy inverse for $d_{bas}$. First, it is obvious that $F\cW_{bas}\subset\cW_{bas}$. Moreover we have
$$
[d_{bas},F]=\sum_{n>0}\left(-b^{\xi_i}(-n)c^{\xi_i'}(n-1)+\beta^{\xi_i}(-n)\gamma^{\xi_i'}(n-1)\right)
$$
This acts diagonalizably on $\cW_{bas}$ with eigenvalues in $\Z_\geq$. This shows that
all basic cohomology is concentrated in eigenspace-zero. This is the subspace of $\cW_{bas}$
annihilated by the $c^{\xi'}(n-1),\gamma^{\xi_i'}(n-1)$, $n>0$, which is just the subalgebra $\bra\gamma\ket$. On the other hand $d_{bas}$ is identically zero on this subalgebra. $\Box$

Note that $\bra\gamma\ket$ is a free polynomial algebra generated
by the commuting vertex operators $\partial^k\gamma^{\xi'}$,
$k\geq0$, which are linear in $\xi'\in\gt^*$. Each of these vertex
operators has cohomology degree 2 and Virasoro weight $k$. As
expected from Theorem \ContainClassical, the chiral equivariant
cohomology $\bra\gamma\ket$ contains the polynomial subalgebra
generated by the weight zero vertex operators $\gamma^{\xi'}$,
which is a copy of the classical equivariant cohomology
$S(\gt^*)$.


\subsec{A spectral sequence for $\H_T^*(\cA)$}

Recall that $\cW_{hor}=\bra b,\beta,\gamma\ket$ has a monomial basis given by interated Wick
products of the $b^{\xi_i},\beta^{\xi_i},\gamma^{\xi_i'}$ and
their derivatives. In particular, there is a $\Z_\geq$ valued
grading on $\cW_{hor}$, which we shall call the $b\#$, which is
given by the eigenvalues of the diagonalizable operator on
$\cW_{hor}$:
$$
B=\sum_{n>0}b^{\xi_i}(-n)c^{\xi_i'}(n-1)
$$
The idea is that the vertex operators $b$ are non-classical
(because they have weight one), and we should first ``crop'' them
from the chiral Cartan complex. Likewise for the $\beta$.

Let $(\cA,d_\cA)$ be a given $O(\gs\gt)$-algebra. Recall that the Cartan model for the chiral $T$-equivariant cohomology of $\cA$ is the cohomology of the complex
$$
C_T(\cA)=(\cW_{hor}\otimes\cA)^{\gt_\geq}=\bra b,\beta,\gamma\ket\otimes\cA^{\gt_\geq}.
$$
The second equality follows from the important fact that the $O(\gt)$-structure on $\cW(\gt)$ is trivial because $T$ is abelian. The differential is
$$
d_T=D(0)\otimes 1+1\otimes d_\cQ-(\gamma^{\xi_i'}\otimes \iota_{\xi_i})(0)+(c^{\xi_i'}\otimes L_{\xi_i})(0).
$$
Again, $D(0)=K(0)$ because $T$ is abelian. Consider the vertex subalgebra
$$
\cC_T(\cA)=\bra\gamma\ket\otimes\cA^{\gt_\geq}.
$$
We claim that
$$
d_T|\cC_T(\cA)=1\otimes d_\cA-(\gamma^{\xi_i'}\otimes\iota_{\xi_i})(0).
$$
By Theorem \AbelianCase, $D(0)|\bra\gamma\ket=0$. We have
$(c^{\xi_i'}\otimes L_{\xi_i})(0)=\sum_{n\in\Z}c^{\xi_i'}(-n-1) L_{\xi_i}(n)$. Since $L_{\xi}(n)\cA^{\gt_\geq}=0$ for $n\geq0$, and $c^{\xi'}(-n-1)\bra\gamma\ket=0$ for $n<0$,
we see that $d_T$ reduces to the desired form. Thus we obtain an inclusion of complexes
$$
\varphi:(\cC_T(\cA),d_T)\hra (C_T(\cA),d_T).
$$
We shall call the first complex the {\it small} chiral Cartan complex of $\cA$.

\lemma{The map induced by $\varphi$ on cohomology is surjective.}
\proof Let $a\in\bra b,\beta,\gamma\ket\otimes\cA^{\gt_\geq}$ be
nonzero and $d_T$-closed. We will show that $a$ is
$d_T$-cohomologous to an element in
$\bra\gamma\ket\otimes\cA^{\gt_\geq}$. Let $a_{max}$ be the
component of $a$ with the maximum $b\#$. Suppose this $b\#$ is
positive. We can write $a_{max}=\sum p_j\otimes\omega_j$ where the
$\omega_j\in\cA^{\gt_\geq}$ are linearly independent elements and
the $p_j\in\bra b,\beta,\gamma\ket$ have the same maximum $b\#$.

We look at the effects of each of the four terms in $d_T$ on the
$b\#$ in $\bra b,\beta,\gamma\ket\otimes\cA^{\gt_\geq}$. We have
\eqn\dumb{\matrix{ operators & b\# &\cr K(0)\otimes
1=\sum_{n\geq0}\gamma^{\xi_i'}(n)b^{\xi_i}(-n-1)& +1\cr 1\otimes
d_\cA & 0\cr -(\gamma^{\xi_i'}\otimes\iota_{\xi_i})(0)& 0\cr
(c^{\xi_i'}\otimes L_{\xi_i})(0)=\sum_{n<0}c^{\xi_i'}(-n-1)
L_{\xi_i}(n) &  -1. }} It follows that the term $(K(0)\otimes
1)a_{max}$ is the component in $d_T a$ with the highest $b\#$.
Since $d_Ta=0$, this highest term must be zero, hence $\sum
K(0)p_j\otimes\omega_j=0$, implying that $K(0)p_j=0$. By Theorem
\AbelianCase, $p_j=K(0)q_j$ for some $q_j$ with $b\#$ one less
than that of $p_j$. Put
$$
a'=a-d_T\sum q_j\otimes\omega_j=a-a_{max}+\cdots
$$
where the terms in $\cdots$ have only components with lower $b\#$
than $a_{max}$. Thus by induction, we see that $a$ is
$d_T$-cohomologous to an element with $b\#$ zero i.e. in
$\bra\beta,\gamma\ket\otimes\cA^{\gt_\geq}$. So we may
assume that $a$ does not depend on the $b$.

Next we show that $a$ does not depend on the $\beta$ either. Since $a$ is independent of the $b$, it follows that the $b^{\xi_i}(-n-1)a$, $n\geq0$, $i=1,..,dim~\gt$, are linearly independent. Again by \dumb, $(K(0)\otimes 1)a=0$. It follows that $\gamma^{\xi_i'}(n)b^{\xi_i}(-n-1)a=0$ for all $n\geq0$.
This shows that $a$ is independent of the $\beta$. This completes the proof. $\Box$

\lemma{The map induced by $\varphi$ on cohomology is injective.}
\proof Given an element $\omega$ in the small complex
$\bra\gamma\ket\otimes\cA^{\gt_\geq}$ which is $d_T$-exact in the big
complex $\bra b,\beta,\gamma\ket\otimes\cA^{\gt_\geq}$, we want to
show that $\omega$ is $d_T$-exact in the small complex. Write
$\omega=d_Ta$. We want to find $a'$ in the small complex so that
$d_Ta=d_Ta'$. Again, let $a_{max}$ be the component of $a$ with
the maximum $b\#$. By \dumb, $(K(0)\otimes 1)a_{max}$ is the
component of $d_Ta$ with the maximum and positive $b\#$. Since
$d_Ta$ is in the small complex, it does not depend on the $b$,
implying that $(K(0)\otimes 1)a_{max}=0$. Thus we get the shape
$a_{max}=\sum K(0) q_j\otimes\omega_j$ and that
$$
a-d_T\sum q_j\otimes\omega_j=a-a_{max}+\cdots
$$
where $\cdots$ have only components with lower $b\#$ than
$a_{max}$, as before. Thus we may as well assume that $a$ does not
depend on the $b$. Following the same argument as in the preceding
lemma, we see that $a$ does not depend on the $\beta$ either.
$\Box$

\theorem{For any $O(\gs\gt)$-algebra $(\cA,d_\cA)$, we have $\H_T^*(\cA)\cong H^*(\cC_T(\cA),d_T)$.}
\proof
The two preceding lemmas show that $\varphi$ induces an isomorphism on cohomology. $\Box$

From now on, we specialize to the case $\cA=\cQ(M)$ where $M$ is a $T$-manifold.

\theorem{For any $T$-manifold $M$, we have $\H_T^*(\cQ(M))\cong H^*(\cC_T(\cQ(M)),d_T)$.}

Observe that the small complex on the right side is a double complex with the differentials
$$\eqalign{
d&=1\otimes d_\cQ=\sum_{n\in\Z}\beta^i(n)c^i(-n-1)\cr
\delta&=-(\gamma^{\xi_i'}\otimes\iota_{\xi_i})(0)=-\sum_{n\geq0}\gamma^{\xi_i'}(-n-1)\iota_{\xi_i}(n).
}$$
where the expression for $d$ makes use of a choice of local
coordinates on $M$. This is an analogue of the double complex
structure in the classical Cartan model (see Chap. 6 \GS). Note
however that in the classical case, the double complex structure
is on the Cartan complex itself, whereas in our case, it is on a
much smaller subcomplex of the Cartan complex. It is also clear
that the weight zero piece of the small complex coincides with the
classical Cartan complex.

As usual, associated to the double complex structure on the small Cartan complex, there are two filtrations and two spectral sequences \HS. We shall consider the following one:
$$
F^n_k=\oplus_{p+q=n,~p\geq k}~ \cC_T^{p,q},~~~~\cC_T^{p,q}=\bra\gamma\ket^p\otimes
(\cQ(M)^{q-p})^{\gt_\geq}.
$$
Here $p$ denotes the $\gamma$-number on $\bra\gamma\ket$, and $q-p$ is the $bc$-number on $\cQ(M)$. Let $(E_r,\delta_r)$ be the spectral sequence associated with this filtration.

\theorem{In each weight, the spectral sequence $(E_r,\delta_r)$ converges to the graded object associated with $H^*(\cC_T(\cQ(M)),d+\delta)$. In fact, in each weight, the spectral sequence collapses at $E_r$ for some $r$.}
\proof The first statement follows
immediately from the fact that both $d,\delta$ are operators of
weight zero, and the filtering spaces $F_k$ and the terms in the
spectral sequence are all graded by the weight. In a given weight
$m$, we have $\cQ(M)^{q-p}[m]=0$ for $|q-p|>>0$. For if $|q-p|$ is
not bounded then either the vertex operators $\partial^k c^i$ or
the $\partial^k b^i$ would have to be present on some coordinate
open set of $M$, with unbounded $k$, because these operators are
fermionic. But $wt~\partial^kc^i=k$ and $wt~\partial^kb^i=k+1$,
violating that $m$ is fixed and that the weights in $\cC_T$ are
bounded below by zero. This shows that $E_r[m]=E_{r+1}[m]=\cdots$
for all $r>>0$ (cf. p66 \GS). $\Box$

Note that
$$
E_1^{p,q}=H^q(\cC_T^{p,*},d),~~~~\delta_1=\delta:E_1^{p,q}\ra E_1^{p+1,q}.
$$

\newsec{Non-Abelian Case}

We shall now proceed to construct cohomology classes in the chiral
equivariant cohomology of the trivial $O(\gs\gg)$-algebra $\C$.
As before, $\gg$ is the complexified Lie algebra of the compact group $G$. We also
choose an orthonormal basis $\xi_i$ with respect to a fixed $G$-invariant pairing on $\gg$.
We shall often identify $\gg,\gg^*$ via this pairing, for convenience

\theorem{For any (connected) compact group $G$, the vertex operator $\gamma^{\xi_i}\partial\gamma^{\xi_i}$ represents a nonzero class in $\H_G^4(\C)[1]$.}
\thmlab\gammagamma
\proof It is straightforward to show that this vertex operator is basic and closed. Since the basic complex $\cW_{bas}$ contains no vertex operators involving the $c^{\xi_i}$, for any given homogeneous element $\mu\in\cW_{bas}$, $K(0)\mu=(\gamma^{\xi_i}b^{\xi_i})(0)\mu$ is either zero or it will contain some $b^{\xi_i}$. Thus if
$$
(K(0)+J(0))\mu=\gamma^{\xi_i}\partial\gamma^{\xi_i},
$$
then $J(0)\mu$ must have components (in the standard basis) of the form $\gamma^{\xi_i}\partial\gamma^{\xi_i}$. We will show that this leads to a contradiction.

First note that for $B=\beta^{\xi_i}\beta^{\xi_i}$, we have
$$
J(z)B(w)\sim (c^{\xi_i}\beta^{[\xi_i,\xi_j]}\beta^{\xi_j})(w)(z-w)^{-1}=0
$$
because $\beta^{[\xi,\xi_i]}\beta^{\xi_i}=0$. In particular $B(2)$
commutes with $J(0)$. Hence
$$
J(0)B(2)\mu=B(2)J(0)\mu
=B(2)(\gamma^{\xi_i}\partial\gamma^{\xi_i}-K(0)\mu)
=dim~\gg-B(2)K(0)\mu.
$$
But since $K(0)\mu$ is either zero or contains some $b^\xi$, the same is true of the second term in the last expression above. But since that second term has weight zero and $b^\xi$ has weight one, this forces that second term to be zero. Since $J(0)$ has $bc\#$ 1 and weight zero, and $dim~\gg$ has $bc\#$ 0 and weight zero, there must be a component of $B(2)\mu$ having $bc\#$ -1 and weight zero. This is absurd because weight zero elements cannot have negative $bc\#$. $\Box$

\remark{The vertex operator
$B=\beta^{\xi_i}\beta^{\xi_i}\in\cW(\gg)_{bas}$ turns out to be part
of a current algebra $O(sl_2,-\frac{dim(\gg)}{8}
\kappa)$-structure which plays a fundamental role in the
description of the full chiral equivariant cohomology of $\C$. This will
be explained in a future follow-up paper.} 

\subsec{Weight one classes}

Theorem \ContainClassical~gives a complete description of $\H^*_G(\C)[0]$, i.e. it
coincides with the classical equivariant cohomology $H^*_G(\C)$. We now give a complete description
of the weight one piece.

{\it Notations.} We identify $Sym(\gg)$ with the algebra $\C[\gamma^{\xi_1},...,\gamma^{\xi_n}]$, $n=dim~\gg$. For $P\in Hom_\gg(\gg,Sym(\gg))$, we write $P:\xi\mapsto P_\xi$. Throughout this subsection, $P$ shall denote such a map.

\lemma{Let $Q:\gg\ra Sym(\gg)$, $\xi\mapsto Q_\xi$, be any linear map such that $Q_{\xi_i}b^{\xi_i}\in\cW[1]$ is $\gg$-invariant. Then $Q$ is a $\gg$-module map.
Likewise the same is true under the assumption that $Q_{\xi_i}\beta^{\xi_i}$ or $Q_{\xi_i}\partial\gamma^{\xi_i}$ is $\gg$-invariant.}
\thmlab\QgInvariant
\proof
We will prove one case, the other two being similar. We have
$$
0=\Theta_\cW^\xi(0)Q_{\xi_i}b^{\xi_i}=(\Theta_\cS^\xi(0)Q_{\xi_i})b^{\xi_i}+
Q_{\xi_i}b^{[\xi_i,\xi]}.
$$
But the second term on the right is equal to $-Q_{[\xi,\xi_j]}b^{\xi_j}$. By linear independence of the $b^{\xi_i}$, it follows that $\Theta_\cS^\xi(0)Q_{\xi_i}=Q_{[\xi,\xi_j]}$. This says that $Q$ is a $\gg$-module map. $\Box$

\lemma{For $\omega\in\cS(\gg)$, $\omega$ is $\gg_\geq$ invariant iff $J(0)\omega=0$.}
\proof
Since $\omega$ has no  $b,c$,
$$
J(0)\omega=(c^{\xi_i}\Theta^{\xi_i}_\cS)(0)\omega=\sum_{n\geq0} {1\over n!} \partial^nc^{\xi_i}\Theta^{\xi_i}_\cS(n)\omega.
$$
Since all nonzero terms on the right side are independent,
$J(0)\omega=0$ iff $\Theta^{\xi_i}_\cS(n)\omega=0$ for $n\geq0$. $\Box$

\lemma{$P_{\xi_i}b^{\xi_i}\in\cW_{bas}$.}
\thmlab\PbBasic
\proof
The vertex operator $P_{\xi_i}b^{\xi_i}$ is clearly $\gg$-invariant, i.e. killed by the $\Theta_\cW^\xi(0)$.
Since it has weight one, it suffices to show that it is killed by $\Theta_\cW^\xi(1)$. Now
$\Theta_\cW^\xi(1)P_{\xi_i}b^{\xi_i}$ is the term with second order pole in the OPE
of $\Theta_\cW^\xi(z)(P_{\xi_i}b^{\xi_i})(w)$. Since $\Theta_\cW^\xi(z)$ has the shape $\beta\gamma+bc$, there is no double contraction by Wick's theorem. So there is no second order pole in the OPE. $\Box$

\lemma{$P_{\xi_i}\partial\gamma^{\xi_i}\in\cW_{bas}$, hence it lies in $Ker~J(0)$.}
\thmlab\PgammaJ
\proof
Again, it is clear that $P_{\xi_i}\partial\gamma^{\xi_i}$ is horizontal and $\gg$-invariant. We have
$$
\Theta_\cW^\xi(1)P_{\xi_i}\partial\gamma^{\xi_i}=
\Theta_\cS^\xi(1)P_{\xi_i}\partial\gamma^{\xi_i}=\gamma^{\xi_i}P_{[\xi_i,\xi]}
=\gamma^{\xi_i}(-\gamma^{\xi_j}\beta^{[\xi_i,\xi_j]})(0)P_\xi=0.
$$
Since $P_{\xi_i}\partial\gamma^{\xi_i}$ has weight one, this shows that it is $\gg_\geq$-invariant, hence basic. $\Box$

\lemma{$P_{\xi_i}\beta^{\xi_i}\in\cW_{bas}$, hence it lies in $Ker~J(0)$.}
\thmlab\PbetaBasicJ
\proof
Recall that $B=\beta^{\xi_j}\beta^{\xi_j}$ commutes with $J$.
Thus by the preceding lemma, we have
$$
0=B(1)J(0)P_{\xi_i}\partial\gamma^{\xi_i}=
J(0)B(1)P_{\xi_i}\partial\gamma^{\xi_i}=J(0)(B(1)P_{\xi_i})\partial\gamma^{\xi_i}
+2J(0)P_{\xi_i}\beta^{\xi_i}.
$$
For the last equality, we have used the second identity in Lemma \Associator, and the fact that
$B(1)\partial\gamma^{\xi_i}=2\beta^{\xi_i}$ and $(B(0)P_{\xi_i})(0)\partial\gamma^{\xi_i}=0$, which follow from Wick's theorem.
Since $B$ is $\gg$-invariant, it follows that the $B(1)P_{\xi_i}$ transform in the adjoint module, i.e.
$\xi\mapsto B(1)P_{\xi}$ defines element of $Hom_\gg(\gg,Sym(\gg))$. By the preceding lemma,
$J(0)(B(1)P_{\xi_i})\partial\gamma^{\xi_i}=0$. This shows that $J(0)P_{\xi_i}\beta^{\xi_i}=0$. $\Box$

\lemma{$J(0)P_{\xi_i}b^{\xi_i}=0$.}
\thmlab\PJ
\proof
Note that $K(0)P_{\xi_i}\beta^{\xi_i}=P_{\xi_i}b^{\xi_i}$.
Since $J(0),K(0)$ commute, it follows from the preceding lemma that $J(0)P_{\xi_i}b^{\xi_i}=0$. $\Box$

\theorem{The weight one subspace $\H_G^*(\C)$[1] is canonically
isomorphic to $Hom_G(\gg,Sym(\gg))$.}
\thmlab\WeightOneCohomology
\proof
Any element $a\in \cW_{hor}[1]$ can be uniquely written as
\eqn\dumb{
a=P_{\xi_i}b^{\xi_i}+Q_{\xi_i}\beta^{\xi_i}+R_{\xi_i}\partial\gamma^{\xi_i}
}
where $P,Q,R$ are linear maps $\gg\ra Sym(\gg)$. Suppose $a$ is $\gg$-invariant.
Since the $\beta^\xi,\gamma^{\xi},b^{\xi}$ form three copies of the adjoint module,
each of the three terms in $a$ above must be separately $\gg$-invariant. By Lemma \QgInvariant, it follows that $P,Q,R$ are $\gg$-module maps. By the preceding four lemmas, the three terms in \dumb~are separately basic and killed by $J(0)$.

We have
$$
D(0)P_{\xi_i}\beta^{\xi_i}=P_{\xi_i}b^{\xi_i}.
$$
It follows that the first term in \dumb~ is $D(0)$-exact.
Suppose, in addition, that $a$ is $D(0)$-closed, i.e. $a$ represents a class in $\H_G^*(\C)[1]$.
Then $\dumb$ represents the same class if we drop the first term, so we may assume that $P=0$. Then
$$
0=D(0)a=Q_{\xi_i}b^{\xi_i}+D(0)R_{\xi_i}\partial\gamma^{\xi_i}.
$$
The second term is zero because $R_{\xi_i}\partial\gamma^{\xi_i}$ is obviously killed by $K(0)$ and is killed by $J(0)$ by Lemma \PgammaJ. It follows that $Q=0$. This shows that we have a canonical surjective linear map
\eqn\dumb{
Hom_\gg(\gg,Sym(\gg))\ra \H_G^*(\C)[1],~~~R\mapsto [R_{\xi_i}\partial\gamma^{\xi_i}].
}
Suppose a given $R$ is killed by this map. Then for some $P',Q',R'\in Hom_\gg(\gg,Sym(\gg))$, we have
$$
R_{\xi_i}\partial\gamma^{\xi_i}=D(0)(P_{\xi_i}'b^{\xi_i'}+Q_{\xi_i}'\beta^{\xi_i}+R_{\xi_i}'\partial\gamma^{\xi_i}) =Q_{\xi_i}'b^{\xi}.
$$
This implies that $R=Q'=0$. This shows that \dumb~is injective. $\Box$

\subsec{Weight two classes}

For simplicity, we shall assume that $G$ is simple throughout this subsection.
The Virasoro algebra will be playing a crucial role here.

\lemma{Let $A$ be a vertex algebra in which $L_1,L_2\in A$ are
Virasoro elements of central charges $c_1,c_2$. Suppose that $L_2$
is quasi-primary of conformal weight 2 with respect to $L_1$, i.e.
$$
L_1(z)L_2(w)\sim \half c_3 (z-w)^{-4}+2L_2(w)(z-w)^{-2}+\partial L_2(w)(z-w)^{-1}
$$
for some scalar $c_3$. Then $L_1-L_2$ is a Virasoro element of central charge $c_1+c_2-2c_3$.}
\proof
This is a well-known trick borrowed from physics \GKO\Pol. By Lemma \Commutator~again, we find
$$
L_2(z)L_1(w)\sim \half c_3 (z-w)^{-4}+2L_2(w)(z-w)^{-2}+\partial L_2(w)(z-w)^{-1}.
$$
Now combining the four OPEs of $L_i(z)L_j(w)$, $i,j=1,2$, we get
the desired OPE of $L_1-L_2$ with itself. $\Box$

Since $\cS(\gg)$ is an $O(\gg,-\kappa)$-module, by Lemma \OgModule, it has a Virasoro element given by the Sugawara-Sommerfield formula
$$
L_\cS=-:\Theta_\cS^{\xi_i}\Theta_\cS^{\xi_i}:
$$
where $\xi_i$ is an orthonormal basis of $(\gg,\kappa)$, as in Example \ExampleIII.

\lemma{$L_\cS$ is a quasi-primary of conformal weight 2 with respect to $\omega_\cS$. In fact,
$$
\omega_\cS(z)L_\cS(w)\sim dim~\gg~(z-w)^{-4}+2L_\cS(w)~(z-w)^{-2}+\partial L_\cS(w)(z-w)^{-1}.
$$}
\proof By Lemma \Associator, we find that if $a,b$ are primary of
conformal weight 1 with respect to a Virasoro element $\omega$, then we have
$$
\omega(z) c(w)\sim\sum_{k\geq0}(a\circ_k b)(w)(z-w)^{-k-3}+2 c(w)(z-w)^{-2}+\partial c(w)(z-w)^{-1}.
$$
where $c=:ab:$. By Lemma \ThetaPrimary, the $\Theta_\cS^\xi$ are primary of conformal weight 1 with respect to $\omega_\cS$, and so we can apply this to the case $\omega=\omega_\cS$, $a=-b=\Theta_\cS^{\xi_i}$, in $\cS(\gg)$, so that $c=L_\cS$ when we sum over $i$. We find that $a\circ_kb=\delta_{k,1}~\kappa(\xi_i,\xi_i)$. Summing over $i$, this becomes $\delta_{k,1}dim~\gg$, which yields the desired OPE of $\omega_\cS(z)L_\cS(w)$. $\Box$

\corollary{$\omega_\cS-L_\cS$ is a Virasoro element in $\cW(\gg)$ of central charge 0.}
\proof
By Lemma \FMSVirasoro, $\omega_\cS$ has central charge $2dim~\gg$. By Example \ExampleIII,
$L_\cS$ has central charge $2dim~\gg$ also. Now our assertion follows from the two preceding lemmas. $\Box$

\lemma{$\omega_\cS-L_\cS$ commutes with the $\Theta_\cS^\xi$, and
lives in the basic subalgebra $\cW(\gg)_{bas}$.}
\proof\thmlab\BasicProperty By Lemma \ThetaPrimary, the
$\Theta_\cS^\xi$ are primary of conformal weight 1 with respect to
$\omega_\cS$. By Example \ExampleIII, the same is true with
respect to $L_\cS$. It follows that
$(\omega_\cS(z)-L_\cS(z))\Theta_\cS^\xi(w)\sim0$. Now Lemma
\Commutator~implies our first assertion. Since
$\omega_\cS-L_\cS\in\cS(\gg)$, it also commutes with the
$b^\xi\in\cE(\gg)$. $\Box$

Note however that $\omega_\cS-L_\cS$ is not $D(0)$-closed. The
idea is to try to ``correct'' it, so that it becomes $D(0)$-closed
without destroying its {\it basic} property. Since
$\omega_\cS-L_\cS$ is basic, any correction would need to be basic
as well. In particular, it must be $G$-invariant and horizontal.
We examine the simplest nontrivial case first: $G=SU(2)$. Let's
ask: what is the simplest possible (say, lowest degree as a
polynomial in $\bra b\ket\otimes\cS(\gg)$) horizontal
$G$-invariant vertex operator $C$ such that $\omega_\cS-L_\cS+C$ is
$D(0)$-closed? Choose the standard basis $x,h,y$ for the
complexified Lie algebra $\gg=sl_2$.

\lemma{$C=-\gamma^{h'}b^xb^y+\half\gamma^{x'}b^xb^h-\half\gamma^{h'}b^yb^h$
is the unique homogeneous lowest degree element which is
$G$-invariant and horizontal, and makes $\L=\omega_\cS-L_\cS+C$,
$D(0)$-closed. Moreover, this $\L$ is also a Virasoro element with
central charge 0.} \proof\thmlab\CorrectionTerm  Note that all the
vertex operators appearing in $C$ commute with each other. It is
clear that $C$ is horizontal since it does not depend on the
$c^{\xi}$. Note also that $C$ comes from a cubic trace polynomial
in $\gg^*\otimes\gg\otimes\gg$, hence is $G$-invariant \Procesi.
We have verified the uniqueness assertion and $D(0)$-closed
condition by direct computations.

Clearly $C(z)C(w)\sim0$ and $C$ is primary of conformal weight 0 with
respect to $\omega_\cS$. By Lemmas \AdjointCoadjoint, \Associator,
we find that
$$
-L_\cS(z)\gamma^{\xi'}(w)=(:\Theta_\cS^{\xi_i}\Theta_\cS^{\xi_i}:)(z)\gamma^{\xi'}(w)\sim \gamma^{ad^*(\xi_i)ad^*(\xi_i)\xi'}(w)(z-w)^{-2}+\cdots=\gamma^{\xi'}(w)(z-w)^{-2}+\cdots
$$
where $\cdots$ are lower order poles, and the $\xi_i$ form an
orthonormal basis of $\gg$. Hence $-L_\cS(z)C(w)\sim
C(w)(z-w)^{-2}+\cdots$. This shows that
$(\omega_\cS-L_\cS)(z)C(w)\sim C(w)(z-w)^{-2}+\cdots$. That $\L$
is Virasoro element with central charge 0 now follows from the
next lemma. $\Box$

\lemma{Let $A$ be a vertex algebra and $L\in A$ be a Virasoro
element. Suppose that $a,b\in A$ are vertex operators with the
property that
$$
a(z)a(w)\sim~0,~~~L(z)a(w)\sim a(w)(z-w)^{-2}+b(w)(z-w)^{-1}.
$$
Then $L+a$ is Virasoro element with the same central charge.}
\proof By \VOcom, we find that $a(z)L(w)\sim
a(w)(z-w)^{-2}+(\partial a(w)+b(w))(z-w)^{-1}$. It follows that
$$\eqalign{
&(L(z)+a(z))(L(w)+a(w))\sim \half c(z-w)^{-4}+2L(w)(z-w)^{-2}+\partial L(w)(z-w)^{-1}\cr
&~~~+a(w)(z-w)^{-2}+(\partial a(w)-b(w))(z-w)^{-1}+ a(w)(z-w)^{-2}+b(w)(z-w)^{-1}
}$$
This shows that $L+a$ has the OPE of a Virasoro element of the same central charge $c$. $\Box$

Now comes a crucial observation. A computation using Wick's
theorem shows that
$$
C=K(0)(\Theta_\cS^{\xi_i}b^{\xi_i})=(K(0)\Theta_\cS^{\xi_i})b^{\xi_i}
$$
where the $\xi_i$ is an orthonormal basis of $(\gg,\kappa)$. This
computation suggests that for a general $G$, we should consider
the vertex operator
$$
\L{\br
def\over=}\omega_\cS-L_\cS+(K(0)\Theta_\cS^{\xi_i})b^{\xi_i}.
$$
From the shape of the vertex operator $K$ and the $\Theta_\cS^{\xi}$,
one finds that the last term of the right side above is a priori a $G$-invariant
vertex operator having the shape $bb\gamma$.

\theorem{$\L$ is a Virasoro element with central charge 0, which is basic, $D(0)$-closed, and satisfies
$$
\L(0)a=\partial a,~~~~\L(1)a=(wt~a)a
$$
for any homogeneous element
$a\in\bra\gamma\ket\cap\cW(\gg)_{bas}$, where $\bra\gamma\ket$ is
the vertex algebra generated by the $\gamma^{\xi'}$.}
\proof\thmlab\CohomologyVirasoro Applying Wick's theorem, we find
that the last term of $\L$ is
$$
C=(K(0)\Theta_\cS^{\xi_i})b^{\xi_i}=b^{\xi_i}b^{\xi_j}\gamma^{ad^*(\xi_i)\xi_j'}.
$$
This implies that $C(z)C(w)\sim0$. To apply the preceding lemma to show the Virasoro property for $\L$,
it remains to verify the OPE
$$
(\omega_\cS-L_\cS)(z)~\gamma^{\xi'}(w)\sim\gamma^{\xi'}(w)(z-w)^{-2}+\cdots
$$
where $\cdots$ means lower order poles. This is verbatim as in Lemma \CorrectionTerm.

To show the basic property, by Lemma \BasicProperty, it suffices to check it for $C$.
Clearly $C$ commutes with the $b^\xi$. Applying Wick's theorem again, we get easily
$(\Theta^\xi_\cS+\Theta^\xi_\cE)(z)~C(w)\sim0$. This shows that $C$ is basic.

Next, we claim that
$$
\L=D(0)(\Theta_\cS^{\xi_i}b^{\xi_i}+\beta^{\xi_i}\partial c^{\xi_i'}),
$$
hence it is automatically $D(0)$-closed. By Wick's theorem, we get
$$\eqalign{
K(0)(\beta^{\xi_j}\partial c^{\xi_j'})&=-:b^{\xi_i}\partial c^{\xi_i'}:+:\beta^{\xi_i}\partial\gamma^{\xi_i'}:=\omega_\cW\cr
J(0)(\beta^{\xi_i}\partial c^{\xi_i'})&=0\cr
(\Theta_\cS^{\xi_i}c^{\xi_i'})(0)(\Theta_\cS^{\xi_j}b^{\xi_j})&=:(\Theta_\cS^{\xi_i}c^{\xi_i'}\circ_0\Theta_\cS^{\xi_j})b^{\xi_j}:+
:\Theta_\cS^{\xi_j}(\Theta_\cS^{\xi_i}c^{\xi_i'}\circ_0b^{\xi_j}):\cr
&=:(\Theta_\cS^{[\xi_i,\xi_j]}c^{\xi_i'})b^{\xi_j}:-:\partial c^{\xi_j}b^{\xi_j}:+:\Theta_\cS^{\xi_i}\Theta_\cS^{\xi_i}:,
~~~~Lemma~\Associator\cr
&=:\Theta_\cS^{[\xi_i,\xi_j]}c^{\xi_i'}b^{\xi_j}:+:\Theta_\cS^{\xi_i}\Theta_\cS^{\xi_i}:-\omega_\cE\cr
&=-:\Theta_\cS^{\xi_i}\Theta_\cE^{\xi_i}:+:\Theta_\cS^{\xi_i}\Theta_\cS^{\xi_i}:-\omega_\cE\cr
\half(:\Theta_\cE^{\xi_i}c^{\xi_i'}:)(0)(\Theta_\cS^{\xi_j}b^{\xi_j})
&=:\half:\Theta_\cS^{\xi_j}(\Theta_\cE^{\xi_i}c^{\xi_i'}\circ_0b^{\xi_j}):
=:\Theta_\cS^{\xi_i}\Theta_\cE^{\xi_i}:
}$$
Applying Lemma \ThetaC, we find that the sum of the four left sides plus $K(0)(\Theta_\cS^{\xi_i}b^{\xi_i})=C$ yields $D(0)(\Theta_\cS^{\xi_i}b^{\xi_i}+\beta^{\xi_i}\partial c^{\xi_i'})$ on the one hand, and $\L$ on the other hand.

Finally, let $a\in\bra\gamma\ket\cap\cW(\gg)_{bas}$. Then in terms of the generators of $\cW$, $a$
does not depend on the vertex operators $b^\xi,c^{\xi'},\beta^\xi$. In particular $a$ commutes with $C$ above and with the $\Theta_\cE^\xi$.  Since $a$ is assumed basic, it commutes with
the $\Theta_\cW^\xi$, hence with the $\Theta_\cS^\xi$ as well. In particular, $a$ commutes with $L_\cS$.
It follows that
$$
\L(z)a(w)\sim\omega_\cS(z)a(w)\sim \cdots+(wt~a)a(w)(z-w)^{-2}+\partial a(w)(z-w)^{-1}
$$
where $\cdots$ here means terms with higher order poles. Here we have used Lemma \FMSVirasoro~to get the left side. This completes the proof. $\Box$

\corollary{$\L$ represents a nontrivial weight two class in $\H_G^0(\C)$.}
\proof
It is easy to verify that $\L$ has cohomology degree zero and weight two.
Since $a=\gamma^{\xi_i'}\partial\gamma^{\xi_i'}$ is a nonzero weight one class by Theorem \gammagamma, it follows that $\L(1)a=\L\circ_1a=a$ by the preceding theorem. Since circle products are cohomological operations, it follows that $\L$ cannot represent the zero class. $\Box$

\corollary{$\H_G^*(\C)$ is a non-abelian vertex algebra.}
\proof
Since $\L$ represents a nonzero Virasoro element in the vertex algebra $\H_G^*(\C)$, it does not commute with itself. $\Box$

\remark{This indicates that the departure of the chiral theory from the classical theory in the non-abelian case is quite dramatic.}

\lemma{Let $A$ be any vertex algebra and $a\in A$ such that
$a(m)=0$ for some $m<0$. Then $a$ commutes with $A$.} \proof By
Lemma \Commutator~we have 
$$
[a(m),b(q)]=\sum_p\left(\matrix{m\cr p}\right)(a\circ_p b)(m+q-p).
$$
Consider the maximum $m<0$ such that $a(m)=0$. Suppose $a$ does not commute with $b$, so that there exists a largest $N\geq0$ such that $a\circ_N b\neq0$. Pick $q=N-m-1$. Then
$$\eqalign{
0=[a(m),b(q)]\One&=\sum_{p=0}^N\left(\matrix{m\cr
p}\right)(a\circ_p b)(m+q-p)\One\cr &=\left(\matrix{m\cr
N}\right)(a\circ_N b)(-1)\One =\left(\matrix{m\cr
N}\right)a\circ_N b\neq0 }$$ a contradiction. $\Box$

\corollary{Any positive weight nonzero class $a\in\H_G^*(\C)$ with $\L(1)a=(wt~a)a$
cannot be killed by $\partial$.}
\proof
If $a$ is killed by $\partial$, then $a(m)=0$ for some $m<0$. The preceding lemma says that
$a$ must be in the center of $\H_G^*(\C)$. But $\L(1)a=(wt~a)a$ and $wt~a>0$ imply that $a$ does not commute with $\L$, a contradiction. $\Box$

\subsec{A general spectral sequence in the Cartan model}

We now generalize the spectral sequence for computing $\H^*_G(\cA)$ to
non-abelian $G$ in the Cartan model. In a future paper, we will
give an example to show that unlike in the classical case, the
spectral sequence in the chiral case does not collapse at $E_1$,
in general.

Recall that the chiral Cartan model for a $O(\gs\gg)$-algebra $(\cA,d_\cA)$ is
$$
C_G(\cA)=(\cW_{hor}\otimes\cA)^{\gg_\geq}
$$
equipped with the chiral Cartan differential
$$
d_G=D(0)\otimes 1+1\otimes d_\cQ-(\gamma^{\xi_i'}\otimes\iota_{\xi_i})(0)+(c^{\xi_i'}\otimes L_{\xi_i})(0).
$$f
The key observation here is that $d_G$ can be broken up into two commuting differentials as follows. Write $D(0)=K(0)+J(0)$ as before, and put
$$
d=K(0)\otimes 1+1\otimes d_\cA,~~~\delta=J(0)\otimes 1-(\gamma^{\xi_i'}\otimes\iota_{\xi_i})(0)+(c^{\xi_i'}\otimes L_{\xi_i})(0)
$$
so that $d_G=d+\delta$. As usual, $\cA^s$ denotes the subspace of $\cA$ consisting of elements $a$ with $deg_\cA a=s$. Let $C_G(\cA)^{p,q}$ be the subspace of $C_G(\cA)$ consisting of elements with
$$
\beta\gamma\#+deg_\cA=q,~~~\beta\gamma\#-b\#=p.
$$
Note that the vertex operators $\Theta_\cS^\xi\otimes 1+1\otimes L_\xi\in\cW_{hor}\otimes\cA$ are homogeneous of degrees $(p,q)=(0,0)$, and the $b^{\xi}\otimes 1$ have degree $(-1,0)$. It follows that $C_G(\cA)$ is graded by the number $(p,q)$.

\lemma{We have that
\item{a.} $d,\delta$ are $O(\gs\gg)$-invariant;
\item{b.} $d,\delta$ preserve $C_G(\cA)$;
\item{c.} $d^2=\delta^2=[d,\delta]=0$;
\item{d.} $d:C_G(\cA)^{p,q}\ra C_G(\cA)^{p,q+1},~\delta:C_G(\cA)^{p,q}\ra C_G(\cA)^{p+1,q}$.}
\proof
a. Recall that $[J(0),b^\xi]=\Theta_\cW^\xi$. Since $J(0)^2=0$, it follows that $[J(0),\Theta_\cW^\xi]=0$. Likewise $[D(0),\Theta_\cW^\xi]=0$, hence $[K(0),\Theta_\cW^\xi]=0$. Likewise $[d_\cA,L_\xi]=0$. It follows that $[d,\Theta_\cW^\xi\otimes 1+1\otimes L_\xi]=0$. Since $d_G=d+\delta$ is $O(\gs\gg)$-invariant, so is $\delta$.

b. $K=\gamma^{\xi_i'}b^{\xi_i}$ obviously preserves $\cW_{hor}=\bra b,\beta,\gamma\ket$. Since $K(0),d_\cA$ are both $O(\gs\gg)$-invariant, it follows that they both, and hence $d$ too, preserve $C_G(\cA)$.

c. Since $d_G^2=d^2+[d,\delta]+\delta^2$, assertion c. follows from d., which we show next.

d. Recall that $J=:c^{\xi_i'}(\Theta_\cS^{\xi_i}+\half\Theta_\cE^{\xi_i}):$.
Thus we can further break up $J(0)$ into two terms, the first being $(c^{\xi_i'}\Theta_\cS^{\xi_i})(0)$. If we add to this the term $(c^{\xi_i'}\otimes L_{\xi_i})(0)$
appearing in $\delta$, the sum acting on $C_G(\cA)$ takes the form
$$
(c^{\xi_i}L_\xi^{tot})(0)=\sum_{n\geq0}c^{\xi_i'}(n)L_\xi^{tot}(-n-1)
$$
where $L_\xi^{tot}=\Theta_\cS^\xi\otimes 1+1\otimes L_\xi$, because $C_G(\cA)$ is $\gg_\geq$-invariant. (This is also consistent with the fact that $C_G(\cA)$ has no $c^{\xi'}$.) We now list the effects of all terms appearing in $d_G$ on the various gradings on $C_G(\cA)$:
\eqn\dumb{\matrix{
operators & b\# & \beta\gamma\# &deg_\cA\cr
K(0)\otimes 1=\sum_{n\geq0}\gamma^{\xi_i'}(n)b^{\xi_i}(-n-1)& +1&+1&0\cr
1\otimes d_\cA & 0 & 0 & +1\cr
-(\gamma^{\xi_i'}\otimes\iota_{\xi_i})(0)& 0 & +1 & -1\cr
(c^{\xi_i'}\otimes L_{\xi_i}^{tot})(0) &  -1 & 0 & 0\cr
\half(c^{\xi_i'}\Theta_\cE^{\xi_i})(0) & -1 & 0 & 0.
}}
The first two operators add up to $d$ and the rest add up to $\delta$. From the table,
it follows that $d,\delta$ have the right effects on $C_G(\cA)^{p,q}$ as claimed. $\Box$

It is also clear that the weight zero piece of the complex coincides with the
classical Cartan complex with differentials reduced to $d=1\otimes d_\cA$ and $\delta=-\gamma^{\xi_i'}(-1)\otimes\iota_{\xi_i}(0)$. In particular in weight zero, the double complex structure
above agrees with the classical one.

As usual, associated to the double complex structure on the chiral Cartan complex, there are two filtrations and two spectral sequences. We shall consider the following one:
$$
F^n_k=\oplus_{p+q=n,~p\geq k}~ C_G(\cA)^{p,q}.
$$
Let $(E_r,\delta_r)$ be the spectral sequence associated with this filtration. Let's specialize to the case
$$
\cA=\cQ(M).
$$

\theorem{In each weight, the spectral sequence $(E_r,\delta_r)$
converges to the graded object associated with
$H^*(C_G(\cQ(M)),d+\delta)$. In fact, in each weight, the
spectral sequence collapses at $E_r$ for some $r$.} \proof The
first statement follows immediately from the fact that both
$d,\delta$ are operators of weight zero, and the filtering spaces
$F_k$ and the terms in the spectral sequence are all graded by the
weight. In a given weight $m$, we claim that there are no nonzero
elements in $C_G^{p,q}$ for $|q-p|>>0$. Note that
$|q-p|=|deg_\cQ+b\#|$.  So if $|q-p|$ is not bounded then either
the vertex operators $\partial^k b^{\xi_i}$, or the $\partial^k
c^i$, or the $\partial^k b^i$ would have to be present on some
coordinate open set of $M$, with unbounded $k$, because these
operators are fermionic. But $wt~\partial^k
b^{\xi_i}=wt~\partial^kc^i=k$ and $wt~\partial^kb^i=k+1$,
violating that $m$ is fixed and that the weights in
$C_G(\cQ(M))$ are bounded below by zero. This shows that for a
given weight $m$, $E_r[m]=E_{r+1}[m]=\cdots$ for all $r>>0$.
$\Box$

Note that
$$
E_1^{p,q}=H^q(C_G(\cQ(M))^{p,*},d),~~~~\delta_1=\delta:E_1^{p,q}\ra E_1^{p+1,q}.
$$

\subsec{A general spectral sequence in the Weil model}

Let $(\cA,d_\cA)$ be a $O(\gs\gg)$-algebra. The Weil model for the chiral equivariant cohomology of $(\cA,d_\cA)$ is the cohomology of the complex
$$
\cD_G(\cA)=((\cW(\gg)\otimes\cA)_{bas},~K(0)\otimes 1+J(0)\otimes 1+1\otimes d_\cA).
$$
The three terms in the differential have the following gradings:
\eqn\dumb{\matrix{
operators & bc\# & \beta\gamma\# &deg_\cA\cr
K(0)\otimes 1& -1&+1&0\cr
J(0)\otimes 1&+1&0&0\cr
1\otimes d_\cA & 0 & 0 & +1.
}}
Let $\cD_G(\cA)^{p,q}$ be the subspace of $\cD_G(\cA)$ of homogeneous degree
$$
p=bc\#+\beta\gamma\#+deg_\cA,~~~q=\beta\gamma\#.
$$
The vertex operators $\Theta_\cW^\xi\otimes 1+1\otimes L_\xi$ are
homogeneous of degrees $(p,q)=(0,0)$, and the $b^\xi\otimes
1+1\otimes\iota_\xi$ have degrees $(p,q)=(-1,0)$. The operators
$$
d=K(0)\otimes 1,~~~\delta=J(0)\otimes1+1\otimes d_\cA
$$
have degrees $(p,q)=(0,1)$ and $(p,q)=(1,0)$ respectively.

\lemma{We have that
\item{a.} $d,\delta$ preserve $\cD_G(\cA)$;
\item{b.} $d^2=\delta^2=[d,\delta]=0$;
\item{c.} $d:\cD_G(\cA)^{p,q}\ra \cD_G(\cA)^{p,q+1},~\delta:\cD_G(\cA)^{p,q}\ra \cD_G(\cA)^{p+1,q}$.}
\proof
a. We have seen that $[K(0),\Theta_\cW^\xi]=0$, implying that $d$ commutes with $\Theta_\cW^\xi\otimes 1+1\otimes L_\xi$. We also have that $K(0)$ commute with $b^\xi\otimes 1+1\otimes\iota_\xi$. Thus $K(0)$ preserves $\cD_G(\cA)$. Since $d+\delta$ preserves $\cD_G(\cA)$, so does $\delta$.

b. Since $0=(d+\delta)^2=d^2+[d,\delta]+\delta^2$, assertion b. follows from c., which we have shown above.  $\Box$

As in the Cartan model, we have the filtration
$$
F^n_k=\oplus_{p+q=n,~p\geq k}~ \cD_G(\cA)^{p,q}.
$$
Let $(E_r,\delta_r)$ be the spectral sequence associated with this filtration. Let's specialize to the case
$$
\cA=\cQ(M).
$$

\theorem{In each weight, the spectral sequence $(E_r,\delta_r)$ converges to the graded object associated with $H^*(\cD_G(\cQ(M)),d+\delta)$. In fact, in each weight, the spectral sequence collapses at $E_r$ for some $r$.}
\proof
The argument is similar to the case of the Cartan model above. The only difference is that
we now have $|q-p|=|deg_\cQ+bc\#|$. So if $|q-p|$ is
not bounded then either the vertex operators $\partial^k b^{\xi_i}$, or the $\partial^kc^{\xi_i'}$, or the $\partial^k c^i$, or the $\partial^k b^i$ would have to be present on some coordinate
open set of $M$, with unbounded $k$, because these operators are
fermionic. The rest of the argument is verbatim. $\Box$

Note that
$$
E_1^{p,q}=H^q(\cD_G(\cQ(M))^{p,*},d),~~~~\delta_1=\delta:E_1^{p,q}\ra E_1^{p+1,q}.
$$

\subsec{Abelianization?}

Recall that classically if $T$ is a closed subgroup of $G$ then
$\gg^*\ra\gt^*$ induces a map $W(\gg^*)\ra W(\gt^*)$. Since every $G^*$-algebra
is canonically a $T^*$-algebra, it follows that for any given $G^*$-algebra $A$
one has a canonical map
$$
A\otimes W(\gg^*)\ra A\otimes W(\gt^*).
$$
This induces on cohomology a map $H_G^*(A)\ra H_T^*(A)$. In fact, when $T$ is a maximal torus of $G$, then a spectral sequence argument shows that this map yields an isomorphism
$$
H_G^*(A)\cong H_T^*(A)^W
$$
where $W=N(T)/T$ is the Weyl group of $G$. See Chap. 6 \GS.

One might expect that there would be a similar construction in the vertex algebra setting. Unfortunately, this cannot be expected to go through, at least not in a naive way. Here is why.

\lemma{$\cW(\gg)$ is simple. In other words, it has no nontrivial
ideal.} \proof As before, we can regard $\cW(\gg)$ as a polynomial
(super) algebra with generators given by the negative Fourier
modes $b^{\xi}(n),c^{\xi'}(n)$ (odd) and
$\beta^{\xi}(n),\gamma^{\xi'}(n)$ (even), $n<0$, which are linear
in $\xi\in\gg$ and $\xi'\in\gg^*$. In this polynomial
representation, each of the non-negative Fourier modes act by
formal differentiation; for example
$b^{\xi}(m)c^{\xi'}(-n)=\bra\xi',\xi\ket\delta_{m-n+1,0}$ for
$m\geq0, n>0$. From this, it is clear that any nonzero polynomial
in $\cW(\gg)$ can be reduced to a nonzero scalar by a suitable
repeated application of these derivations. Translated into vertex
algebra operations, it says that any nonzero vertex operator in
$\cW(\gg)$ can be reduced to a nonzero multiple of 1 by taking
suitable repeated circle products with the generators of
$\cW(\gg)$. In other words, any nonzero ideal of the circle
algebra $\cW(\gg)$ contains 1.  $\Box$

{\it Warning.} It is not true that a positive Fourier mode of a generator of $\cW$ acts as a derivation of the circle products in $\cW$. For example $\beta^\xi(1):\cW\ra\cW$ is not a derivation of circle products. But if we represent $\cW$ as a polynomial space, then this is a derivation with respect to the usual polynomial products.
This follows immediately from the construction of $\cE,\cS$ as induced modules over a Lie algebra.
See Example \ExampleIV.



The preceding lemma shows that there exists a vertex algebra homomorphism $\cW(\gg)\ra\cW(\gt)$ extending the classical map $W(\gg^*)\ra W(\gt^*)$ only if $\gt=\gg$. Next, suppose $T\subset G$ is a maximal torus. {\it Is there a vertex algebra homomorphism $\H_G^*(\C)\ra\H_T^*(\C)$ that extends the classical map $H_G^*(\C)\ra H_T^*(\C)$?}

\theorem{The answer is negative.} 
\proof Suppose there were such a
map $f:\H_G^*(\C)\ra\H_T^*(\C)$. Since $\H_T^*(\C)$ is an abelian
vertex algebra by Theorem \AbelianCase, it follows that $f$ must
kill $\L\in\H^*_G(\C)$, hence the ideal generated by $\L$. In
particular by Theorem \CohomologyVirasoro, the image of $\partial$
on $\H_G^*(\C)[0]=H_G^*(\C)=S(\gg^*)^\gg\subset\bra\gamma\ket$ must also be killed. On the other hand, since
$\partial$ is a vertex algebra operation, it follows that
\eqn\dumb{
0=f(\partial a)=\partial f(a).
}
By assumption, the restriction $f:H_G^*(\C)\ra H_T^*(\C)$ is the classical map. Consider for example
$a=\gamma^{\xi_i'}\gamma^{\xi_i'}\in H_G^*(\C)$ which is obviously
nonzero, so that $f(a)$ is nonzero. In fact $f(a)$ has the same
shape as $a$ but we sum over only an orthonormal basis of $\gt$.
But according to the description of $\H_T^*(\C)$ given by Theorem
\AbelianCase, we have that $\partial f(a)\neq0$, contradicting
\dumb. $\Box$

This shows that one cannot hope to get new information via an abelianization that extends the classical
isomorphism $H_G^*(\C)\cong H_T^*(\C)^W$, at least in the case of a point. This suggests that in the new theory, the chiral equivariant cohomology for non-abelian groups may be far more interesting that in the classical case.


\newsec{Concluding Remarks}

We have constructed a cohomology theory $\H^*_G(\cA)$ for $O(\gs\gg)$-algebras $\cA$, which is the vertex algebra analogue of the classical equivariant cohomology of $G^*$-algebras. A principal example we give is when $\cA=\cQ(M)$, the chiral de Rham complex of a $G$-manifold $M$. It turns out that there are other similar differential vertex algebras associated to a $G$-manifold which give rise to interesting chiral equivariant cohomology. For example, we can consider the vertex subalgebra $\cQ'(M)$ generated by the weight zero subspace $\Omega(M)$ of $\cQ(M)$. It turns out that $\cQ'(M)$ is an abelian differential vertex algebra that belongs to an appropriate category on which the functor $\H^*_G$ is defined. Moreover, $\H^*_G(\cQ'(M))$ is also a degree-weight graded vertex algebra containing the classical equivariant cohomology $H^*_G(M)$ as the weight zero subspace. We can prove the following

\theorem{If the $G$-action on $M$ has a fixed point, then the chiral Chern-Weil map $\kappa_G:\H^*_G(\C)\ra\H^*_G(\cQ'(M))$ is injective. If, furthermore, $G$ is simple, then $\H^*_G(\cQ'(M))$ is a conformal vertex algebra with the Virasoro element $\kappa_G(\L)$.}

Details of this and other related results will appear in a forthcoming paper.

\footatend\vfill\supereject\immediate\closeout\rfile\writestoppt
\baselineskip=14pt\centerline{{\bf References}}\bigskip{\frenchspacing%
\parindent=20pt\escapechar=` \input refs.tmp\vfill\eject}\nonfrenchspacing

\bs
\item{} Bong H. Lian, Department of Mathematics, National University of Singapore, 2 Science Drive 2, Singapore 117543. On leave of absence from Department of Mathematics, Brandeis University, Waltham MA 02454.
\bs
\item{} Andrew R. Linshaw, Department of Mathematics, Brandeis University, Waltham MA 02454.

\end